\documentclass[times]{elsarticle}

\usepackage[utf8]{inputenc}
\usepackage{etex}
\usepackage{secteqns}
\usepackage{float}
\usepackage[english]{babel}
\usepackage[dvips]{epsfig}
\usepackage{rotating}
\usepackage{enumerate}
\usepackage{setspace}
\usepackage[section]{placeins}
\usepackage[normalem]{ulem}
\usepackage{pictex}
\usepackage{latexsym}
\usepackage{amsfonts,amsmath,amssymb,amsthm,mathtools}
\usepackage{rawfonts}
\usepackage{multirow}
\usepackage{color}
\usepackage{cancel}
\usepackage{ifthen}
\usepackage{etoolbox}
\usepackage{stmaryrd}
\usepackage{mathrsfs}
\usepackage{xcolor}
\usepackage{tikz,tikz-3dplot}
\usepackage{pgfplots}
\pgfplotsset{compat=newest}
\usepackage{algorithm}
\usepackage{algorithmicx}
\usepackage{algpseudocode}
\usepackage{subcaption}
\usepackage{changepage}
\usepackage[colorinlistoftodos,prependcaption,textsize=normalsize]{todonotes}
\usepackage{rotating}
\usepackage{comment}
\usepackage{hyperref}
\usepackage[T1]{fontenc}
\usepackage{graphicx,adjustbox}
\usepackage{type1cm}
\usepackage{eso-pic}
\usepackage{colortbl}
\usepackage{booktabs}
\usepackage[margin=2.5cm]{geometry}
\usepackage[separate-uncertainty = true,multi-part-units = repeat]{siunitx}

\hypersetup{urlcolor=blue, colorlinks=true}

\newcommand{\MB}[1]{\todo[author=\textbf{MATIAS},color=red,inline]{#1}}

\newcommand{\stefano}[2][]{\todo[author=\textbf{STEFANO},color=yellow,inline,#1]{#2}}

\newcommand{\jump}[1]{\ensuremath{[\![#1]\!]} }



\newcommand{\ie}[1] {\textit{i.e.}}
\newcommand{\eg}[1] {\textit{e.g.}}
\newcommand{\ea}[1] {\textit{et al.}}
\newcommand{\cf}[1] {\textit{cf.}}

\newcommand{\bs}[1]{\boldsymbol{#1}}

\usetikzlibrary{calc}
\definecolor{RYB1}{RGB}{69,117,180}
\definecolor{RYB2}{RGB}{145,191,219}
\definecolor{RYB3}{RGB}{0,153,76}
\definecolor{RYB4}{RGB}{173,255,47}
\definecolor{RYB5}{RGB}{252,141,89}
\definecolor{RYB6}{RGB}{215,48,39}

\newtheorem{theorem}{Theorem}

\newtheorem*{remark}{Remark}

\newcommand{\mockalph}[1]{}

\begin{document}

\begin{frontmatter}


\title{An arbitrary order Mixed Virtual Element formulation for coupled multi-dimensional flow problems}

\author[UQ]{Benedetto, M. F.*}
\author[POLI]{Borio, A.}
\author[ITBA]{Kyburg F.}
\author[ITBA]{Mollica, J.}
\author[POLI]{Scial\`o, S.}

\address[UQ]{School of Mathematics and Physics, The University of Queensland, Brisbane, QLD, Australia}
\address[ITBA]{Computational Mechanics Center, Mechanical Engineering Department, Instituto Tecnológico de Buenos Aires (ITBA), Argentina}
\address[POLI]{Department of Mathematical Sciences “G.L. Lagrange”, Politecnico di Torino, Italia}

\begin{abstract}
Discrete Fracture and Matrix (DFM) models describe fractured porous media as complex sets of 2D planar polygons embedded in a 3D matrix representing the surrounding porous medium. The numerical simulation of the flow in a DFM requires the discretization of partial differential equations on the three dimensional matrix, the planar fractures and the one dimensional fracture intersections, and suitable coupling conditions between entities of different dimensionality need to be added at the various interfaces to close the problem. The present work proposes an arbitrary order implementation of the Virtual Element method in mixed formulation for such multidimensional problems. Details on effective strategies for mesh generation are discussed and implementation aspects are addressed. Several numerical results in various contexts are provided, which showcase the applicability of the method to flow simulations in complex multidimensional domains.
\end{abstract}

\begin{keyword}
Mixed VEM, DFN, Subsurface flow, Inter-dimensional coupling
\end{keyword}

\end{frontmatter}

\section{Introduction}
\label{Intro} 
There are many practical contexts where effective flow simulations in underground fractured media are
strategic, including  geothermal applications, protection of water
resources, Oil\&Gas enhanced production and geological waste storage. Taking advantage from an increased and easily available computational power, several
problems not considered in the past have been tackled. Regardless of the application, they all share the demand for high accuracy and reliability in the results. On the other hand, due to the complexity of the typical domains of interest and to the high uncertainty of the data, effective simulations of underground phenomena are still extremely challenging, such that the research for robust and efficient numerical methods for underground flow simulations still attracts great interest from the scientific community.

In this work, the computation of the hydraulic head distribution in the subsoil is considered. The physical components of the problem are a rock matrix with an embedded network of fractures. Fractures are thin regions of the soil with different properties from the surrounding bulk material, and have one spatial dimension, the thickness, that can be orders of magnitude smaller than the domain size. For numerical simulations, the simultaneous representation of the fracture-thickness scale and of the domain-scale is unfeasible as would result in an extremely large number of unknowns, such that models are introduced to represent subsoil. Next to homogenization techniques \cite{Qi2005}, dual and multy-porosity models \cite{DPbook}, or embedded discrete fracture matrix (EDFM) models \cite{Li2008,Moinfar2014}, Discrete Fracture and Matrix (DFM) models aim at an explicit representation of the underground fractures, which are dimensionally reduced to planar interfaces into the porous matrix. Fractures are generated randomly following probability distribution functions concerning their geometrical (position, orientation, density) and hydraulic properties. The quantity of interest is the hydraulic head distribution, which is governed by the Darcy law in the porous matrix, and by an averaged-across-thickness Darcy law in the fracture planes, plus additional coupling conditions at fracture/matrix interfaces an at the intersections between fractures, \cite{ModelingFractures}. Despite the dimensional-reduction operated on the fractures, DFM models are still highly complex and multi-scale: this is a consequence of the random orientation of the fractures that usually form an intricate system of intersections, with the presence of fractures with different sizes and forming intersections that might span various orders of magnitude. In fact one planar fracture in a DFM might have an intersection of few centimeters length with one fracture and of several kilometers with another fracture, with also the smaller intersection having a relevant impact on the flow pattern. This geometrical severity, combined to the stochastic nature of simulation data, demands numerical tools robust to complex geometries and highly efficient, thus allowing to perform repeated simulations on random geometries necessary to obtain statistics on the output quantities of interest.

DFM models are widely used for underground flow simulations \cite{Angot2009,Ahmed2015,Brenner2016,Antonietti2016,BERRONE2017768}, the major complexity being the generation of a conforming mesh of the domain. The generation of a conforming mesh for the imposition of the matching conditions at the various interfaces, might in fact result in an impossible task, for the extremely high number of geometrical constraints for fracture networks of practical interest. Further the mesh generation process with conventional strategies is a global process for the whole domain, which usually requires an iterative process that might not converge.
In some cases the rock matrix can be neglected, as almost impervious compared to the fractures, with the subsurface flow mostly determined by the fracture distribution. These are called Discrete Fracture Network (DFN) \cite{DF99,X2,Neuman05,NTTDA92} problems, with an only partially mitigated geometrical complexity.

Over the last decade, there has been a great development of numerical methods to tackle the problem of efficient flow simulations of realistic DFM/DFNs. An efficient algorithm for conforming discretizations have been proposed by \cite{NGO2017}. The complexity of DFN flow simulations is reduced in \citep{NOETJCP12,NOETJCP15} by removing the unknowns in the interior of the fractures, reducing the dimension of the problem and rewriting it at the interfaces. 
The mesh conformity requirement at the interfaces in DFMs can be relaxed by using eXtended Finite Elements (XFEM) \cite{XFEMreview} as in \cite{FS13,formaggia2014}. In \cite{BSV,BFPSV18,BERRONE2017768}, an optimization approach is proposed for both DFN and DFM problems which avoids any need for mesh conformity at the interfaces and instead seeks the solution as the constrained minimum of a functional representing the error in the fulfillment of interface conditions. 
In recent times, techniques as the Mimetic Finite Difference method (MFD) \cite{Lipnikov2014} have been
used for flow simulations in DFMs by \cite{Wheeler2015,Antonietti2016}, or as Hybrid High Order (HHO) methods by \cite{Chave2018}, where a partial non-conformity is allowed between the mesh of the porous medium and of the fractures. Other approaches use two or multi-point flux approximation based techniques, \cite{Sandve2012,Hajibeygi2011,Gable2015bTransp,Faille2016} for DFN and DFM problems, or gradient schemes \cite{Brenner2016}. A survey on various conforming and non-conforming discretization strategies for flow simulations in networks of fractures can be found in \cite{FKS}.

The recently developed Virtual Element Method (VEM) \cite{VEMmixedBasic,Beirao2016,VEMbasic} is gaining increasing interest in the field of the numerical simulation of underground phenomena as it allows to handle mesh of polygonal/polyhedral elements and is robust to badly shaped and elongated elements \cite{BBorth}, thus allowing to easily generate conforming polygonal/polyhedral meshes in complex geometries. This method was applied, e.g., in \cite{BBPS,DFNvem1, BBBChapter,Fumagalli2018} for DFN simulations, and in DFM problems in \cite{Fumagalli2019} for flow computation, in \cite{Berrone2018} coupled to the Boundary Element method and in \cite{Coulet2019} coupled to finite volumes for poro-elasticity problems in DFMs.

The present work proposes an arbitrary order mixed VEM-based approach for the computation of the flow in poro-fractured media, following the DFM model proposed by \cite{Nordbotten2019}. A mixed formulation is a widely common choice for underground flow simulation for its mass conservation properties \cite{ModelingFractures,MixedDarcyNonMatching,FS13,P2,VMS07,AlHinai2015,Wheeler2015,Antonietti2016}. The approach is an extension of the work proposed in \cite{DFNvem1} and \cite{VEMmixedDFN}, as now the VEM-based conforming approach includes the rock matrix in the problem domain, and differs from other VEM based approaches for the strategy used to obtain the computational mesh. Further, to the best of authors knowledge, this is the first arbitrary-order implementation of mixed virtual elements for flow simulations in DFM problems.

The manuscript is organized as follows: in Section~\ref{prob_for} the formulation for the problem at hand is presented. 
Section~\ref{MeshConform} described the mesh generation process, whereas Section~\ref{VEM} is devoted to providing a description of the mixed formulation of the Virtual Element Method in the present context. Implementation is discussed in Section~\ref{impl}. Next, numerical results are described in Section~\ref{num_res}, where convergence analysis of the method is proposed and problems on increasingly complex configurations are solved and analyzed. The work ends with some concluding remarks in Section~\ref{Conclusions}.

\section{Problem formulation} 
\label{prob_for}

\subsection{DFN problem formulation}
\label{problemformulation}
\newcommand{\ellrange}{\ell=1,\ldots,N^d}
\newcommand{\gammadd}{\gamma^{(d,d-1)}}
\newcommand{\avg}[1]{\{\!\{#1\}\!\}}
\newcommand{\xx}{\bs{x}}
\newcommand{\nn}{\bs{n}}
\newcommand{\uu}{\bs{u}}
\newcommand{\vv}{\bs{v}}

The present section is devoted to briefly recall the formulation of the flow problem in fractured porous media with the Discrete Fracture and Matrix model, referring to \cite{DFNBoon} for a more detailed description and for well posedness results. 

Let us consider a three dimensional block of porous material $\Omega^3$ crossed by a network of fractures. According to the DFM model, fractures are represented as planar polygons $\Omega^2_\ell$, $\ell=1,\ldots,N^2$, which might intersect forming intersection segments $\Omega^1_\ell$, $\ell=1,...,N^1$, also called traces. Further, traces can meet at intersection points $\Omega^0_\ell$, $\ell=1,\ldots,N^0$. For uniformity of notation, a subscript will be indicated also for $\Omega^3$, with the assumption that $\Omega^3_1\equiv \Omega^3$, and $N^3=1$. We then denote by $\Omega^2=\bigcup_{\ell=1}^{N^2} \Omega^2_\ell$ the union of all fractures, by $\Omega^1=\bigcup_{\ell=1}^{N^1} \Omega^1_\ell$ the union of all traces, and by $\Omega^0=\bigcup_{\ell=1}^{N^0} \Omega^0_\ell$ the union of all trace intersections. Problem domain is thus mixed dimensional, as it involves a 3D problem in $\Omega^3$, 2D problems on the fractures $\Omega^2$, 1D problems on the traces $\Omega^1$ and 0D problems at trace intersections $\Omega^0$. Each domain $\Omega^d$, $d=1,\ldots,3$ does not include lower dimensional domains, i.e. $\Omega^d\cap\left(\bigcup_{j=0}^{d-1}\Omega^j\right)=\emptyset$, such that \textit{internal} boundaries are present. For $d=1,\ldots,3$, the internal boundary of domain $\Omega_\ell^d\subset\mathbb{R}^d$, $\ellrange$, is the portion of boundary that matches a sub-dimensional domain, the remaining being instead the external boundary. 

For each domain $\Omega_\ell^d$, $d=0,\ldots,2$, $\ell=1,\ldots,N^d$, let us introduce the index set $\mathcal{O}_\ell^{d^+}$, containing indexes $j=1,\ldots,N^{d+1}$ such that,if $j\in\mathcal{O}_\ell^{d^+}$, $\Omega^{d}_\ell$ coincides with a portion of the boundary of $\Omega^{d+1}_j$, i.e. $\bar{\Omega}^{d+1}_j\cap \Omega^{d}_\ell\neq \emptyset$. It is further set $\mathcal{O}^{3^+}=\emptyset$. Similarly, for $d=1,\ldots,3$, $\mathcal{O}_\ell^{d^-}$ contains indexes of domains $\Omega_j^{d-1}$, such that $\Omega_j^{d-1}$, for $j\in\mathcal{O}_\ell^{d^-}$ has a non empty intersection with the boundary of $\Omega^d_\ell$. 

\begin{figure}
\includegraphics[width=0.99\textwidth]{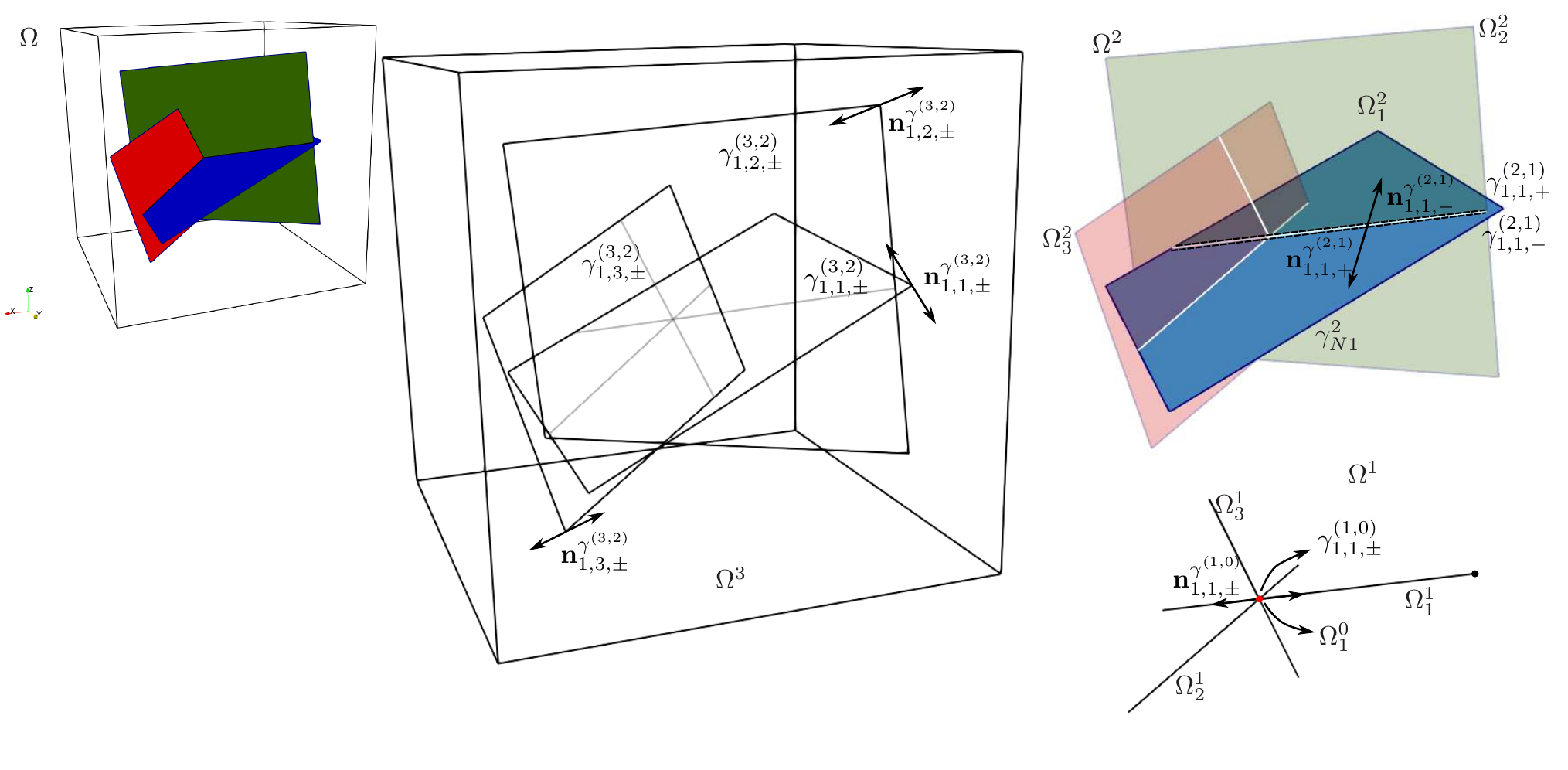}
\caption{Exemplification of the nomenclature for a sample domain}
\label{fig:nomenclature}
\end{figure}

For $d=1,\ldots,3$, we denote by $\gamma^d_{D,\ell}$ and $\gamma^d_{N,\ell}$ the Dirichlet and Neumann part, respectively, of the external boundary of $\Omega^d_\ell\subset\mathbb{R}^d$, $\ell=1,\ldots,N^d$, $\gamma^d_{D\ell}\cap\gamma^d_{N\ell}=\emptyset$, and by $\gamma^{(d,d-1)}_{\ell,j,+}$, $\gamma^{(d,d-1)}_{\ell,j,-}$ the boundary of $\Omega^d_\ell$ around the lower dimensional domain $\Omega^{d-1}_j$, $j\in\mathcal{O}_\ell^{d-}$, with fixed arbitrarily chosen $\pm$ sign for each of the two sides. For $d=2,3$, $\nn_{\ell,j,\pm}^{\gamma^{(d,d-1)}}$ is used to denote the unit normal vector to $\gamma^{(d,d-1)}_{\ell,j,\pm}$ in $\Omega_\ell^d$, outward pointing, and, for $d=1$, $\nn_{\ell,j,\pm}^{\gamma^{(d,d-1)}}$ is the unit vector tangential to $\Omega_\ell^d$ at the boundary points $\gamma^{(d,d-1)}_{\ell,j,\pm}$, outward pointing. An exemplification of the used nomenclature is proposed in Figure~\ref{fig:nomenclature}, for a simple DFN counting three fractures, three traces and one trace intersection.

A Darcy-type equation governs the flow problem on each geometrical dimension and suitable matching conditions couple the problems: for $\ell=1,\ldots,N^d$ we have
\begin{eqnarray*}
&\uu^d_\ell(\xx)=\mathfrak{a}_\ell^d(\xx) \nabla p^d_\ell(\xx), & d=1,\ldots,3 \\
&\nabla \cdot \uu^d_\ell(\xx) - \sum_{j\in \mathcal{O}_\ell^{d^+}}\jump{\uu^{d+1}_j(\xx)\cdot \nn_{j,\ell,\pm}^{\gamma^{(d+1,d)}}} = f^d_\ell(\xx), & d=1,\ldots,3 \\
&-\sum_{j\in \mathcal{O}_\ell^{d^+}}\jump{\uu^{d+1}_j(\xx)\cdot \nn_{j,\ell,\pm}^{\gamma^{(d+1,d)}}}=f^d_\ell(\xx), & d=0,
\end{eqnarray*}
with boundary conditions, for $d=1,\ldots,3$, $\ellrange$,
\begin{eqnarray*}
&\uu^d_\ell\cdot \nn_\ell^d=0, & \text{on } \gamma_{N\ell}^d,\\
& p^d_\ell=0, & \text{on } \gamma_{D\ell}^d,
\end{eqnarray*}
and coupling conditions for $d=0,\ldots,3$
\begin{eqnarray*}
&\uu^{d+1}_j(\xx)\cdot \nn_{j,\ell,\pm}^{\gamma^{(d+1,d)}}=-\eta_\ell^d \left(p^d_\ell(\xx)-\left(p^{d+1}_j(\xx)\right)_{|\gamma^{(d+1,d)}_{j,\ell,\pm}}\right), & j \in \mathcal{O}_\ell^{d^+}. \\
\end{eqnarray*}
In the previous equations, for $d=0,\ldots,3$, $\ellrange$, $f_\ell^d$ represents a source term, the operator $\jump{\cdot}$ denotes the jump across the interface, i.e.
\begin{displaymath}
\jump{\uu^{d+1}_j(\xx)\cdot \nn_{j,\ell,\pm}^{\gamma^{(d+1,d)}}}=\uu^{d+1}_j(\xx)\cdot \nn_{j,\ell,+}^{\gamma^{(d+1,d)}}+\uu^{d+1}_j(\xx)\cdot \nn_{j,\ell,-}^{\gamma^{(d+1,d)}}, \quad j \in \mathcal{O}_\ell^{d^+},
\end{displaymath}
whereas, for $d\neq 0$, $\mathfrak{a}_\ell^d(\xx)$ is the fracture transmissivity in $\Omega_\ell^d$ and $\eta_d^\ell$ is the transmissivity in the direction normal to $\Omega_\ell^d$, with clear extension to the case $d=0$. Homogeneous Neumann and Dirichlet boundary conditions are used in order to simplify the exposition. The choice of homogeneous Neumann boundary conditions for external boundaries of domains $\Omega_\ell^d$ $d=1,2$, $\ellrange$, not touching the boundary of $\Omega^3$ is widely adopted, see e.g. \cite{ModelingFractures,Antonietti2016,Nordbotten2019}, other generalizations being straightforward.

Let us now move to the weak formulation of the previous problem. Let us introduce, for $d=1,\ldots,3$, $\ellrange$, the functional spaces $V^d_\ell=\mathrm{H}_0(\nabla\cdot,\Omega^d_\ell):=\left\lbrace \vv\in \mathrm{H}(\nabla\cdot,\Omega^d_{\ell}) : \vv\cdot\nn^d_{\ell} = 0\right\rbrace$, and for $d=0,\ldots,3$ the spaces $Q^d_\ell=\mathrm{L}^2(\Omega^d_\ell)$. It is then possible to write the problem in mixed weak formulation as: for $d=1,\ldots,3$, $\ellrange$, find $\uu^d_\ell\in V_\ell^d$, and for $d=0,\ldots,3$, $\ellrange$ find $p_\ell^d\in Q_\ell^d$ such that, for all $\vv^d_\ell\in V_\ell^d$, $q_\ell^d\in Q_\ell^d$:
\begin{eqnarray}
  \left(\left(\mathfrak{a}_\ell^d\right)^{-1}\uu^d_\ell,\vv^d_\ell
  \right)_{\Omega_\ell^d} -\left(p_\ell^d,\nabla\cdot \vv_\ell^d
  \right)_{\Omega_\ell^d} - \sum_{j\in\mathcal{O}^{d^-}} \sum_{\xi=+,-}
  \left(\eta_j^{d-1}\right)^{-1} \left(\uu^{d}_\ell\cdot \nn^{\gammadd}_{\ell,j,\xi},
  \vv_\ell^d \cdot \nn^{\gammadd}_{\ell,j,\xi}\right)_{\gammadd_{\ell,j,\xi}}
  +&& \nonumber
  \\
  \sum_{j\in\mathcal{O}^{d^-}} \left(p^{d-1}_j,\jump{\vv^{d}_\ell\cdot
  \nn_{\ell,j,\pm}^{\gamma^{(d,d-1)}}}\right)_{\gammadd_{\ell,j,\pm}}=0,
   && \text{if } d\neq 0 \label{weak:1}
  \\
  \left(\nabla \cdot \uu_\ell^d,q_\ell^d\right)_{\Omega_\ell^d} -
  \sum_{j\in \mathcal{O}^{d^+}}\left(\jump{\uu^{d+1}_j
  \cdot \nn_{j,\ell,\pm}^{\gamma^{(d+1,d)}}},q_\ell^d\right)_{\Omega_\ell^d}
  = \left(f^d_\ell,q_\ell^d\right)_{\Omega_\ell^d},
   && \text{if } d\neq 0 \label{weak:2}
  \\
  -\sum_{j\in \mathcal{O}^{d^+}}\left(\jump{\uu^{d+1}_j\cdot
  \nn_{j,\ell,\pm}^{\gamma^{(d+1,d)}}},q_\ell^d\right)_{\Omega_\ell^d}
  = \left(f^d_\ell,q_\ell^d\right)_{\Omega_\ell^d}
   && \text{if } d=0. \label{weak:3}
\end{eqnarray}
The used coupling equations state that the drop of pressure is proportional to $\eta^{d-1}$, which is a sort of Darcy's law for interdimensional flux exchange. In the limit when $\eta^{d-1} \rightarrow \infty$, and for geometrical parameters as fracture width, trace diameter that are very small compared to the other dimensions, these coupling conditions can be modeled as requiring that the pressure be continuous across fractures, traces and trace intersections; further, the third term in \eqref{weak:1} vanishes. This means, for instance, that there would be no jump in pressure in a path that leaves the matrix, goes through a fracture, enters the trace network and eventually re-enters the matrix. In this case, the discrete solution should converge to a pressure field that is globally continuous across all domains. 


\section{Mesh generation}
\label{MeshConform}
A key aspect of a numerical tool for the simulation of flow problems in mixed-dimensional domains with complex geometries is the generation of the computational mesh, as usually this is a non trivial and computationally expensive procedure. The meshing for standard discretization methods, as Finite Volumes or Finite Elements, is, indeed an iterative process, aiming at building good quality elements from the three dimensional down to the one dimensional domains which are perfectly matching, i.e. conforming, at the various interfaces. For non trivial geometries this process is likely to fail or to produce an extremely large number of elements, independently of the required level of accuracy, due to the need to honor the geometrical constraints.

\begin{figure}
\centering
\includegraphics[width=0.3\textwidth]{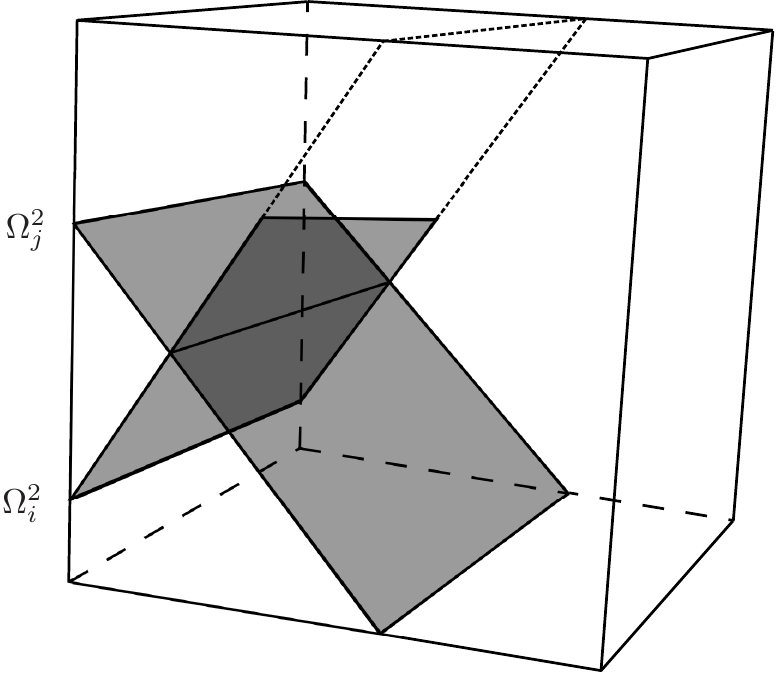}
\caption{Example of cutting of an element $E$ crossed by two fractures for mesh generation}
\label{meshing}
\end{figure}

The use of polyhedral/polygonal meshes, in conjunction with the robustness of virtual elements to badly shaped elements \cite{BERRONE201714}, allows instead to define an easy meshing process, which can be performed in a non iterative manner. The starting point is a general polyhedral mesh of the closure of the three dimensional domain $\Omega^3$, built independently of any lower dimensional domain contained in it, and thus not conforming to the interfaces. Let us then consider a generic element $E$ of such mesh, crossed by some two dimensional domains, possibly intersecting in it or ending in its interior. Figure~\ref{meshing} shows an example where a cubic element $E$ is crossed by two fractures, $\Omega^2_i$ and $\Omega^2_j$. The element can be easily cut into sub-elements not crossing the fractures, eventually prolonging the cut, for fractures ending in the interior of the element, up to element boundary. Please observe that the geometry of the fractures is not altered in any way, as the original fracture boundaries are preserved and, when the cut is prolonged over fracture boundaries, co-planar ``hanging'' faces are introduced. Also ``hanging'' faces appear on elements neighboring cut elements. Once this process is completed, the mesh on the 2D domains is obtained simply collecting the faces of the 3D elements laying on each fracture plane as a ``patchwork'', and similarly for the 1D domains. The resulting mesh is thus conforming.

\section{Mixed Virtual Elements}
\label{VEM}
In the following, an outline of the main definitions and features of
the mixed formulation of VEM is presented. Its initial introduction is
given in \cite{VEMmixedBasic}, with a follow-up work generalizing the
method \cite{Beirao2016} and a related work on virtual
$\mathrm{H}(\nabla\cdot)$ and $H(\nabla\times)$ spaces in
\cite{divcurlVEM}. A thorough description of mixed VEM spaces for the
Stokes, Darcy and Navier–Stokes equations is given in
\cite{DASSI2019}. Despite its recent introduction, a variety of
applications can be found in the scientific literature. Namely, Stokes
flow in \cite{caceresmixed1} and \cite{caceresmixed2}, the Brinkman
problem \cite{VEMbrinkman}, plane elasticity \cite{VEMmixedStress} and
flow in networks of fractures in an impervious matrix
\cite{VEMmixedDFN}.

Let us consider a domain $\Omega_\ell^d$, for $d=2,3$ and
$\ell=1,\ldots,N^d$, whose mesh, comprised of arbitrary polyhedra with
mesh parameter $h$, is indicated by $\mathcal{T}_{h,\ell}^d$, and
satisfies basic regularity conditions as in \cite{VEMmixedBasic}. In
the following we will drop the domain index $\ell$, for simplicity of
notation, as this plays no role in the discussion.


Let us then introduce a local VEM space for the velocity variable on
an element $E \in \mathcal{T}_h^d$, thus a polyhedron for $d=3$ or a
polygon for $d=2$, and, to this end let us denote by $\mathbb{P}_k^d$
the space of polynomials of maximum order $k\geq 0$, in
$\mathbb{R}^d$, with the additional conventional notation
$\mathbb{P}_{-1}=0$, for $k=-1$. The local VEM space on element
$E\in \mathcal{T}_h^d$ is denoted as:
\begin{equation}
  \label{VEMspace}
  \mathrm{V}_{k,k_\nabla}^{E,d} =\left \lbrace \bs{v}_h \in \mathrm{H}
    (\nabla\cdot,E): (\bs{v}_h \cdot \bs{n})_{|_{f}} \in \mathbb{P}^{d-1}_{k}(f) \  \forall f \subset  \partial E, \nabla\cdot\bs{v}_h  \in \mathbb{P}^d_{k_\nabla}(E), \text{ and rot}(\bs{v}_h)\in \mathbb{P}^d_{k}(E)\right \rbrace,
\end{equation}
where $f\subset\partial E$ denotes a face for $d=3$ or an edge for
$d=2$. Depending on the choice of $k_\nabla$, this space might
represent an extension of BDM elements to general elements, for
$k_\nabla=k-1$, $k\geq 1$, termed BDM$k$-VEM, or an extension of
Raviart-Thomas elements, for $k_\nabla=k$, $k\geq 0$, labeled
RT$k$-VEM. Other choices are also possible, even if not considered
here \cite{divcurlVEM}.  The local space for the pressure variables on
element $E\in \mathcal{T}_h^d$, $d=2,3$, is
$\mathrm{Q}_{k_\nabla}(E) :=\mathbb{P}_{k_\nabla}(E)$.  We remark that
BDM elements will be used only for the discretization of the 3D
equations, to have the same polynomial accuracy $k$ for the $d$
pressure and the $d+1$ face-normal fluxes that are coupled by
equation~\eqref{weak:2}, for $d=1,2$.

The discrete global space on $\Omega^d$, $d=2,3$, is
\begin{align*}
  \mathrm{V}_{k,k_\nabla}^d(\Omega^d) := \left\lbrace \bs{v}_h \in \mathrm{H}
  (\nabla\cdot,\Omega^d): \bs{v}_{h|E} \in \mathrm{V}_{k,k_\nabla}^{E,d} \  \forall E \in \mathcal{T}_h^d \right\rbrace,
\end{align*}
resulting in a $\mathrm{H}(\nabla\cdot,\Omega^d)$ conforming space. The
global space for the pressure variable is
\begin{equation*}
  \mathrm{Q}_{k_\nabla}^d(\Omega^d) :=
  \left\lbrace v_h \in \mathrm{L}^2(\Omega^d) \colon
    v_{h|_E} \in \mathbb{P}_{k_\nabla}^d(E) \,
    \forall E \in \mathcal{T}_h^d  \right\rbrace
  \quad \text{for $d\in\{1,\ldots,3\}$}\,,
\end{equation*}
while for $d=0$ the degrees of freedom of $p$ are the values at domain
$\Omega^0_\ell$, $\forall \ell\in\{1,\ldots,N^0\}$. The dimension of
the space $\mathbb{P}_k^d(E)$ is $n_k^d = {{k+d}\choose{d}}$, and a
basis for this space can be chosen as the monomial base
$\mathfrak{M}_k^d(E)$:
\begin{displaymath}
  \mathfrak{M}_k^d(E)=\left\lbrace \frac{\left(\mathbf{x}-\mathbf{x}_E\right)^\alpha}{h_E^{|\alpha|}}, \ \forall \alpha=\left(\alpha_1,\alpha_2,\ldots,\alpha_d\right), \ 0\leq|\alpha|\leq k \right\rbrace
\end{displaymath}
where $\mathbf{x}_E\in \mathbb{R}^d$ is the centroid of element $E$
and $h_E$ its diameter.

We then define the space \begin{equation} \nabla \mathbb{P}_{k+1}^d(E)
  := \left\lbrace \bs{g} \in \left[ \mathbb{P}_k^d(E) \right]^d \text{
      such that } \bs{g} = \nabla \hat{m} \text{ for some } \hat{m}
    \in \mathbb{P}_{k+1}^d(E) \right\rbrace,
\end{equation}
with dimension $n_{k,\nabla}^d := n_{k+1}^d - 1$, and by
$\left(\nabla\mathbb{P}_{k+1}^d(E) \right)^\oplus$ the $\mathrm{L}^2$
orthogonal complement of $\nabla\mathbb{P}_{k+1}^d(E)$ in
$\left[\mathbb{P}_k^d(E) \right]^d$ so that
$\left[ \mathbb{P}_k^d(E) \right]^d =
\left(\nabla\mathbb{P}_{k+1}^d(E) \right) \oplus \left(
  \nabla\mathbb{P}_{k+1}^d(E) \right)^\oplus$, whose dimension is
$n_{k,\oplus}^d := d n_{k}^d - n_{k,\nabla}^d$.


The DOFs of a function $\bs{v}_h$ in $V_{k,k_\nabla}^{E,d}$, following
\cite{Beirao2016} are:
\begin{equation}
  \begin{aligned}
    \label{DOF}
    i.) \   & \frac{1} {|f|} \int_{f}  (\bs{v}_h \cdot \bs{n}_f) m \text{ dV } \quad & \forall m \in \mathbb{P}_{k}^{d-1}(f), \quad & \forall f \in \partial E \\
    ii.) \  & \frac{1} {|E|} \int_{E}  \bs{v}_h \cdot \bs{g}  \text{ dV } \quad & \forall \bs{g} \in (\nabla \mathbb{P}_{k_\nabla}^d(E)), \\
    iii.) \ & \frac{1} {|E|} \int_{E} \bs{v}_h \cdot \bs{g} \text{ dV
    } \quad & \forall \bs{g} \in (\nabla\mathbb{P}_{k+1}^d(E))^\oplus,
  \end{aligned}
\end{equation}
A proof of unisolvence can be seen in \cite{VEMmixedBasic,Beirao2016}
for BDM- and RT-VEM respectively.
The first set of DOF can be replaced by any other way to fix a
polynomial of degree $k$ on a face. A listing of the dimensions of
some of the polynomial spaces involved in the definition of the DOF is
provided in the \textit{Supplementary Material} to the present
manuscript, along with a graphical exemplification of the DOFs for a
convex hexagon and convex
polyhedron. 


For the pressure space any set of DOFs that univocally determines a
polynomial of order $k$ in dimension $d$ could be adopted, as, for
example, $n_{k}^d$ distinct point values. However, for computations,
it is advantageous to take as DOFs the $n_{k}^d$ moments with respect
to the monomial basis of order $k$, since element geometry can be
arbitrary.

Let us now introduce the projection operator
$\boldsymbol\Pi^0_k : \ V_{k,k_\nabla}^{E,d} \to
\left[\mathbb{P}_{k}^d(E) \right]^d$ as:
\begin{align}
  &\int_{E} \boldsymbol\Pi^0_k \bs{v}_h\cdot \bs{g} \text{ dV } = \int_{E}  \bs{v}_h\cdot \bs{g} \text{ dV } \qquad \forall \bs{g} \in \left[\mathbb{P}^d_{k}(E) \right]^d,
    \label{projector}
\end{align}
It can be observed that knowledge of the functions at the DOFs is
enough to compute the projector. The left hand side of
\eqref{projector} is an integral between polynomials in dimension $d$,
and can be explicitly computed by a suitable quadrature rule.
For the left hand side, since
$\left[\mathbb{P}_k^d(E) \right]^d = \left(\nabla\mathbb{P}_{k+1}^d(E)
\right) \oplus \left( \nabla\mathbb{P}_{k+1}^d(E) \right)^\oplus$, we
can find
$\tilde{\bs{g}} \in \left( \nabla\mathbb{P}^d_{k+1}(E) \right)$ and
$\bs{g}^\oplus \in \left( \nabla\mathbb{P}^d_{k+1}(E) \right)^\oplus$
such that $\bs{g} = \tilde{\bs{g}} + \bs{g}^\oplus$.
Thus, 
\begin{align}
  \label{proy1}
  \int_{E}  \bs{v}_h\cdot \bs{g} \text{ dV } = \int_{E}  \bs{v}_h\cdot \tilde{\bs{g}} \text{ dV }  + \int_{E}  \bs{v}_h\cdot \bs{g}^\oplus \text{ dV }.
\end{align}
The second term on the right hand side of this equation can be
obtained directly from the set of DOFs of type $iii$, and, for the
other term, we have that there is a polynomial
$\hat{m} \in \mathbb{P}_{k+1}(E)$ such that
$\nabla \hat{m} = \tilde{\bs{g}}$ so that applying integration by
parts we obtain
\begin{align}
  \label{proy2}
  \int_{E}  \bs{v}_h\cdot \tilde{\bs{g}} \text{ dV } = \int_{E}  \bs{v}_h\cdot \nabla \hat{m}\text{ dV } = - \int_{E} 
  \left( \nabla\cdot\bs{v}_h \right)\hat{m} \text{ dV }  + \sum_{f\subset\partial E}\int_{f} \left(\bs{v}_h \cdot \bs{n}\right)_{|_{f}} \hat{m} \text{ dS}.
\end{align}
Once again, the second term on the right hand side can be computed
directly as an integration on the faces/edges of the 3D/2D element, by
using the DOFs of type $i$. For the first term, this can be computed
once $\nabla\cdot\bs{v}_h \in \mathbb{P}_{k_\nabla}^d$ is defined, as
follows:
\begin{align}
  \label{proy3}
  \int_{E} (\nabla\cdot\bs{v}_h)q \text{ dV } = - \int_{E} \bs{v}_h \cdot \nabla q \text{ dV }  + \sum_{f\subset\partial E}\int_{f} \left(\bs{v}_h \cdot \bs{n}\right)_{|_{f}} q \text{ dS } \qquad \forall q \in \mathbb{P}_{k_\nabla}^d(E).
\end{align}
using the set of DOFs of type $i$ and $ii$.

Let us now introduce a discrete counterpart for the bi-linear form
$a^E:=\left(\nu^d\uu^d,\vv^d\right)_E$ in equation~\eqref{weak:1},
restricted on an element $E\subset\mathbb{R}^d$, for $d=2,3$, with
$\nu^d=(\mathfrak{a}^d)^{-1}$, which, for
$\bs{u}_h, \bs{v}_h \in V_{k,k_\nabla}^{E,d}$ reads as:
\begin{align}
  \label{bilform}
  a_h^E (\bs{u}_h,\bs{v}_h) := (\nu^d \boldsymbol\Pi^0_k \bs{u}_h,\boldsymbol\Pi^0_k \bs{u}_h)_E + S^E(\bs{u}_h-\boldsymbol\Pi^0_k \bs{u}_h,\bs{v}_h-\boldsymbol\Pi^0_k \bs{v}_h).
\end{align}
where $S^E$ stands for any symmetric and definite positive bilinear
form that scales like $a^E$ on the kernel of
$\boldsymbol\Pi^0_k$. Specifically, there exist two positive constants
$\alpha_*$ and $\alpha^*$ independent of the mesh and data such that
\begin{equation}
  \alpha_* a^E(\bs{u}_h,\bs{u}_h) \leq S^E(\bs{u}_h,\bs{u}_h) \leq \alpha^* a^E(\bs{u}_h,\bs{u}_h) \quad \quad \forall E \in \mathcal{T}_h^d \quad \forall \bs{u}_h \in \text{Ker} (\boldsymbol\Pi^0_k). 
  \label{Shyp}
\end{equation}
$S^E$ is usually taken as the standard Euclidean product of the vector
of values at the DOFs scaled by the measure of $E$ and an
approximation of $\nu^d$ at the barycentre or its average on the
element. More precisely,
\begin{equation}
	S^E(\bs{u},\bs{v}) = \overline{\nu^d} |E| \sum\limits_{i=1}^{n_{dof}^E}
  \text{DOF}_i(\bs{u}) \text{DOF}_i(\bs{v}) \quad \quad
  \bs{u},\bs{v} \in V_{k,h}^E
	\label{Sdef}
\end{equation}
where $\text{DOF}_i$ stands for evaluating at the i-th DOF on the
element, i.e. the linear operator
$\text{DOF}_i: \ V_{k,k_\nabla}^{E,d} \to \mathbb{R}$, defined as
$\text{DOF}_i(\cdot):=$ evaluating $(\cdot)$ at the $i$-th DOF.

Finally, the global discrete bilinear form is defined as
\begin{equation}
  \label{globalbilform}
  a_h (\bs{u}_h,\bs{v}_h) := \sum_{E \in \mathcal{T}_h^d} a_h^E (\bs{u}_h,\bs{v}_h).
\end{equation}

The discrete version of the weak form of the problem specified in
Eqs. \eqref{weak:1}-\eqref{weak:3} is then obtained replacing the
continuous variables with their discrete counterparts, using the
projection of the VEM shape functions for \eqref{bilform}. No
projection is required for terms involving the discrete pressure
variable.  With classical assumptions on the data, well posedness of
this problem is due to continuity and coercivity of $a(\cdot,\cdot)$
and the assumptions on $S$ as well as the satisfaction of an inf-sup
condition \cite{VEMmixedBasic,Beirao2016,DFNBoon}.
\section{Implementation}
\label{impl}
The details for an implementation of the three dimensional mixed
formulation of the VEM method is provided. The main reference for this
section is \cite{VEMimplementation}, from where much of the procedures
provided here are inspired. Further insight for more general problems
is given in \cite{DASSI2019}.

\subsection{Computation of the local stiffness matrix}
\label{Computation}

In this section we address the computation of the local stiffness
matrix on a general 2D or 3D polytope $E \in \mathcal{T}_h^d$, where
Virtual Element spaces are used.

In the following, it is assumed that numerical computations of
integrals of known functions over 2D and 3D polytopes can be
performed. The most straightforward approach is to divide a polygon
into triangles, or a polyhedra into tetrahedrals and use standard
Gaussian integration. Alternatives are to consider cubature
\cite{cubature}, algebraic integration by parts or heuristic methods
(like Montecarlo integration).  The procedure for computing the local
matrices needed to obtain the discrete linear system is explained
next, with emphasis on its implementation in 3D.  First we consider a
basis of $\left( \mathbb{P}_k(E) \right)^d$ denoted by
$\left\lbrace \bs{g}^E_\alpha \right\rbrace_{\alpha=1,...,d n_{k}^d}=
\left\lbrace \bs{g}_\beta^{\nabla,E}
\right\rbrace_{\beta=1,\ldots,n_{k,\nabla}^d} \cup \left\lbrace
  \bs{g}_\gamma^{\oplus,E} \right\rbrace_{\gamma=1,\ldots,
  n_{k,\oplus}^d}$, where the basis functions
$\bs{g}_\beta^{\nabla,E}$ are chosen to be the gradients of monomials
in $\mathfrak{M}_{k+1}^d(E)$, such that:
\begin{equation*}
  \bs{g}_\beta^{\nabla,E} = \nabla m_{\beta+1}^E
  \quad \forall \beta \in \left\{1,\ldots,n_{k,\nabla}^d
    = n_{k+1}^d-1\right\} ,\, m^E_{\beta+1} \in \mathfrak{M}_{k+1}(E)\,.
\end{equation*}
The local basis functions of the space $V_{k,k_\nabla}^{E,d}$ will be
denoted by
$\left\lbrace \bs{\varphi}_\alpha^E
\right\rbrace_{\alpha=1,...,n_{dof}^E}$, where $n_{dof}^E$ is the
number of DOFs for the flux variable of the element (see
\eqref{DOF}). Furthermore, we denote by
$\{\mu_\alpha^E\}_{\alpha=1,\ldots,n_{k_\nabla}}$ the basis functions
chosen for the pressure variable, that are locally a basis of
$\mathbb{P}_{k_\nabla}(E)$. These are chosen to be piecewise scaled
monomials in $\mathfrak{M}_{k_\nabla}(E)$ :
$\mu^E_\alpha = m^E_\alpha$, $\alpha \in\{1,\ldots,n_{k_\nabla}^d\}$.

\subsubsection{Local auxiliary matrices}
\label{sec:auxMat}
Firstly, several matrices will be defined, whose usefulness will
become apparent later.

Matrix $G^E \in \mathbb{R}^{d n_{k}^d \times d n_{k}^d}$ is defined
component-wise as the product of the elements in the basis of
$\left( \mathbb{P}_k(E) \right)^d$,
\begin{equation*}
  [G^E]_{\alpha\beta} = \int_E \bs{g}^E_\alpha \cdot
  \bs{g}^E_\beta \,\mathrm{d}E
  \quad \forall \alpha,\beta \in\{1,\ldots,dn_k^d\}\,,
\end{equation*}
and can be computed directly. Using the basis for
$\left( \mathbb{P}_k(E) \right)^d$, $G^E$ can be split into
\begin{equation}
  \label{MatG}
  G^E = \begin{bmatrix}
    G^{\nabla \nabla,E} & G^{\nabla \oplus,E}\\\
    G^{\oplus \nabla,E} &    G^{\oplus \oplus,E}
  \end{bmatrix} \,.
\end{equation}
In the case when $\nu^d_k \not\propto \mathcal{I}_{d\times d}$, \ie{}
when the operator is not the Laplacian, we define $G^{\nu,E}$ with the
same size as $G$ by
\begin{equation*}
  [G^{\nu,E}]_{\alpha\beta} = \int_E
  \nu^d_k \bs{g}^E_\alpha \cdot \bs{g}^E_\beta
  \,\mathrm{d}E  \quad \forall\alpha,\,\beta\in\{1,\ldots,dn_k^d\}\,.
\end{equation*}

We denote by
$H^E \in \mathbb{R}^{{n_{k_\nabla}^d}\times n_{k_\nabla}^d}$ the mass
matrix of the basis of monomials in $\mathfrak{M}_{k_\nabla}(E)$:
\begin{equation*}
  [H^E]_{\alpha\beta} = \int_E m^E_\alpha m^E_\beta\, \mathrm{d}E
  \quad \forall  \alpha,\beta\in\{1,\ldots,n_{k_\nabla}^d\}\,.
\end{equation*}
Similarly,
$H^{\#,E} \in \mathbb{R}^{n_{k,\nabla}^d \times n_{k_\nabla}^d}$ is
defined as
\begin{equation}
  \label{Hhash}
  [H^{\#,E}]_{\alpha\beta} = \int_E m^E_{\alpha+1} m^E_\beta\,\mathrm{d}E
  \quad\forall\alpha\in\{1,\ldots,n_{k,\nabla}^d\},\;
  \forall\beta \in\{1,\ldots,n_{k_\nabla}^d\} \,.
\end{equation}
Both these matrices can be computed directly. Some of the entries in
$H^E$ are repeated in $H^{\#,E}$ .

Matrix $W^E \in \mathbb{R}^{n_{k^\nabla}^d \times n_{dof}^E}$ involves
computations with ``virtual'' shape functions and is defined as
\begin{equation}
  \label{eq:MatW}
  [W^E]_{\alpha\beta} = \int_E m^E_\alpha (\nabla \cdot
  \bs{\varphi}_\beta) \, \mathrm{d}E = - \int_{E} \nabla m^E_\alpha
  \cdot \bs{\varphi}^E_\beta \,\mathrm{d}E + \int_{\partial E}
  (\bs{\varphi}^E_\beta \cdot \hat{\bs{n}}_{\partial E}) m^E_\alpha
  \text{ dS}
  \quad
  \forall \alpha \in\{1,\ldots,n_{k_\nabla}^d\},\;
  \forall\beta\in\{1,\ldots,n^E_{dof}\}\,,
\end{equation}
where integration by parts was used. We define,
$\forall \alpha\in\{1,\ldots,n_{k_{\nabla}}^d\}$ and
$\forall \beta \in\{1,\ldots,n_{dof}^E\}$
\begin{align*}
  [W_1]_{\alpha\beta} &= - \int_{E}  \nabla m^E_\alpha \cdot \bs{\varphi}^E_\beta
                        \,\mathrm{d}E \,,
  \\
  [W_2]_{\alpha\beta} &= \int_{\partial E} (\bs{\varphi}^E_\beta \cdot
                        \hat{\bs{n}}_{\partial E}) m^E_\alpha \,\mathrm{d}s \,,
\end{align*}
and $W^E = W^E_1 + W^E_2$. Since
$\nabla m^E_\alpha \in (\nabla \mathbb{P}_{k_\nabla}^E)$, $W_1$ can be
obtained immediately: recalling type $ii$ DOFs in \eqref{DOF}, we have
that, $\forall\alpha>1$,
\begin{equation*}
  [W^E_1]_{\alpha\beta} = -|E| \int_E
  \bs{g}^{\nabla,E}_{\alpha-1} \cdot \bs{\varphi}^E_\beta =-|E|\,
  \text{DOF}_{\left(n^f_E n_{k}^{d-1}+\alpha-1\right)}(\bs{\varphi}^E_\beta)
  = -|E|\delta_{\beta,(n^f_E n_{k-1}^d+\alpha-1)}
\end{equation*}
Regarding matrix $W^E_2$, once again the term is computable recalling
that DOFs of type $i$ completely define
$(\bs{\varphi}_\beta \cdot \hat{\bs{n}}_{\partial E})$ on $\partial E$
and $m^E_\alpha$ is known.

It is useful to store the degrees of freedom of
$\nabla \cdot \bs{\varphi}^E_\alpha \in \mathbb{P}_{k_\nabla}(E)$. We define
the matrix
\begin{equation*}
  V^E = \left(H^E\right)^{-1} W^E \in
  \mathbb{R}^{n_{k_\nabla}^d\times n_{dof}^E} \,,
\end{equation*}
whose columns contain the coefficients of the polynomial decomposition
of $\nabla \cdot \bs{\varphi}^E_\alpha \in \mathbb{P}_{k_\nabla}(E)$,
$\alpha =1,\ldots, n^E_{dof}$.

$B^E \in \mathbb{R}^{dn_{k}^{d} \times n_{dof}^E}$ is crucial for the
computation since it involves integrating ``virtual'' shape functions,
which is a priori not possible. Its definition is
\begin{equation}
  [B^E]_{\alpha \beta} = \int_E \bs{g}_\alpha \cdot
  \bs{\varphi}_\beta^E \,\mathrm{d}E
  \quad \forall \alpha \in \{1,\ldots,dn_k^d\},\,
  \forall \beta \in\{1,\ldots,n_{dof}^E\}\,,
\end{equation}
which can be split into
$ B^E = \begin{bmatrix} B^{\nabla,E} \\ B^{\oplus,E} \end{bmatrix}$
where, $\forall \beta \in\{1,\ldots,n_{dof}^E\}$,
\begin{equation}
  \label{MatB}
  \begin{split}
    [B^{\nabla,E}]_{\alpha\beta} &= \int_E \bs{g}_\alpha^{\nabla,E}
    \cdot \bs{\varphi}_\beta^E\, \mathrm{d}E \quad \forall\alpha\in
    \{1,\ldots,n_{k,\nabla}^d\}\,,
    \\
    [B^{\oplus,E}]_{\alpha\beta} &= \int_E \bs{g}_\alpha^{\oplus,E}
    \cdot \bs{\varphi}_\beta^E\, \mathrm{d}E \quad \forall\alpha\in
    \{1,\ldots,n_{k,\oplus}^d\}\,,
  \end{split}
\end{equation}
with $B^{\nabla,E} \in \mathbb{R}^{n_{k,\nabla}^d \times n_{dof}^E}$
and $B^{\oplus,E} \in \mathbb{R}^{n_{k,\oplus}^d \times n_{dof}^E}$.
$B^{\oplus,E}$ can be computed using the DOFs of type $iii$ in
\eqref{DOF}, so that
\begin{equation*}
  [B^{\oplus,E}]_{\alpha\beta} = |E|
  \delta_{\beta\,,\,n^f_E n_{k}^{d-1}+n_{k,\nabla}^d+\alpha}
  \quad \forall \alpha\in\{1,\ldots,n^d_{k,\oplus}\},\,
  \beta\in\{1,\ldots,n^E{dof}\}\,.
\end{equation*}
$B^{\nabla,E}$ cannot be computed directly since
$\bs{g}_\alpha \in \left( \mathbb{P}_k(E) \right)^d$. Recalling
$\bs{g}^{\nabla,E}_{\alpha} = \nabla m^E_{\alpha+1}$ by definition and
using integration by parts,
\begin{equation}
  [B^{\nabla,E}]_{\alpha\beta} = - \int_E m^E_{\alpha +1} (\nabla
  \cdot \bs{\varphi}^E_\alpha)\, \mathrm{d}E + \int_{\partial E} m^E_{\alpha +1}
  (\bs{\varphi}^E_\beta \cdot \hat{\bs{n}}_{\partial E}) \,\mathrm{d} s :=
  [B_1^{\nabla,E}]_{\alpha\beta} + [B_2^{\nabla,E}]_{\alpha\beta} \,.
\end{equation}
$B_1^{\nabla,E}$ is computable using the known polynomial expression
of the divergence of the basis functions, represented by the matrix
$V^E$ already computed in the previous paragraph. The following
relationship is obtained:
\begin{equation*}
  B_1^{\nabla,E} = - H^{\#,E} V^E = - H^{\#,E} \left(
    \left( H^E \right)^{-1} W^E\right) = - H^{\#,E}
  \left(    H^E \right)^{-1} \left(W_1^E + W_2^E\right).
\end{equation*}
$B_2^{\nabla,E}$ is directly computable from the DOFs of type $i$
since it involves integration of a known polynomial over the faces of
the element.

Since $\left( \mathbb{P}_k(E) \right)^d \subset V_{k,k_\nabla}^{E,d}$,
it is possible to express the projector computed in next section as an
operator $V_{k,k_\nabla}^{E,d} \rightarrow V_{k,k_\nabla}^{E,d}$,
instead of
$V_{k,k_\nabla}^{E,d} \rightarrow \left( \mathbb{P}_k(E)
\right)^d$. For that purpose, the matrix
$D \in \mathbb{R}^{n_{dof}^E \times d n_{k}^d}$ is defined as
\begin{equation}
  \label{MatD}
  [D^E]_{\alpha\beta} = \text{DOF}_\alpha (\bs{g}^E_\beta)
  \quad \forall \alpha\in\{1,\ldots,n_{dof}^E\},\,
  \forall \beta\in\{1,\ldots,dn_k^d\} \,.
\end{equation}

\subsubsection{Computation of the local projector on polynomials}
\label{sec:proj}
The $\mathrm{L}^2$ projector
$\boldsymbol\Pi^{0,E}_k: V_{k,k_\nabla}^{E,d} \to \left(
  \mathbb{P}_k(E) \right)^d$ from (\ref{projector}) is defined as the
solution of the following linear system:
\begin{equation}
  \label{proysys}
  \int_E \boldsymbol\Pi^{0,E}_k \bs{v}_h \cdot \bs{g} \text{ dV} =
  \int_E \bs{v}_h \cdot \bs{g} \text{ dV} \quad
  \forall \bs{g} \in \left(\mathbb{P}_k(E) \right)^d, \quad
  \bs{v}_h \in V_{k,k_\nabla}^{E,d} \,.
\end{equation}
Specifically, for a basis function
$\bs{\varphi}^E_\alpha \in V_{k,k_\nabla}^{E,d}$ it will be now shown
how to compute this projection. First,
$\boldsymbol\Pi^{0,E}_k \bs{\varphi}_\alpha$ is expressed as a
polynomial using the decompositions of
$\left( \mathbb{P}_k(E) \right)^d$ shown previously:
\begin{equation}
  \label{proyexp}
  \bs\Pi^{0,E}_k \bs{\varphi}_\alpha^E =
  \sum_{\beta=1}^{d n_{k}^d} t^\alpha_\beta \bs{g}^E_\beta  =
  \sum_{\beta=1}^{n_{k,\nabla}^d} t^{\nabla,\alpha}_{\beta}
  \bs{g}_\beta^{\nabla,E}
  + \sum_{\beta=1}^{n_{k,\oplus}^d} t^{\oplus,\alpha}_{\beta}
  \bs{g}_\beta^{\oplus,E} \quad
  \forall \alpha\in\{1,\ldots,n^{E}_{dof}\}\,,
\end{equation}
where
$\bs{t}^\alpha= \begin{bmatrix} t^{\nabla,\alpha}_1 \,, & \ldots \,, &
  t^{\nabla,\alpha}_{n_{k,\nabla}^d} \,, & t^{\oplus,\alpha}_1 \,, &
  \ldots \,, &
  t^{\oplus,\alpha}_{n_{k,\oplus}^d} \end{bmatrix}^\intercal $ is the
column vector containing the coefficients expressing the combination
with respect to the polynomial basis for the projection of
$\bs{\varphi}_\alpha^E$, $\forall
\alpha\in\{1,\ldots,n^E_{dof}\}$. Replacing \eqref{proyexp} in
\eqref{proysys} the following linear system is obtained:
\begin{equation*}
  \begin{cases}
    \displaystyle \sum_{\beta=1}^{n_{k,\nabla}^d}
    t^{\nabla,\alpha}_{\beta} \int_E \bs{g}_\beta^{\nabla,E} \cdot
    \bs{g}^{\nabla,E}_\gamma \,\mathrm{d}E +
    \sum_{\beta=1}^{n_{k,\oplus}^d} t^{\oplus,\alpha}_{\beta} \int_E
    \bs{g}_\beta^{\oplus,E} \cdot \bs{g}^{\nabla,E}_\gamma
    \,\mathrm{d}E = \int_E \bs{\varphi}^E_\alpha \cdot
    \bs{g}^{\nabla,E}_\gamma \,\mathrm{d}E & \forall
    \gamma\in\{1,\ldots,n_{k,\nabla}^d\}\,,
    \\
    \displaystyle \sum_{\beta=1}^{n_{k,\nabla}^d}
    t^{\nabla,\alpha}_{\beta} \int_E \bs{g}_\beta^{\nabla,E} \cdot
    \bs{g}^{\oplus,E}_\gamma \,\mathrm{d}E +
    \sum_{\beta=1}^{n_{k,\oplus}^d} t^{\oplus,\alpha}_{\beta} \int_E
    \bs{g}_\beta^{\oplus,E} \cdot \bs{g}^{\oplus,E}_\gamma
    \,\mathrm{d}E = \int_E \bs{\varphi}^E_\alpha \cdot
    \bs{g}^{\oplus,E}_\gamma \,\mathrm{d}E & \forall
    \gamma\in\{1,\ldots,n_{k,\oplus}^d\}\,,
  \end{cases}
\end{equation*}
which, in view of \eqref{MatG} and \eqref{MatB}, can be rewritten as
\begin{equation*}
  G^E = \begin{bmatrix}
    G^{\nabla \nabla,E} & G^{\nabla \oplus,E}
    \\
    G^{\oplus \nabla,E} & G^{\oplus \oplus,E}
  \end{bmatrix}
  \bs{t}^\alpha= \begin{bmatrix}
    [B^{\nabla,E}]_{. \, \alpha} \\ [B^{\oplus,E} ]_{. \, \alpha}
  \end{bmatrix} \,,
\end{equation*}
so that $\bs{t}^\alpha = \left( G^E \right)^{-1} [B^E]_{. \, \alpha}$
(column $\alpha$ of $B^E$). Collecting all the vectors of coefficients
for $\alpha=1,...,n_{dof}^E$ we can define the projection matrix
$\hat{\Pi}^{0,E}_k \in \mathbb{R}^{d n_{k}^{d} \times n_{dof}^E}$
representing the operator acting from $V_{k,k_{\nabla}}^{E,d}$ to
$\left( \mathbb{P}_k(E) \right)^d$ as:
\begin{equation*}
  \hat{\Pi}^{0,E}_k =
  \begin{bmatrix}  \bs{t}^1 & \cdots & \bs{t}^{n_{dof}^E}
  \end{bmatrix} = \left( G^E \right)^{-1} B^E.
\end{equation*}
In order to obtain the matrix expression of the operator acting from
$V_{k,k_\nabla}^{E,d}$ into itself, we begin by expressing a
polynomial $\bs{g}^E_{\beta}$ as
\begin{equation*}
  \bs{g}^E_{\beta} =  \sum_{\gamma=1}^{n_{dof}^E} \text{DOF}_\gamma
  \left(\bs{g}^E_\beta \right)  \bs{\varphi}^E_\gamma
  = \sum_{\gamma=1}^{n^E_{dof}} [D^E]_{\gamma\beta} \, \bs\varphi^E_{\gamma}
  \quad \forall \beta \in \{1,\ldots,dn_{k}^{d}\}\,,
\end{equation*}
where we used the definition of the matrix $D^E$ given by
$\eqref{MatD}$. Replacing the above equation in \eqref{proyexp} yields
\begin{equation*}
  \bs\Pi^{0,E}_k \bs{\varphi}^E_\alpha = \sum_{\beta=1}^{d n_{k}^d}
  t^\alpha_\beta  \bs{g}^E_\beta  =
  \sum_{\beta=1}^{d n_{k}^d} t^\alpha_\beta
  \left( \sum_{\gamma=1}^{n^E_{dof}} [D^E]_{\gamma\beta} \,
    \bs\varphi^E_{\gamma}  \right)  =
  \sum_{\gamma=1}^{n^E_{dof}}  \sum_{\beta=1}^{d n_{k}^d}
  \left(  [D^E]_{\gamma\beta} \, t^\alpha_\beta
  \right) \bs\varphi^E_{\gamma} \quad
  \forall \alpha\in\{1,\ldots, n^E_{dof}\}\,,
\end{equation*}
thus, $\forall \alpha,\gamma\in\{1,\ldots,n^E_{dof}\}$,
\begin{equation*}
  \text{DOF}_\gamma\left( \bs\Pi^{0,E}_k \bs{\varphi}^E_\alpha \right)
  = [D^E]_{\gamma \, .} [ \hat{\Pi}^{0,E}_k ]_{. \, \alpha} \,.
\end{equation*}
Then we can define the matrix $\Pi^{0,E}_k\in\mathbb{R}^{n^E_{dof}}$
representing the $L^2$ projection seen as an operator from
$V^{E,d}_{k,k_\nabla}$ to itself as
\begin{equation*}
  \Pi^{0,E}_k = 
  \begin{bmatrix}
    \text{DOF}_1\left(\bs\Pi^{0,E}_k \bs{\varphi}^E_1\right) & \ldots
    & \text{DOF}_1\left(\bs\Pi^{0,E}_k
      \bs{\varphi}^E_{n^E_{dof}}\right)
    \\
    \vdots & \ddots & \vdots
    \\
    \text{DOF}_{n^E_{dof}}\left(\bs\Pi^{0,E}_k \bs{\varphi}^E_1\right)
    & \ldots & \text{DOF}_{n^E_{dof}}\left(\bs\Pi^{0,E}_k
      \bs{\varphi}^E_{n^E_{dof}}\right)
  \end{bmatrix} = D^E \hat{\Pi}^{0,E}_k \,.
\end{equation*}

A flow chart describing the interdependence of the matrices described
in sections \ref{sec:auxMat} and \ref{sec:proj} is provided in the
\textit{Supplementary Material} to the manuscript.

\subsubsection{Local stiffness matrices}
We are now ready to establish the matrix implementation of the
discrete equations.

The discrete bilinear form \eqref{bilform} is defined as
\begin{equation*}
  a_h^E \left(\bs{\varphi}^E_\beta,\bs{\varphi}^E_\alpha\right)
  = \left(\nu^d \boldsymbol{\Pi}^{0,E}_k \bs{\varphi}_\beta^E,
    \boldsymbol{\Pi}^{0,E}_k \bs{\varphi}_\alpha^E\right)_E
  + S^E(\bs{\varphi}^E_\beta - \boldsymbol{\Pi}^{0,E}_k \bs{\varphi}^E_\beta ,
  \bs{\varphi}^E_\alpha - \boldsymbol{\Pi}^{0,E}_k \bs{\varphi}^E_\alpha)
  := \left[ K_a^{E} \right]_{\alpha\beta} +
  \left[ K_s^{E} \right]_{\alpha\beta} \,.
\end{equation*}
In terms of the already computed matrices, we have
\begin{align*}
  K_a^E & = \left[\hat{\Pi}^{0,E}_k\right]^\intercal G^{\nu,E}
          \hat{\Pi}^{0,E}_k \,,
  &
    K_s^E & = \overline{\nu^d} |E| \left(\mathcal{I} -
            \Pi^{0,E}_k\right)^\intercal
            \left(\mathcal{I} - \Pi^{0,E}_k\right) \,,
\end{align*}
where $\mathcal{I}$ is the $n_{dof}^E \times n_{dof}^E$ identity
matrix.

Furthermore, the matrix arising from the terms in \eqref{weak:1} and
\eqref{weak:2} that involve the divergence of VEM basis functions has
already been computed as $W^E$, see \eqref{eq:MatW}, since we chose to
represent the pressure variable in the basis of scaled monomials in
$\mathfrak{M}_{k_\nabla}(E)$.

Finally, the local stiffness matrix $K^E$ on an element $E$ is given
by:
\begin{equation}
  \label{LocalK2}
  K^E = \begin{bmatrix}
    K_a^E + K_s^E & -\left(W^E\right)^\intercal
    \\
    W^E &    0
  \end{bmatrix} \,,
\end{equation}
with size
$\left( n^E_{dof} + n_{k_\nabla}^d \right) \times \left( n^E_{dof} +
  n_{k_\nabla}^d \right)$.

\subsection{Assembly of the global matrices}
\label{Imposing}

In this section we describe the structure of the global linear system
that arise from the discrete formulation of the problem. We denote by
$\bs{\varphi}^{d,\ell}_{\alpha}$ the velocity basis function with
index $\alpha$ defined on domain $\Omega^d_\ell$, and with
$\mu^{d,\ell}_\sigma$ the pressure basis function with index $\sigma$
on the same domain, with $\alpha = 1,\ldots, N^{d,\ell}_{u}$ and
$\sigma = 1,\ldots, N^{d,\ell}_p$. The total number of degrees of
freedom on domain $\Omega^{d}_\ell$ is denoted by
$N^{d,\ell}_{dof} = N^{d,\ell}_u + N^{d,\ell}_p$.

\subsubsection{Flux DOFs definition in the presence of fractures and
  traces}

\label{DOFduplication}
\begin{figure}[!h]
  \centering \hspace{\stretch{1}}
  \begin{subfigure}{.4\linewidth}
    \includegraphics[width=\linewidth]{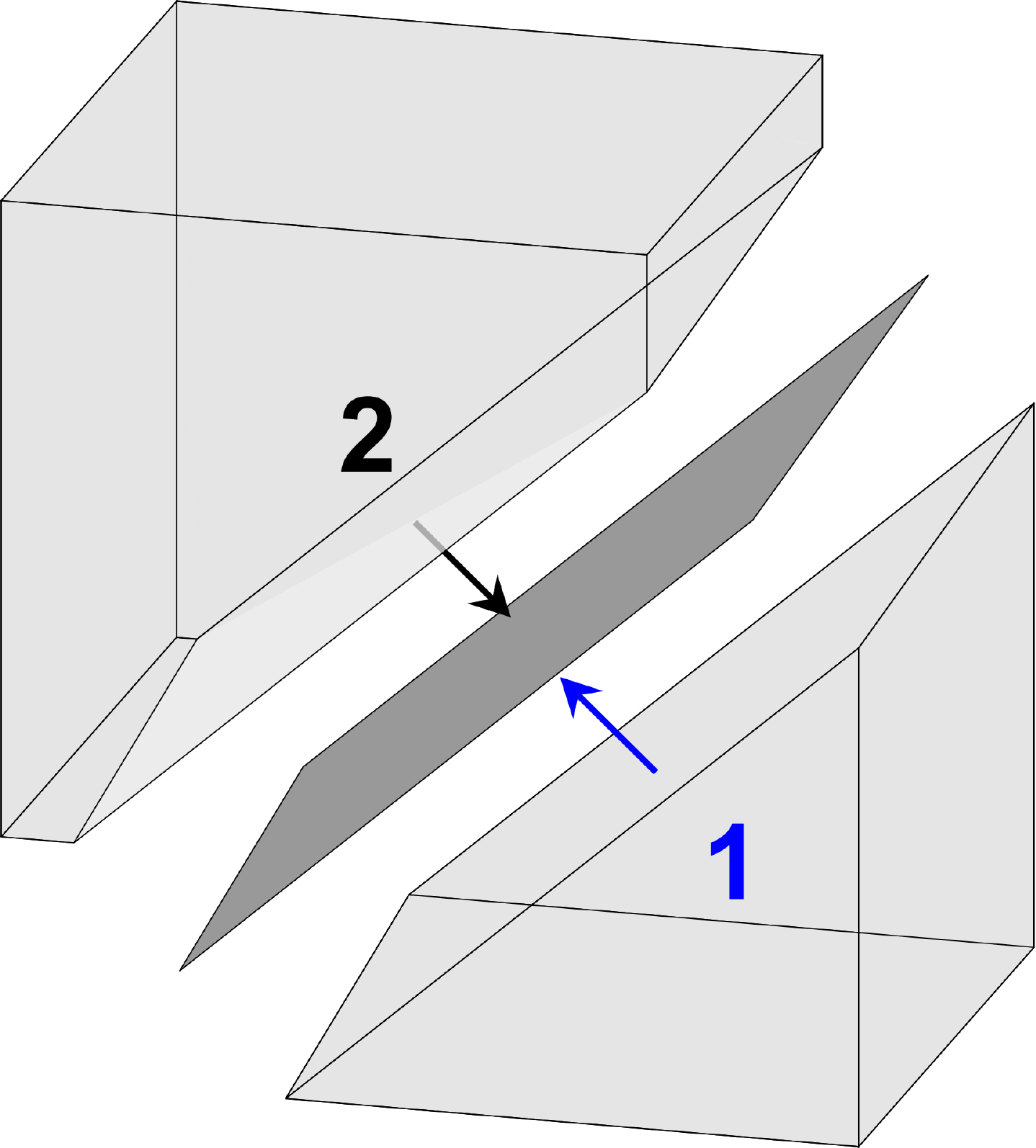}
    \caption{3D-2D coupling}
    \label{DoubleDOF:32}
  \end{subfigure}
  \hspace{\stretch{1}}
  \begin{subfigure}{.4\linewidth}
    \includegraphics[width=\linewidth]{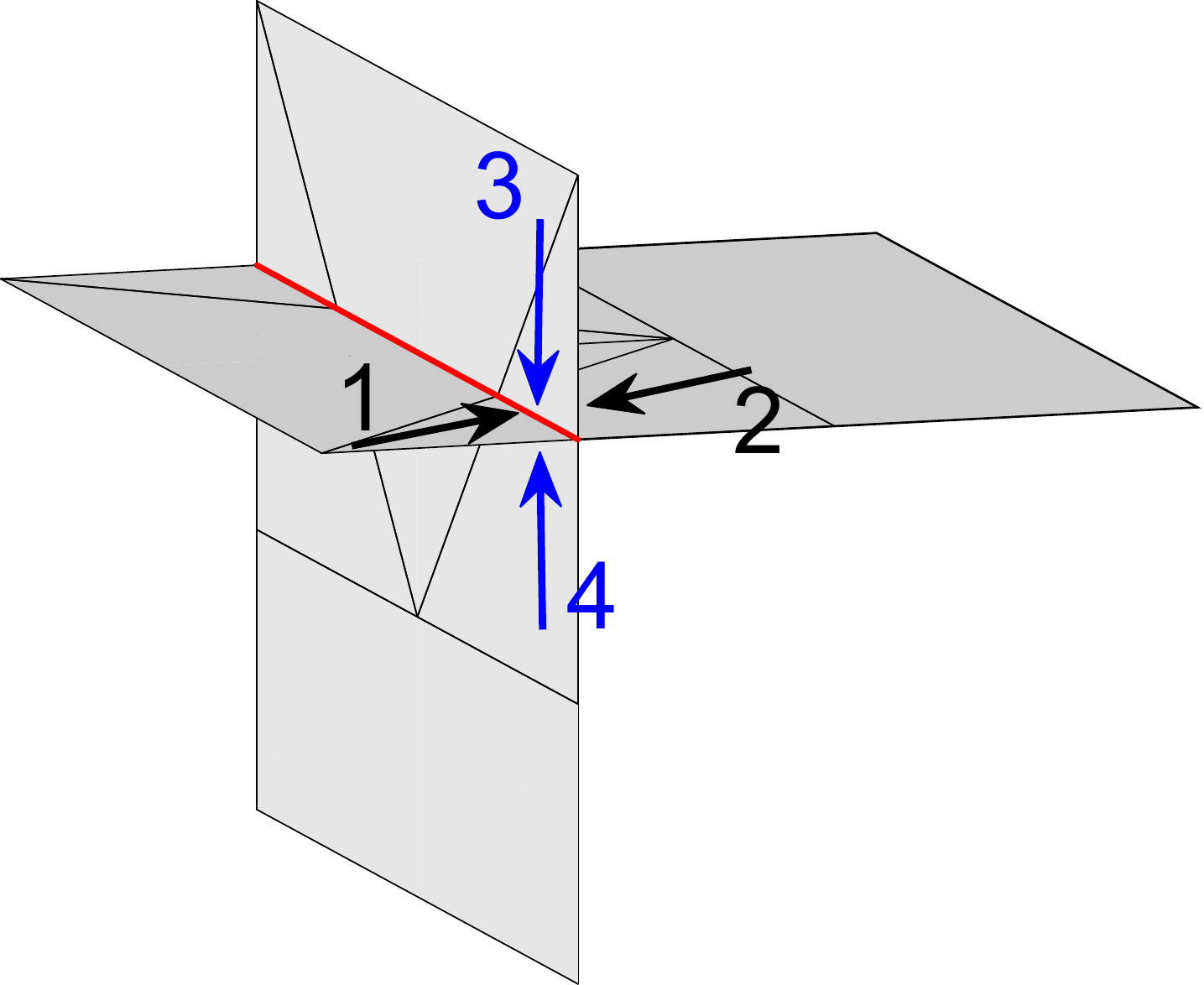}
    \caption{2D-1D coupling}
    \label{DoubleDOF:21}
  \end{subfigure}
  \hspace{\stretch{1}}
  \caption{Duplication of DOFs in the presence of lower dimensional
    objects}
  \label{DoubleDOF}
\end{figure}

For 3D elements that are intersected by planar fractures, a jump will
appear in the flux between adjacent faces that lie on a fracture,
since some of the flux leaving one face enters a 2D element as a
source term. In order to capture this phenomenon, flux DOFs of type
$i$ associated with a face on a fracture must be doubled, so that flux
continuity is no longer enforced and a jump can be represented.
Similarly, in the case of intersecting fractures whose intersection
defines a trace, we double the degrees of freedom on each edge. For
examples clarifying this point see Figure \ref{DoubleDOF}, where
RT0-VEM elements were chosen for simplicity, as they have a single DOF
of type $i$ per face/edge, although conceptually there is no
difference for any order and each DOF is duplicated in the same
way. In the first case (Figure \ref{DoubleDOF:32}), the fracture's
degrees of freedom result in a duplication of the type $i$ DOFs on the
face coinciding with the fracture. In the second case (Figure
\ref{DoubleDOF:21}), a trace segment doubles the flux DOFs on each
fracture, so that 4 flux DOFs are present on each segment that
discretizes the trace.

\subsubsection{3D elements (Matrix)}
Let $\bs{u}^{3} \in \mathbb{R}^{N_u^{3}}$,
$\bs{p}^{3} \in \mathbb{R}^{N_p^{3}}$ be the column vectors of flux
DOFs and pressure DOFs, collected in $\bs{h}^{3} :=
\begin{pmatrix}
  \bs{u}^3 \\ \bs{p}^3
\end{pmatrix}$ and let $\bs{f}^{3} \in\mathbb{R}^{N_{dof}^3}$ be the
vector of load values of the 3D domain (including terms arising from
non-homogeneous boundary conditions). In $\Omega^3$, the stiffness
matrix after assembly of the local matrices is
\begin{equation}
  \label{LocalK3}
  K^{3D} :=
  \left(\phantom{\begin{matrix}
        a_0\\ \ddots\\a_0\\b_0\\ \ddots\\b_0
      \end{matrix}} \right.\hspace{-1.5em} \underbrace{
    \begin{array}{cccc|}
      &  &  & \\
      &  A^{3}&  & \\
      &  &  &  \\ \hline
      &  &  & \\
      &  {D^{3}} &  & \\
      &  &  & \\
    \end{array}
  }_{N_u^3}
  \underbrace{
    \begin{array}{cccc}
      &  &  & \\
      &  -\left({D^{3}}\right)^\intercal&  & \\
      &  &  &  \\ \hline
      &  &  & \\
      & 0 &  & \\
      &  &  &
    \end{array}
  }_{N_p^3}
  \hspace{-1.5em}
  \left.\phantom{\begin{matrix}
        a_0\\ \ddots\\a_0\\b_0\\ \ddots\\b_0 \end{matrix}}\right)\hspace{-1em}
  \begin{tabular}{l}
    $\left.\lefteqn{\phantom{\begin{matrix} a_0\\ \ddots\\ a_0\ \end{matrix}}}\right\}N_u^3$\\
    $\left.\lefteqn{\phantom{\begin{matrix} b_0\\ \ddots\\ b_0\ \end{matrix}}} \right\}N_p^3$
  \end{tabular}
\end{equation}
where $A^{3}$ and $D^{3}$ are the matrices arising from the bilinear
forms. Namely
$[A^{3}]_{\alpha\beta}= (\nu^{3} \boldsymbol\Pi^0_k
\bs{\varphi}_\alpha^3,\boldsymbol\Pi^0_k
\bs{\varphi}_\beta^3)_{\Omega_3}$ and
$[D^{3}]_{\sigma\beta} = (\mu_\sigma^3,\nabla \cdot
\bs{\varphi}_\beta^3)_{\Omega_3}$, with $\alpha,\beta=1,...,N_u^3$ and
$\sigma=1,...,N_p^3$.

\subsubsection{2D elements (DFN)}
\label{imposing2D}
For every fracture $\Omega_\ell^2$, with $\ell=1,...,N^2$, recalling
\eqref{LocalK2} and after assembling the local stiffness matrices, the
stiffness matrix
$K^{2,\ell}\in \mathbb{R}^{N_{dof}^{2,\ell} \times N_{dof}^{2,\ell}}$
has the following structure:
\begin{equation}
  \label{LocalKDFN}
  K^{2,\ell} :=
  \left(\phantom{\begin{matrix}a_0\\ \ddots\\a_0\\b_0\\ \ddots\\b_0 \end{matrix}}
  \right.\hspace{-1.5em}
  \underbrace{
    \begin{array}{cccc|}
      &  &  & \\
      &  A^{2,\ell}&  & \\
      &  &  &  \\ \hline
      &  &  & \\
      &  D^{2,\ell} &  & \\
      &  &  & \\
    \end{array}
  }_{N_{u}^{2,\ell}}
  \underbrace{
    \begin{array}{cccc}
      &  &  & \\
      &  -\left(D^{2,\ell}\right)^\intercal&  & \\
      &  &  &  \\ \hline
      &  &  & \\
      & 0 &  & \\
      &  &  &
    \end{array}
  }_{  N_{p}^{2,\ell}}
  \hspace{-1.5em}
  \left.\phantom{\begin{matrix}a_0\\ \ddots\\a_0\\b_0\\ \ddots\\b_0 \end{matrix}}\right)\hspace{-1em}
  \begin{tabular}{l}
    $\left.\lefteqn{\phantom{\begin{matrix} a_0\\ \ddots\\ a_0\ \end{matrix}}}\right\} N_{u}^{2,\ell}$\\
    $\left.\lefteqn{\phantom{\begin{matrix} b_0\\ \ddots\\ b_0\ \end{matrix}}} \right\} N_{p}^{2,\ell}$
  \end{tabular}
\end{equation}
where $A^{2,\ell}$ and $D^{2,\ell}$ are the matrices arising from the
bilinear forms. Namely
$[A^{2,\ell}]_{\alpha\beta}= (\nu^{2}_\ell \boldsymbol\Pi^0_k
\bs{\varphi}_\alpha^{2,\ell},\boldsymbol\Pi^0_k
\bs{\varphi}_\beta^{2,\ell})_{\Omega^2_\ell}$, and
$[D^{2,\ell}]_{\sigma\beta} = (\mu_\sigma^{2,\ell},\nabla \cdot
\bs{\varphi}_\beta^{2,\ell})_{\Omega^2_\ell}$, with
$\alpha,\beta=1,..., N_{u}^{2,\ell}$ and
$\sigma=1,...,, N_{p}^{2,\ell}$.  The column vectors
$\bs{u}^{2}_\ell \in \mathbb{R}^{ N_{u}^{2,\ell}}$,
$\bs{p}^{2}_\ell \in \mathbb{R}^{ N_{p}^{2,\ell}}$ and
$\bs{f}^2_{\ell}\in\mathbb{R}^{N_{dof}^{2,\ell}}$ are the vectors of
flux DOFs, pressure DOFs and load values (including terms arising from
non-homogeneous boundary conditions) respectively. We define
$\bs{h}_{\ell}^2 :=\begin{pmatrix} \bs{u}^{2}_\ell \\ \bs{p}^{2}_\ell
\end{pmatrix}$ as the vector of values of the DOFs of the complete discrete
solution on fracture $\Omega_\ell^2$. We note that the matrix
$K^{2,\ell}$ is singular for fractures with pure Neumann boundary
conditions.  For the complete DFN we have:

\begin{center}
  $K^{2D} =
  \begin{pmatrix}
    K^{2,1} & 0 & \cdots & 0 \\
    0 & K^{2,2} & \cdots & \vdots \\
    \vdots  & \vdots  & \ddots & \vdots  \\
    0 & \cdots & \cdots & K^{2,N^2}
  \end{pmatrix}$ \quad , \qquad $ \bs{f}^{2D} =
  \begin{pmatrix}
    \bs{f}^2_{1}  \\
    \vdots \\
    \vdots \\
    \bs{f}^2_{N^2}
  \end{pmatrix}$ \quad and \quad $ \bs{h}^{2D} =
  \begin{pmatrix}
    \bs{h}^2_{1}   \\
    \vdots  \\
    \vdots  \\
    \bs{h}^2_{N^2}
  \end{pmatrix}$.
\end{center}

\begin{remark}
  At this point one can decide not to consider flow on traces and use
  a simpler coupling mechanism between fractures that simply relies on
  imposing flux continuity. By reason of the mesh being conforming across all
  traces, a simple linear constraint using Lagrange multipliers on
  DOFs enforces this situation. These multipliers in turn will provide
  an approximation of pressure head on traces. This coupling can be applied regardless of the presence of 3D elements, as it involves only fractures. For an implementation with the primal VEM
  formulation see \cite{DFNvem1} and for mixed VEM see
  \cite{VEMmixedDFN}.
\end{remark}

\subsubsection{1D elements (Trace network)}
For every trace $\Omega_{\ell}^1$, with $\ell=1,...,N^1$, the
stiffness matrix
$K^{1,\ell}\in \mathbb{R}^{N_{dof}^{1,\ell} \times N_{dof}^{1,\ell}}$
has the following structure, similarly to the 2D case:

\begin{equation}
  \label{LocalKT}
  K^{1,\ell} :=
  \left(\phantom{\begin{matrix}a_0\\ \ddots\\a_0\\b_0\\ \ddots\\b_0 \end{matrix}}
  \right.\hspace{-1.5em}
  \underbrace{
    \begin{array}{cccc|}
      &  &  & \\
      &  A^{1,\ell}&  & \\
      &  &  &  \\ \hline
      &  &  & \\
      &  D^{1,\ell} &  & \\
      &  &  & \\
    \end{array}
  }_{N_{u}^{1,\ell}}
  \underbrace{
    \begin{array}{cccc}
      &  &  & \\
      &  -\left(D^{1,\ell}\right)^\intercal&  & \\
      &  &  &  \\ \hline
      &  &  & \\
      & 0 &  & \\
      &  &  &
    \end{array}
  }_{N_{p}^{1,\ell} }
  \hspace{-1.5em}
  \left.\phantom{\begin{matrix}a_0\\ \ddots\\a_0\\b_0\\ \ddots\\b_0 \end{matrix}}\right)\hspace{-1em}
  \begin{tabular}{l}
    $\left.\lefteqn{\phantom{\begin{matrix} a_0\\ \ddots\\ a_0\ \end{matrix}}}\right\}N_{u}^{1,\ell}$\\
    $\left.\lefteqn{\phantom{\begin{matrix} b_0\\ \ddots\\ b_0\ \end{matrix}}} \right\}N_{p}^{1,\ell}$
  \end{tabular}
\end{equation}
where $A^{1,\ell}$ and $D^{1,\ell}$ are the matrices arising from the
bilinear forms. Namely
$[A^{1,\ell}]_{\alpha\beta}= (\nu^{1}_\ell
\bs{\varphi}_\alpha^{1,\ell},\bs{\varphi}_\beta^{1,\ell})_{\Omega^1_\ell}$,
and
$[D^{1,\ell}]_{\sigma\beta} = (\mu_\sigma^{1,\ell},\nabla \cdot
\bs{\varphi}_\beta^{1,\ell})_{\Omega^1_\ell}$, with
$\alpha,\beta=1,...,N_{u}^{1,\ell}$ and
$\sigma=1,\ldots,N^{1,\ell}_{p}$. The column vectors
$\bs{u}^{1}_\ell \in \mathbb{R}^{ N_{u}^{1,\ell}}$,
$\bs{p}^{1}_\ell \in \mathbb{R}^{ N_{p}^{1,\ell}}$ and
$\bs{f}^1_{\ell}\in\mathbb{R}^{N_{dof}^{1,\ell}}$ are the flux DOFs,
pressure DOFs and load values (including terms arising from
non-homogeneous boundary conditions), respectively. They are collected
in
$\bs{h}_{\ell}^1 :=\begin{pmatrix}\bs{u}^{1}_\ell,
  \bs{p}^{1}_\ell\end{pmatrix}$ as the vector of values of the DOFs of
the complete discrete solution on trace $\Omega_\ell^1$. Note that 1D
elements are polynomials in $\mathbb{P}_{k}$ and
$\mathbb{P}_{k_\nabla}$ and there is no polynomial projector involved
in the computations. Once again, matrix $K^{1,\ell}$ is singular for
traces with pure Neumann boundary conditions in Darcy's problem.  For
the complete trace network we have:
\begin{equation*}
  K^{1D} =
  \begin{pmatrix}
    K^{1,1} & 0 & \cdots & 0 \\
    0 & K^{1,2} & \cdots & \vdots \\
    \vdots  & \vdots  & \ddots & \vdots  \\
    0 & \cdots & \cdots & K^{1,N^1}
  \end{pmatrix} \quad , \qquad \bs{f}^{1D} =
  \begin{pmatrix}
    \bs{f}_{1}^1  \\
    \vdots \\
    \vdots \\
    \bs{f}_{N^1}^1
  \end{pmatrix} \quad \text{and} \quad \bs{h}^{1D} =
  \begin{pmatrix}
    \bs{h}_{1}^1   \\
    \vdots  \\
    \vdots  \\
    \bs{h}_{N^1}^1
  \end{pmatrix} \,.
\end{equation*}
Similarly to the DFN case in the previous Section (\ref{LocalKDFN}),
there is still no imposition of coupling conditions between the
different 1D domains. Flux balance at the intersection of traces is
imposed as explained in Section \ref{Coupling}.

\subsubsection{0D elements (Trace intersections)}
For each trace intersections $\Omega_{\ell}^0$ with
$\ell = 1,\ldots,N^0$ there is one pressure DOF $p^{0}_\ell$. The
equation on $\Omega_\ell^0$ is purely algebraic and the matrix arising
from it is treated in Section \ref{Coupling}, being of the same form
of the matrices arising from coupling conditions. Similarly to the
above sections, we define
\begin{equation*}
  \bs{f}^{0D} =
  \begin{pmatrix}
    f_{1}^0 \\ \vdots \\ \vdots \\ f_{N^0}^0
  \end{pmatrix}
  \quad \text{and} \quad
  \bs{h}^{0D} = \begin{pmatrix}
    p_{1}^0   \\
    \vdots  \\
    \vdots  \\
    p_{N^0}^0
  \end{pmatrix} \,.
\end{equation*}
\subsubsection{Final global matrix}
\label{Coupling}
The procedures described in the previous sections for introducing
coupling conditions on pure DFN problems and between 2D and 3D
elements can be readily combined into a single global
problem. Furthermore, the coupling conditions between fractures act on
edge DOFs of the flux variables while 2D-3D coupling establishes
conditions between flux DOFs on faces with internal 2D DOFs. So that,
in fact, the coupling conditions are ``decoupled'' in the final global
system.  The complete system with coupling conditions between 0D, 1D,
2D and 3D elements is as follows:
\begin{equation*}
  \label{COMPLETEsystem}
  \begin{bmatrix}
    K^{3D} + C^{3D/3D} & C^{3D/2D} & 0 & 0
    \\
    -\left(C^{3D/2D}\right)^\intercal & K^{2D} + C^{2D/2D} & C^{2D/1D}
    & 0
    \\
    0 & - \left( C^{2D/1D} \right)^\intercal & K^{1D} + C^{1D/1D} &
    C^{1D/0D}
    \\
    0 &0 & - \left( C^{1D/0D} \right)^\intercal & 0 &
  \end{bmatrix}
  \begin{bmatrix}
    \bs{h}^{3D} \\ \bs{h}^{2D} \\ \bs{h}^{1D} \\ \bs{h}^{0D}
  \end{bmatrix}
  =
  \begin{bmatrix}
    \bs{f}^{3D} \\ \bs{f}^{2D} \\ \bs{f}^{1D} \\ \bs{f}^{0D}
  \end{bmatrix}
  \,,
\end{equation*}
where $C^{xD/yD}$ are the coupling matrices. Specifically, for
$d = 1,\ldots,3$, and for each domain $\Omega_\ell^d$,
$\ell=1,\ldots,N^d$ (with the convention $\Omega^3_1 = \Omega^3$), we
define the matrix
$C^{dD/dD}_\ell \in \mathbb{R}^{N^{d,\ell}_{dof} \times
  N^{d,\ell}_{dof}}$ such that
\begin{gather*}
  \left[ C^{dD/dD}_\ell \right]_{\alpha\beta} =
  \begin{cases}
    -\sum_{j\in\mathcal{O}^{d^-}_\ell} \sum_{s=+,-}
    \left(\eta_j^{d-1}\right)^{-1}
    \left(\bs{\varphi}_\beta^{d,\ell}\cdot
      \nn^{\gammadd}_{\ell,j,s},\bs{\varphi}_\alpha^{d,\ell} \cdot
      \nn^{\gammadd}_{\ell,j,s}\right)_{\gammadd_{\ell,j,s}} & \forall
    \alpha,\,\beta = 1,\ldots,N^{d,\ell}_{u} \,,
    \\
    0 & \text{otherwise}\,,
  \end{cases}
\end{gather*}
and
\begin{equation*}
  C^{dD/dD} =
  \begin{bmatrix}
    C^{dD/dD}_1 & 0 & \cdots & 0
    \\
    0 & C^{dD/dD}_2 & \cdots & \vdots
    \\&&&\\
    \vdots & \ddots & \ddots & 0
    \\&&&\\
    0 & \cdots & 0 & C^{dD/dD}_{N^d}
  \end{bmatrix} \,.
\end{equation*}
Furthermore, let $\Omega^d_\ell$ be given and let
$j\in\mathcal{O}_\ell^{d^-}$. We define the
$N^{d,\ell}_{dof} \times N^{d-1,j}_{dof}$ matrix
$C^{dD/(d-1)D}_{\ell,j}$, coupling the velocity DOFs of $\Omega_\ell$
with the pressure DOFs of $\Omega_\ell$, as
\begin{equation*}
  \left[C^{dD/(d-1)D}_{\ell,j}\right]_{\beta\sigma} =
  \begin{cases}
    \left(\mu_\sigma,\jump{\bs{\varphi}_\beta^{d,\ell}\cdot
        \mathbf{n}_{\ell,j,\pm}^{\gamma^{(d,d-1)}}}\right)_{\Omega_j^{d-1}}
    & \forall \sigma =
    \left(N^{d-1,j}_{u}+1\right),\ldots,N^{d-1,j}_{dof} \,,\,
    \forall\beta = 1,\ldots, N^{d,\ell}_{u} \,,
    \\
    0 & \text{otherwhise}\,,
  \end{cases}
\end{equation*}
and the $N^{d-1}_{dof} \times N^d_{dof}$ matrix coupling d-dimensional
domains with (d-1)-dimensional domains, denoted by $C^{dD/(d-1)D}$, is
defined by its $N^{d-1} \times N^d$ blocks as follows:
\begin{equation*}
  \mathrm{block}_{j,\ell}\left(C^{dD/(d-1)D}\right) =
  \begin{cases}
    C^{dD/(d-1)D}_{j,\ell} & \text{if $j\in\mathcal{O}_\ell^{d^-}$}
    \\
    0 & \text{otherwhise}
  \end{cases} \quad \forall j = 1,\ldots,N^{d-1},\, \ell =
  1,\ldots,N^d \,.
\end{equation*}

\begin{remark}
  It is important to stress at this point that the use of a mixed
  formulation has a considerable edge over a primal formulation for
  coupling flux between dimensions.  Indeed, only local face
  geometrical information is needed to compute the coupling terms in
  the mixed formulation, since the flux DOFs in 3D and the pressure
  DOFs in 2D are internal to the face and do not interact with the
  analogous DOFs of neighbouring elements. This is not the case in the
  primal formulation: since pressure DOFs are assigned to edges and
  vertices, a single point in space may be part of any number of faces
  and fractures. Given the assumption that a trace is only shared by
  two fractures, this number reduces to 3 fractures. Nevertheless,
  these 3 fracture planes determine 8 subdivisions of the space
  according to which side of the fractures are being
  considered. Therefore, a single point in space may have 8 pressure
  DOFs in 3D. Analogously, that same point can be in 4 subdivisions of
  a fracture, which provides 12 more possible DOFs. Additionally,
  traces contribute another 6 possible DOFs and the trace intersection
  provides the final one resulting in 27 different pressure DOFs on a
  single spatial coordinate.  As expected, this introduces several
  implementation complexities that are completely irrelevant for the
  mixed formulation, such as DOF interaction between adjacent
  elements. In the special case of global continuity of pressure head,
  all pressure DOFs on a single location are equal and the coupling in
  the primal formulation becomes very straightforward.
\end{remark}

\begin{figure}[!h]
  \centering \includegraphics[width=0.48\linewidth]{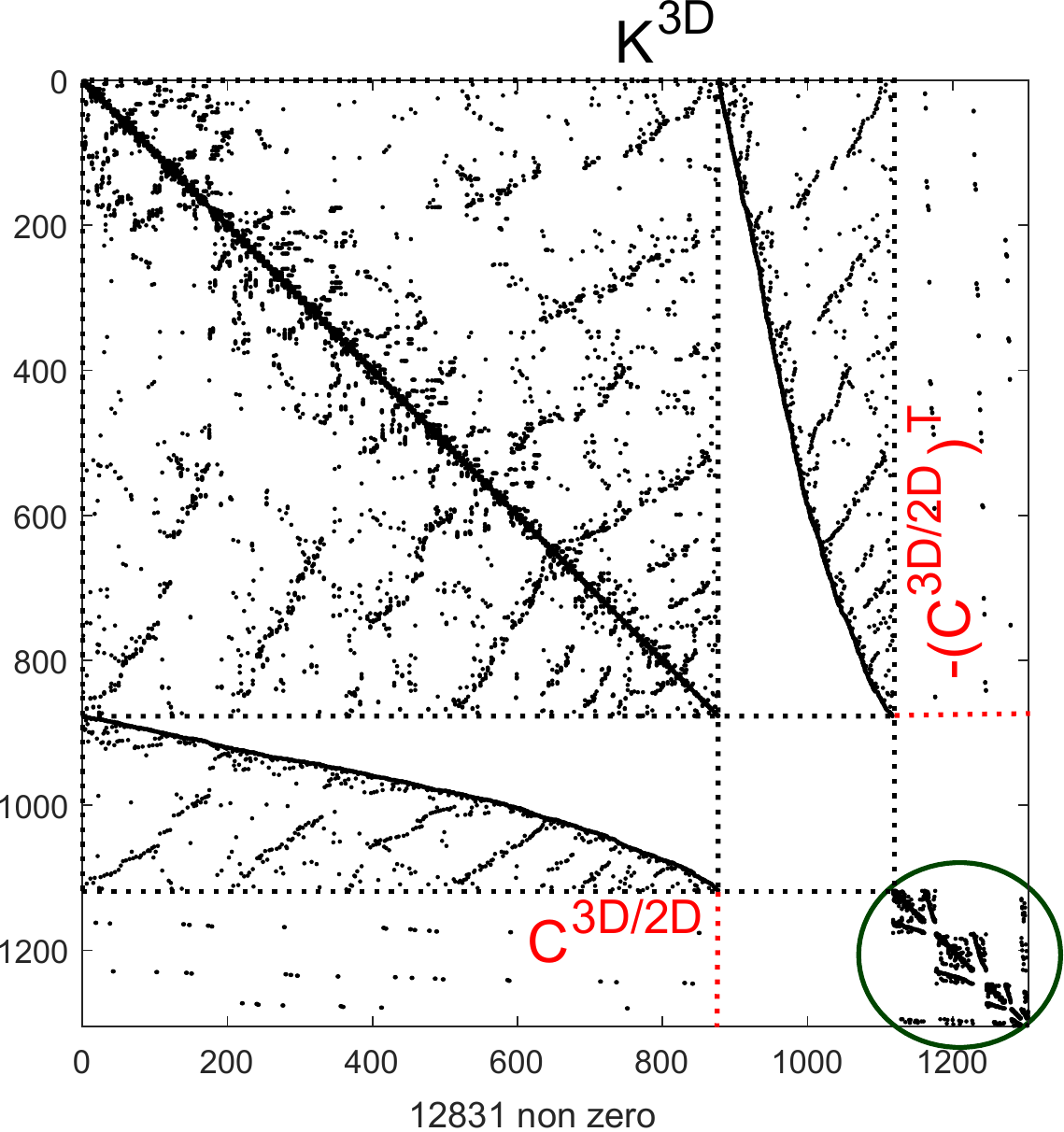}
  \includegraphics[width=0.48\linewidth]{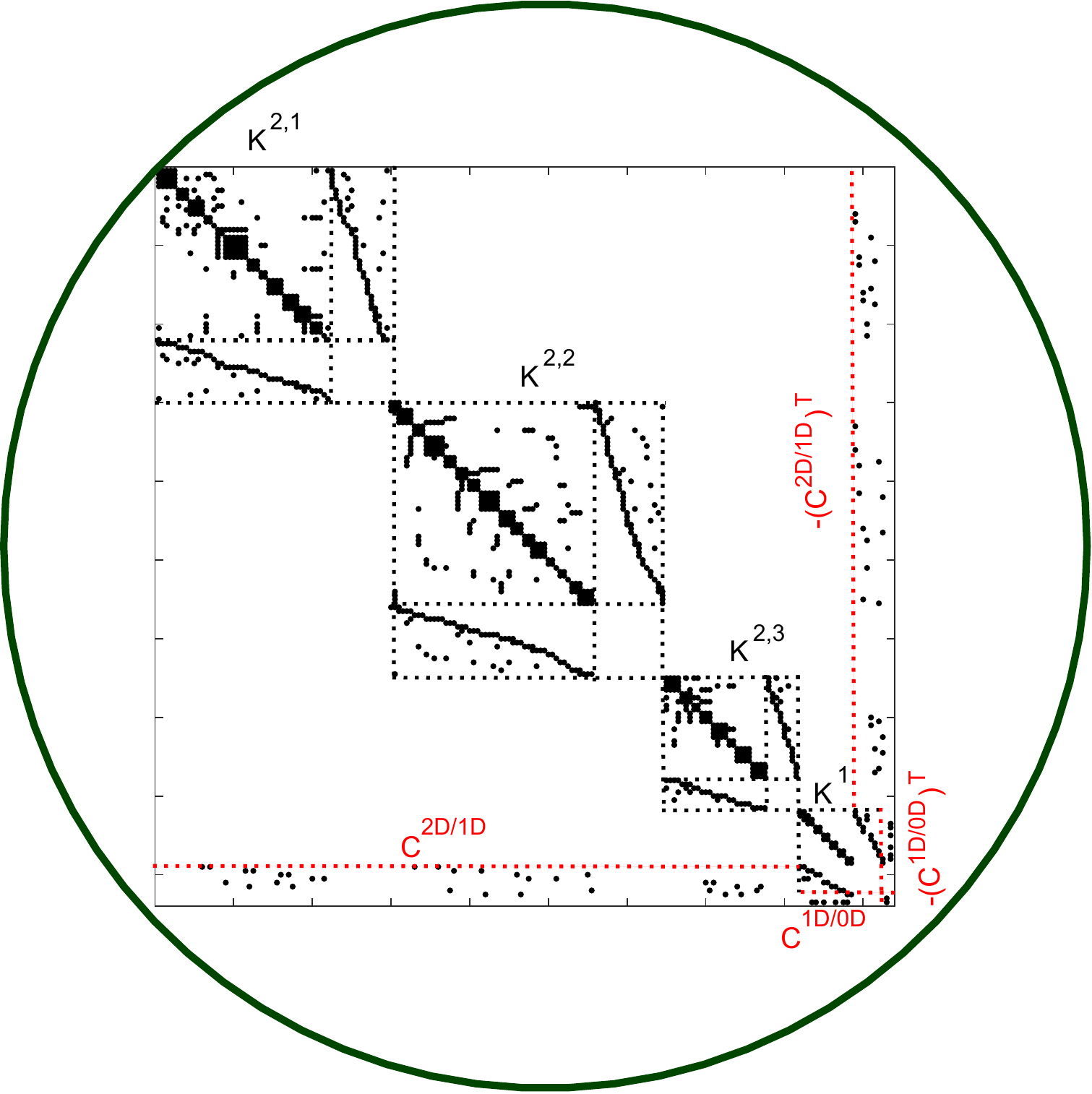}
  \caption{Linear system for the sample problem}
  \label{MatrixDFN3system}
\end{figure}

As an example, we present the complete (sparse) matrix arising for the
discretization with RT0-VEM elements of the previously presented
sample domain (Figure \ref{fig:nomenclature}). Dots in Figure
\ref{MatrixDFN3system} indicate non-null components of the matrix, and
they are associated with the pressure and the flux DOFs.  There are
877 3D flux DOFs, 240 3D pressure DOFs, 122 2D flux DOFs, 41 2D
pressure DOFs and 14 1D pressure DOFs, 8 pressure DOFs and a single
pressure DOF at the trace intersection. This gives 1303 DOFs in total,
but the linear system is very sparse and only 12831 non zero values
($< 1\%$) make up the complete stiffness matrix .


\section{Numerical results} 
\label{num_res}
This section contains a thorough numerical treatment of the proposed methodology. Every result and graphic presented has been produced with an in-house code in the Matlab programming language. Firstly, a benchmark problem  is proposed to assess the capacities of the approach applied to multidimensional problems. Afterwards, a set of complex problems with realistic embedded DFN geometries is analysed.  For the purpose of the numerical simulations we can assume that the coefficients and data parameters have already been homogenized and \textit{nondimensionalized}, meaning that the corresponding input data had their physical units removed and units have been adjusted by a suitable substitution of variables to simplify and parametrize problems.
Results for pure DFN problems of arbitrary order for general elliptic equations are presented in \cite{VEMmixedDFN} while numerical convergence results for pure 3D meshes of the Laplacian problem with constant coefficients for the first 4 orders of accuracy can be found in \cite{MixedDassi}. Several other results of mixed dimensional problems using mixed VEM formulations are presented in \cite{BenedettoThesis}. Also, since the exact solution in 1D networks can be easily obtained, they are not of particular interest. For all these reasons, convergence results and single-dimensional problems are omitted in this work and the focus is put solely on hybrid dimensional problems. Thus the results for this section involve the computation of the complete pressure and velocity fields across the whole domain comprised of 3D elements, 2D planar fractures, 1D traces and 0D points. 

\subsection{Problem 1: Benchmark problem with exact solution}
\label{Problem1}
\begin{figure}[!h]
	\centering
  \begin{subfigure}{.45\linewidth}
    \includegraphics[width=\linewidth]{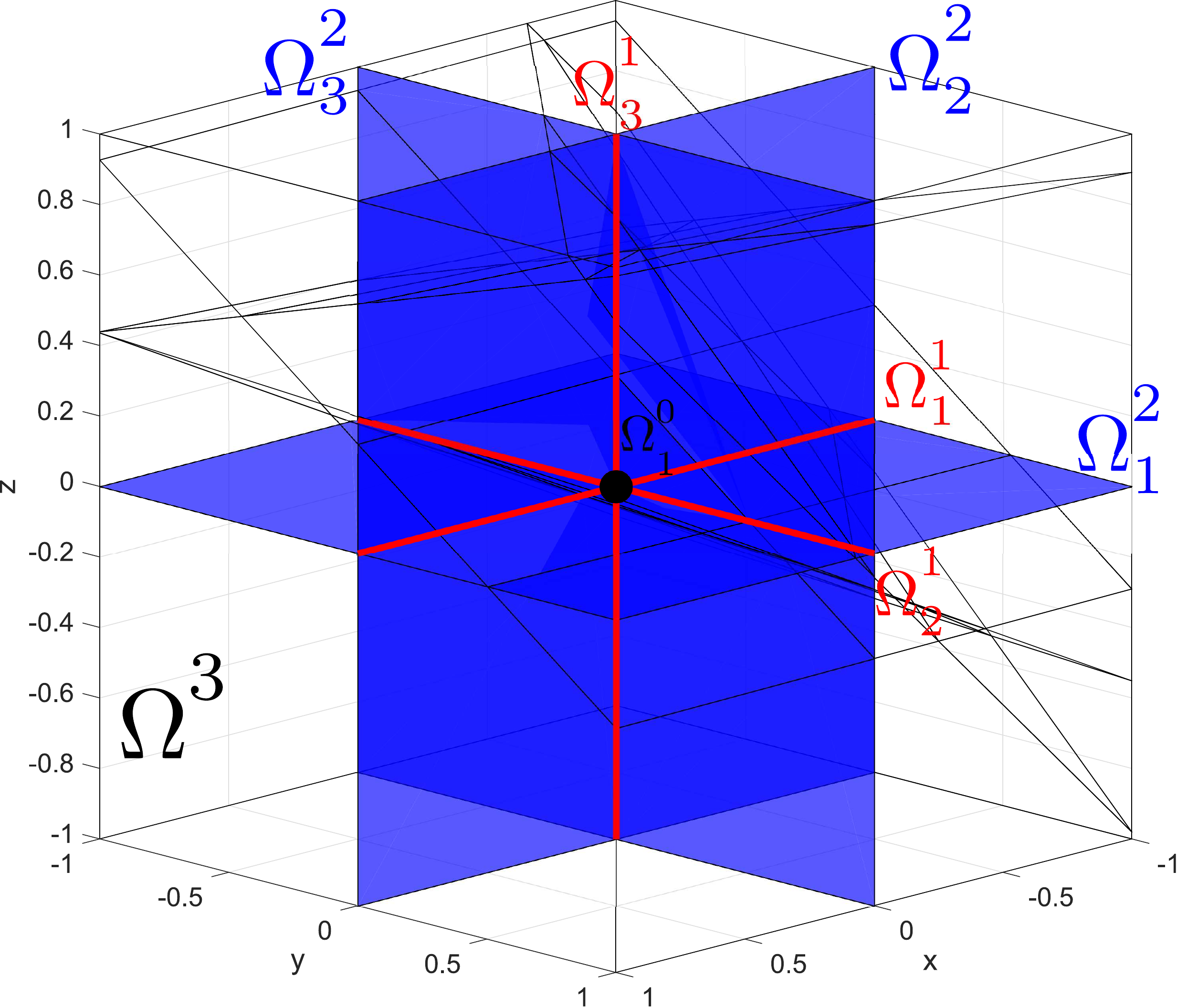}
    \caption{Geometry}
  \end{subfigure}
  \begin{subfigure}{.45\linewidth}
    \includegraphics[width=\linewidth]{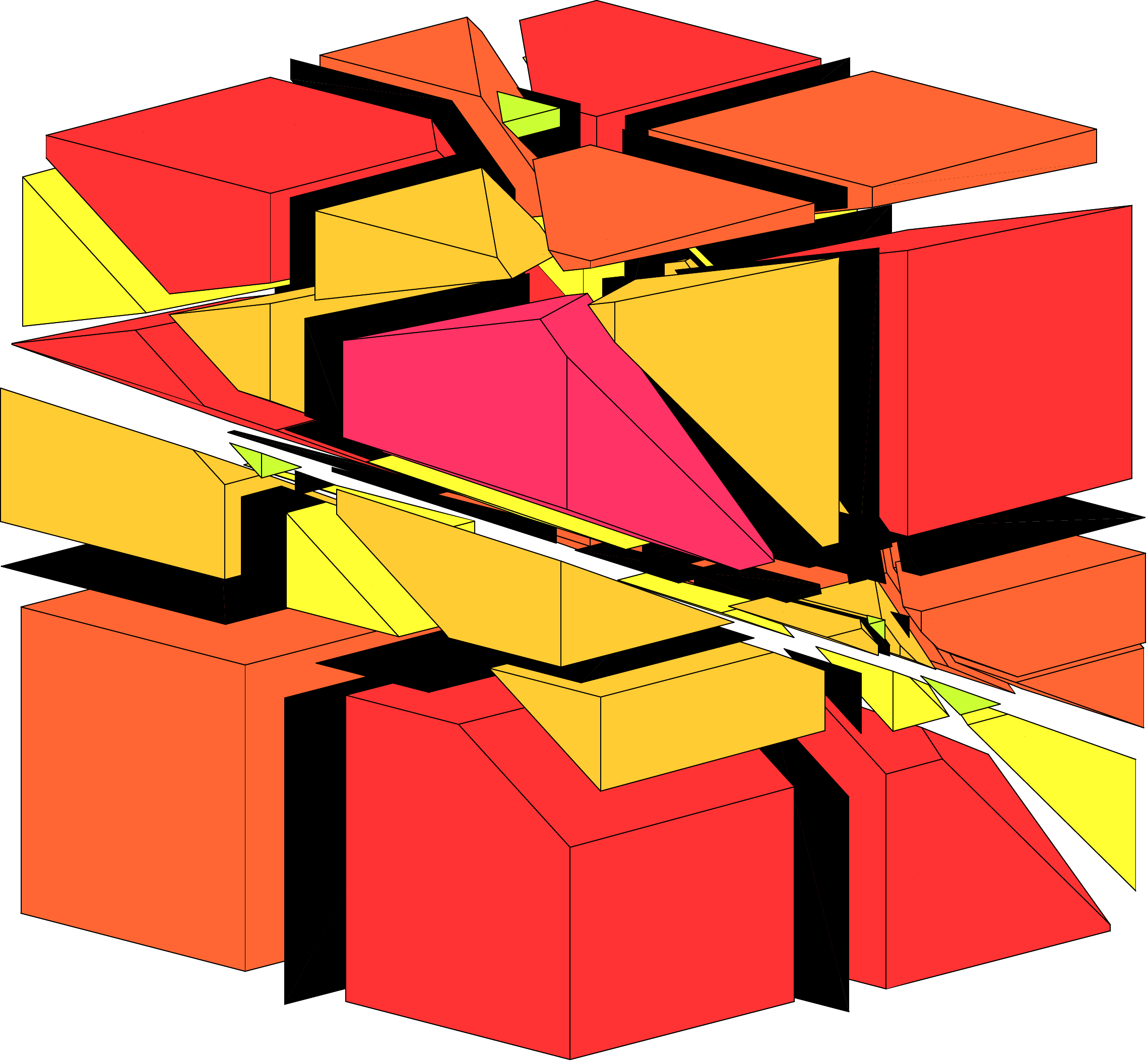}
    \caption{VEM mesh (exploded for visualization purposes)}
  \end{subfigure}
	\caption{Problem 1: Benchmark problem}
	\label{BenchmarkGeom}
\end{figure}

In order to assess the method taking into account the interaction
between different dimensions, a benchmark problem with exact solution
is analysed, which is manufactured to test all the capabilities of the
method and has little physical significance. The 3D domain is
$\Omega^3=\left[-1,1 \right]^3$ while the geometry for the 3 fractures
and traces making up the DFN are
\[ \left.
    \begin{array}{ll}
      \Omega^2_1 = \left\lbrace (x,y,z) : z=0, -1\leq x \leq 1, -1\leq y \leq 1 \right\rbrace \quad \quad \Omega^1_1 = \Omega^2_1 \cap \Omega^2_2 = \left\lbrace (x,y,z) : y=0, z=0 -1\leq x \leq 1 \right\rbrace \\ 
      \Omega^2_2 = \left\lbrace (x,y,z) : y=0, -1\leq x \leq 1, -1\leq z \leq 1 \right\rbrace \quad \quad \Omega^1_2 = \Omega^2_1 \cap \Omega^2_3 = \left\lbrace (x,y,z) : x=0, z=0, -1 \leq y \leq 1 \right\rbrace \\ 
      \Omega^2_3 = \left\lbrace (x,y,z) : x=0, -1\leq y \leq 1, -1\leq z \leq 1 \right\rbrace\quad \quad \Omega^1_3 = \Omega^2_2 \cap \Omega^2_3 = \left\lbrace (x,y,z) : x=0, y=0, -1\leq z \leq 1 \right\rbrace \\ 
    \end{array} 
  \right.
\]
and there is only one trace intersection at $\Omega^0_1 =
[0,0,0]$. Data for this problem is as follows:

\begin{align*}
  &\mathfrak{a}^3 = 1, && \eta^2_l \rightarrow \infty \ (l=1,2,3),\\
  &\mathfrak{a}^2_l = 2 \ (l=1,2,3), && \eta^1_l \rightarrow \infty \ (l=1,2,3),\\
  &\mathfrak{a}^1_l = 4 \ (l=1,2,3), && \eta^0_l \rightarrow \infty \ (l=1).
\end{align*} 

with non-homogeneous Dirichlet boundary conditions imposed in domain
boundaries for all dimensions, \ie, faces for 3D elements, edges for
2D elements, endpoints for 1D elements and point value for the trace
intersection. The loading terms for each domain are obtained from the
exact pressure solution which is globally chosen as a 4th degree
polynomial. Namely,
\begin{equation*}
  P(x,y,z) = (1+|x|)^4 + (1+|y|)^4 + (1+|z|)^4,
\end{equation*}
which is continuous in the whole domain, in agreement with an infinite
inter-dimensional permeability that prevents pressure jumps over
boundaries of elements of different dimensions.  The lowest order
element needed to obtain the exact solution is RT4-VEM, which has a
4th order pressure discretization.  At matrix/fracture intersections
the normal component of the flux variable is not continuous so that
there is flux exchange between matrix and fractures and analogously
between fractures and traces. Finally, there is also outgoing flux
associated with the single trace intersection. These flux exchanges
are accounted for in the forlumation by the coupling terms and
contribute to the loading terms of the lower dimensional entities.
The problem geometry and mesh are presented in Figure
\ref{BenchmarkGeom}, where matrix elements are transparent, 2D
elements on fractures are shown in blue, traces are depicted in red
and the single trace intersection is marked as a black circle. Some
artificial mesh cuts were introduced in the mesh to increase its
complexity. Element colouring indicate number of faces of the
polyhedron and a description of mesh composition is given in Table
\ref{Benchmarktable}, where it can be seen that RT4-VEM elements in 3D
are very computationally expensive due to their high number of DOFs
(node duplication included). Note that although this problem can be
solved with as few as 8 cubes and the exact solution is still
recovered, a more complex mesh is used to highlights the versatility
of a VEM discretization. The problem is pure Darcy flow with constant
diffusion coefficients and a globally continuous pressure, \ie,
infinite mixed-dimensional permeability.

\begin{table}[]
	\centering
	\begin{tabular}[t]{c|cccc|cc|cc|cc|c}
		\toprule
		Element Type & \multicolumn{4}{c|}{\#Elements} & \multicolumn{2}{c|}{\#3D DOFs} & \multicolumn{2}{c|}{\#2D DOFs} & \multicolumn{2}{c|}{\#1D DOFs}  & \multicolumn{1}{c}{\#0D DOFs}\\
		\midrule
		             & 3D & 2D & 1D & 0D & Flux & Pressure & Flux & Pressure & Flux & Pressure & Pressure \\
		\midrule
		RT4-VEM& 46& 47 &14 &1& 7194 & 1610 & 1853 & 705 & 76 & 70 & 1\\
		\bottomrule
	\end{tabular}
	\caption{Problem 1: discretization data}
	\label{Benchmarktable}
\end{table}

\begin{figure}[!ht]
  \centering
  \begin{subfigure}{.45\linewidth}
    \includegraphics[width=\linewidth] {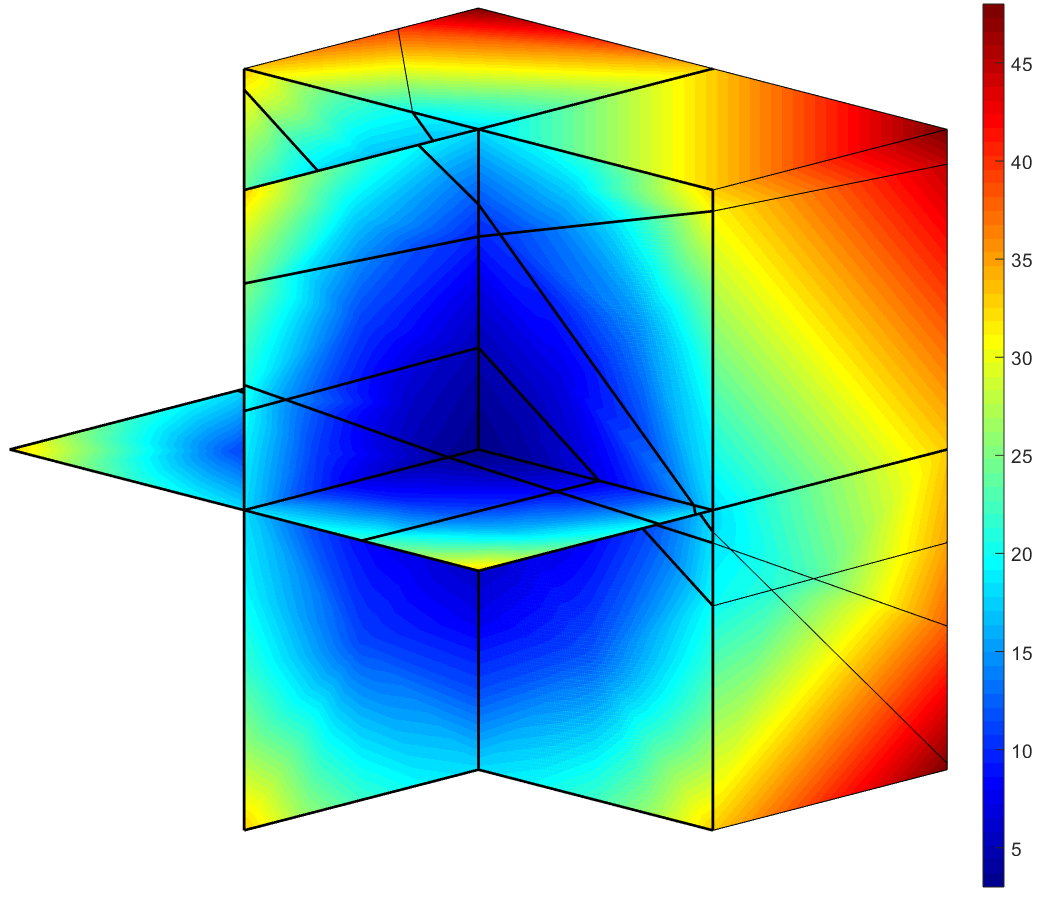}
    \caption{Global pressure head}
    \label{BenchmarkSol:1}
  \end{subfigure}
  \hfill
  \begin{subfigure}{.45\linewidth}
    \includegraphics[width=\linewidth] {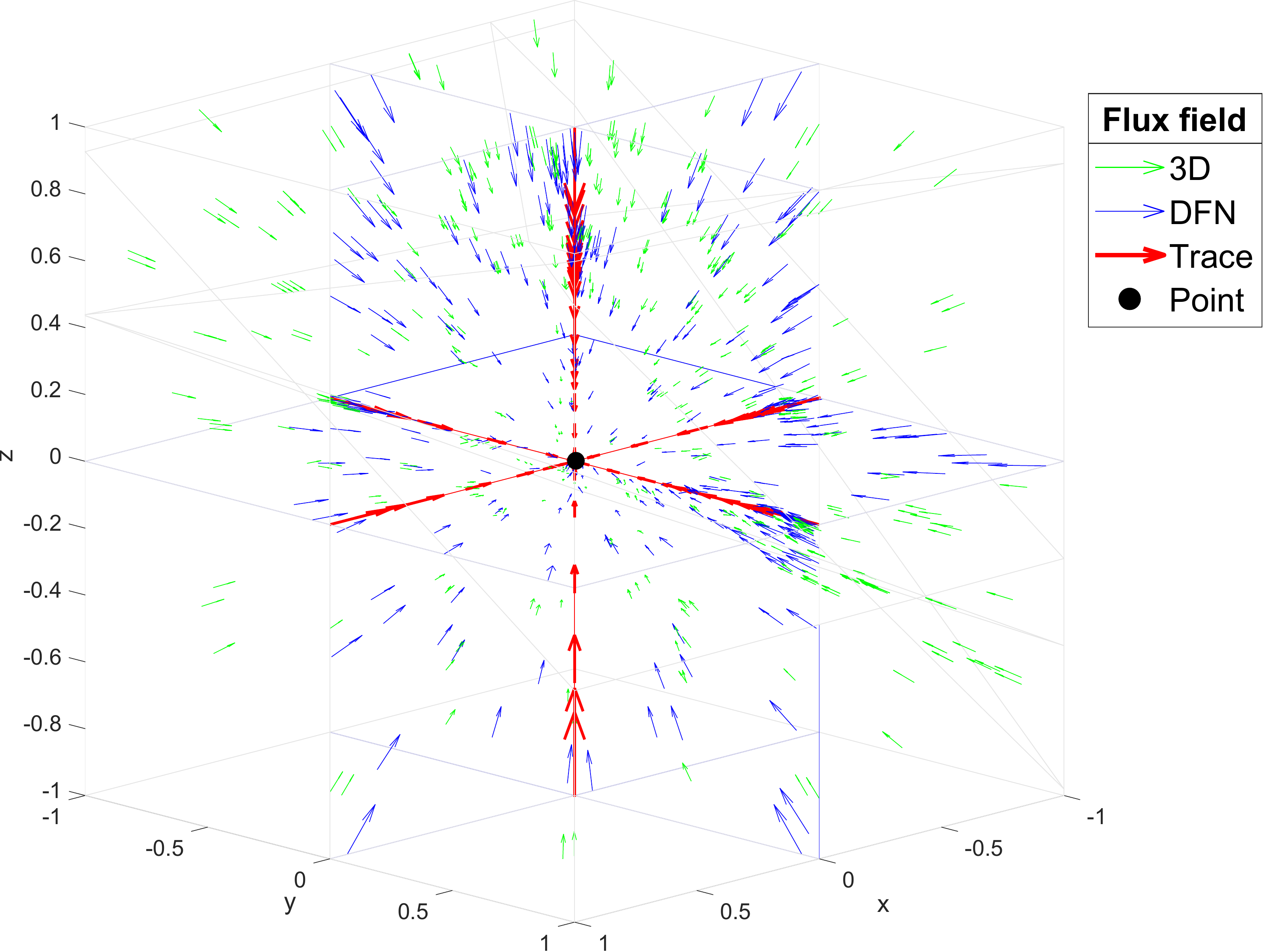}
    \caption{Global velocity field}
    \label{BenchmarkSol:2} %
  \end{subfigure}
  \hfill
  \begin{subfigure}{.45\linewidth}
    \includegraphics[width=\linewidth] {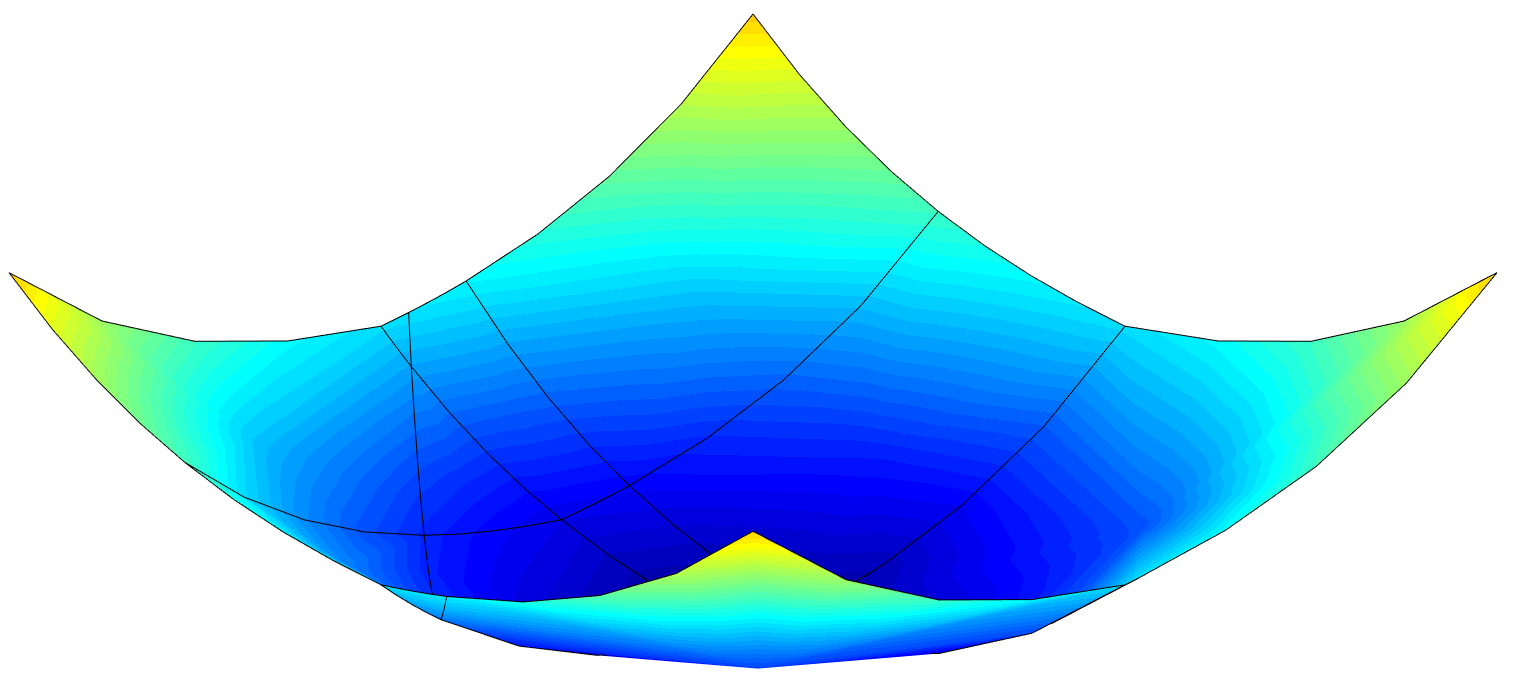}
    \caption{Pressure head on $\Omega^2_1$}
    \label{BenchmarkSol:3}%
  \end{subfigure}%
  \hfill
  \begin{subfigure}{.45\linewidth}
    \includegraphics[width=\linewidth] {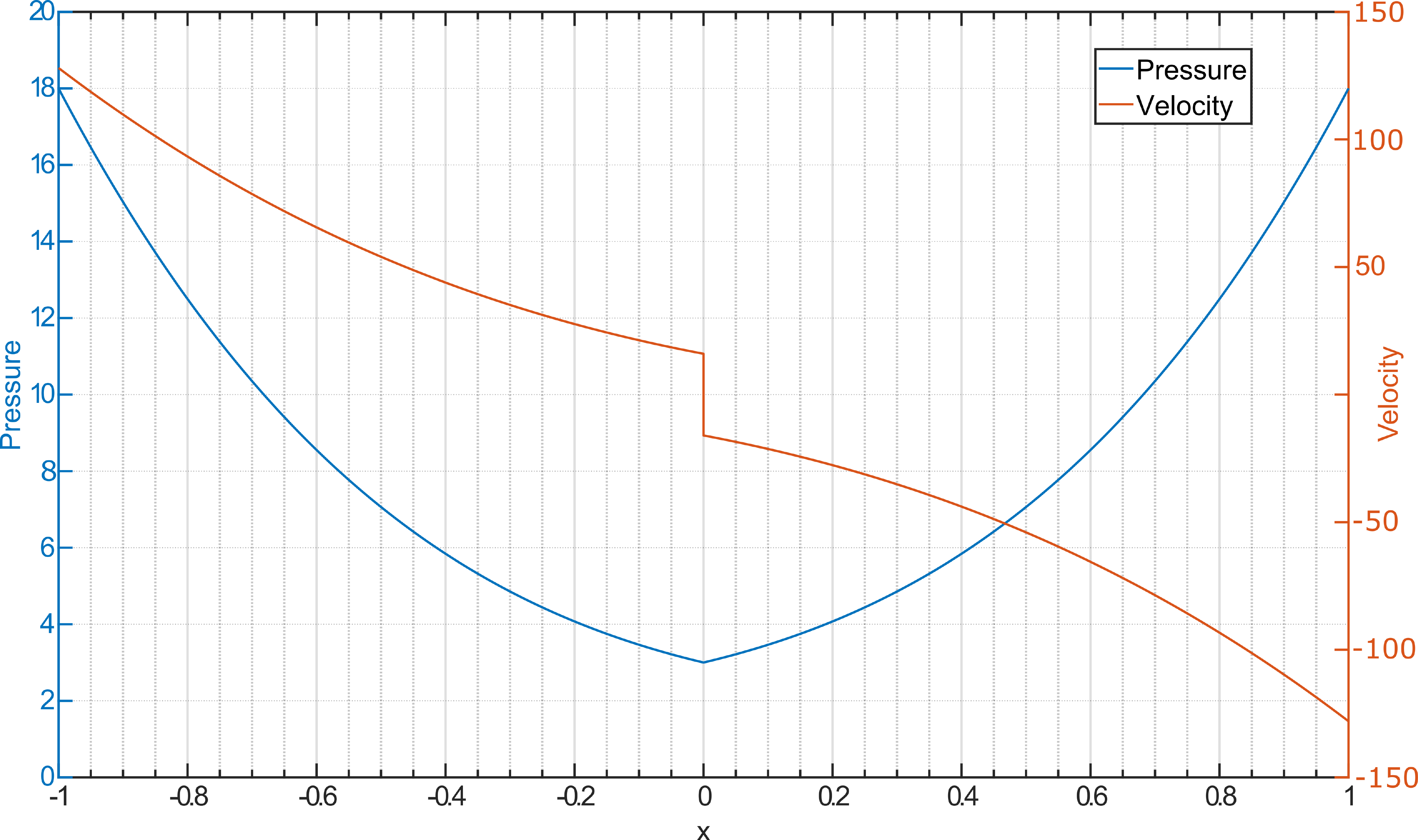}
    \caption{Pressure head and velocity field on $\Omega^1_1$}
    \label{BenchmarkSol:4} %
  \end{subfigure}%
  \caption{Problem 1: Discrete solution}
  \label{BenchmarkSol} 
\end{figure}

The global discrete solutions for the pressure head and the flux are
presented in Figures \ref{BenchmarkSol:1} and \ref{BenchmarkSol:2},
while the local solution on a fracture and a trace are given in
Figures \ref{BenchmarkSol:3} and \ref{BenchmarkSol:4}. In order to
verify the discrete solution, the quantities
$||P-P_h||_{\mathrm{L}^2}$,
$||\bs{u}-\Pi_k^0 \bs{u}_h||_{\mathrm{L}^2}$ and
$||(\nabla\cdot\bs{u})-(\nabla\cdot\bs{u}_h)||_{\mathrm{L}^2}$ are
computed and verified to vanish within numerical accuracy for all
domains, where $P$ and $\bs{u}$ are the exact solutions for any domain
and the subscript $h$ indicates a discrete solution. For the error in
the flux variable the projection $\Pi_k^0 \bs{u}_h$ is required since
the values of the discrete solution $\bs{u}_h$ are not explicitly
known inside a 3D or a 2D element.

\begin{figure}
	\centering
	\begin{tikzpicture}[->,>=stealth',shorten >=1pt,auto,node
    distance=2.5cm, thick,main
    node/.style={circle,draw,font=\sffamily\normalsize\bfseries},secondary
    node/.style={circle,scale=0.65,draw,font=\sffamily\normalsize\bfseries}]
	
    \node[main node] (3D) [color=green] {$\Omega^3$}; \node[secondary
    node] (S3D) [color=brown,node distance=2.5cm,below left of=3D]
    {$S$}; \node[secondary node] (BC3D) [color=brown,node
    distance=2.5cm,above left of=3D] {$BC$};
	
    \node[main node] (F2) [color=blue,right of=3D]{$\Omega^2_2$};
    \node[main node] (F1) [color=blue,above left of=F2]
    {$\Omega^2_1$}; \node[main node] (F3) [color=blue,below left
    of=F2] {$\Omega^2_3$};
	
    \node[secondary node] (SF1) [color=brown,node
    distance=1.75cm,above left of=F1] {$S$}; \node[secondary node]
    (BCF1) [color=brown,node distance=1.75cm,above right of=F1]
    {$BC$};

    \node[secondary node] (SF2) [color=brown,node
    distance=1.75cm,below left of=F2] {$S$}; \node[secondary node]
    (BCF2) [color=brown,node distance=1.75cm,above left of=F2] {$BC$};
	
    \node[secondary node] (SF3) [color=brown,node
    distance=1.75cm,below right of=F3] {$S$}; \node[secondary node]
    (BCF3) [color=brown,node distance=1.75cm,below left of=F3] {$BC$};

    \node[main node] (T2) [color=red,right of=F2]{$\Omega^1_2$};
    \node[main node] (T1) [color=red,above right of=T2]
    {$\Omega^1_1$}; \node[main node] (T3) [color=red,below right
    of=T2] {$\Omega^1_3$};

    \node[secondary node] (ST1) [color=brown,node distance=1.5cm,above
    right of=T1] {$S$}; \node[secondary node] (BCT1) [color=brown,node
    distance=1.5cm,above left of=T1] {$BC$};

    \node[secondary node] (ST2) [color=brown,node distance=1.5cm,above
    right of=T2] {$S$}; \node[secondary node] (BCT2) [color=brown,node
    distance=1.5cm,below right of=T2] {$BC$};
	
    \node[secondary node] (ST3) [color=brown,node distance=1.5cm,below
    left of=T3] {$S$}; \node[secondary node] (BCT3) [color=brown,node
    distance=1.5cm,below right of=T3] {$BC$};
	
    \node[main node] (I) [color=black,below right of=T1] {$\Omega^0_1$};
    \node[secondary node] (BCP) [color=brown, right of=I] {$BC$};

    \path[every node/.style={font=\sffamily\small}] (3D) edge
    [green,bend left] node[right] {32} (F1) (3D) edge [green]
    node[above] {32} (F2) (3D) edge [green, bend right] node[right]
    {32} (F3)

    (BC3D) edge [brown] node[below] {\tiny 768} (3D) (3D) edge [brown]
    node[above] {\tiny 672} (S3D)
	
    (F1) edge [blue] node[above] {32} (T1) (F1) edge [blue]
    node[above] {32} (T2) (F2) edge [blue] node[above] {32} (T1) (F2)
    edge [blue] node[above] {32} (T3) (F3) edge [blue] node[above]
    {32} (T2) (F3) edge [blue] node[above] {32} (T3)
	
    (BCF1) edge [brown] node[below right] {\tiny 512} (F1) (F1) edge
    [brown] node[below left] {\tiny 480} (SF1)
	
    (BCF2) edge [brown] node[above right] {\tiny 512} (F2) (F2) edge
    [brown] node[below right] {\tiny 480} (SF2)
	
    (BCF3) edge [brown] node[above left] {\tiny 512} (F3) (F3) edge
    [brown] node[above right] {\tiny 480} (SF3)
	
    (T1) edge [brown] node[below left] {\tiny 256} (BCT1) (ST1) edge
    [brown] node[below right] {\tiny 224} (T1)
	
    (T2) edge [brown] node[below left] {\tiny 256} (BCT2) (ST2) edge
    [brown] node[above left] {\tiny 224} (T2)
	
    (T3) edge [brown] node[above right] {\tiny 256} (BCT3) (ST3) edge
    [brown] node[above left] {\tiny 224} (T3)
	
    (T1) edge [red,bend left] node[left] {32} (I) (T2) edge [red]
    node[above] {32} (I) (T3) edge [red, bend right] node[left] {32}
    (I)
	
    (I) edge [black] node[above] {96} (BCP);
	\end{tikzpicture}
	\caption{Problem 1: Flux chart}
	\label{BenchmarkFluxChart}
\end{figure}
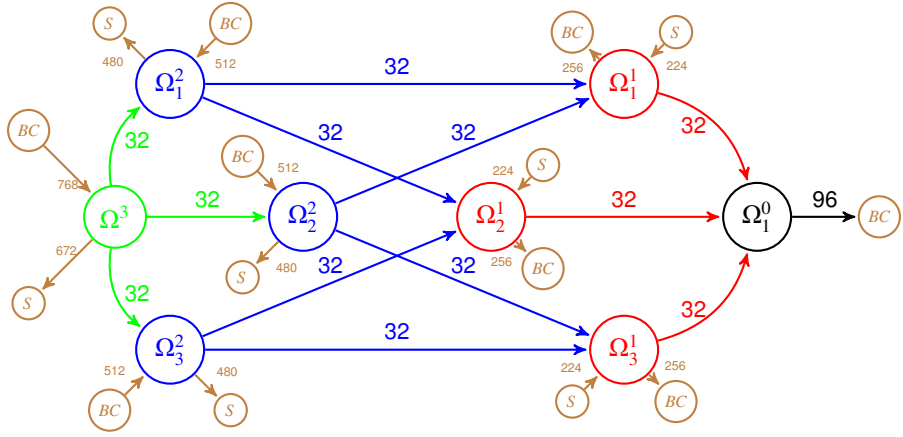

Finally, this benchmark serves to show how a mixed VEM formulation (of
any order) strongly imposes flux conservation and always displays a
null flux mismatch up to machine precision. The flux chart in Figure
\ref{BenchmarkFluxChart} provides a visual example of this and
exhibits how flux travels throughout the whole multidimensional system
including boundary conditions ('BC') and source/sink terms('S').

\subsection{Problem 2: Embedded DFN with trace flow and finite
  inter-dimensional normal permeability}
\label{Problem2}
\begin{figure}[!h]
	\centering
	\begin{subfigure}{.45\linewidth}
    \includegraphics[width=\linewidth] {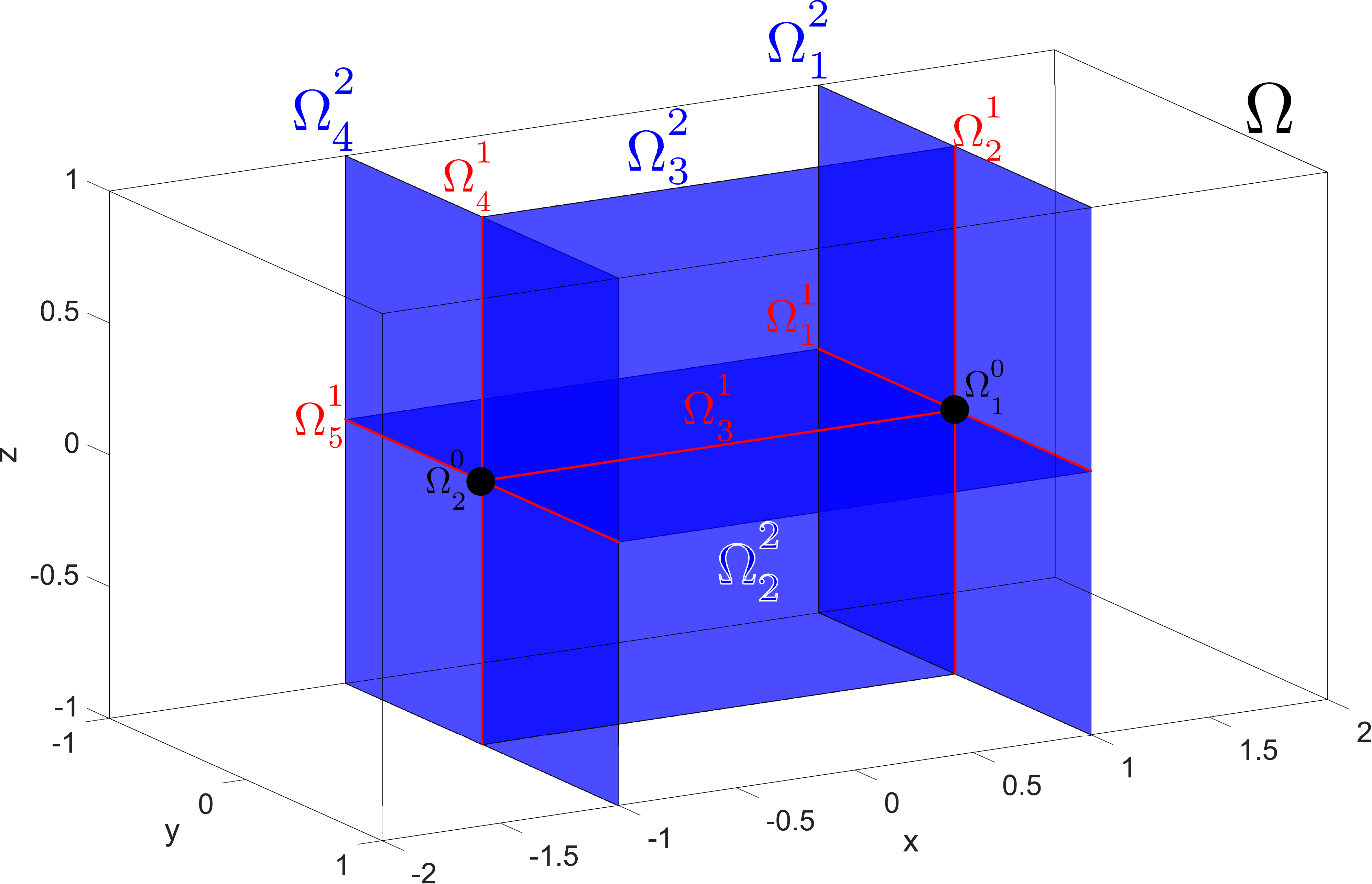}
    \caption{Geometry and notation}
    \label{Problem2:geometry}
	\end{subfigure}
	\begin{subfigure}{.45\linewidth}
    \includegraphics[width=\linewidth] {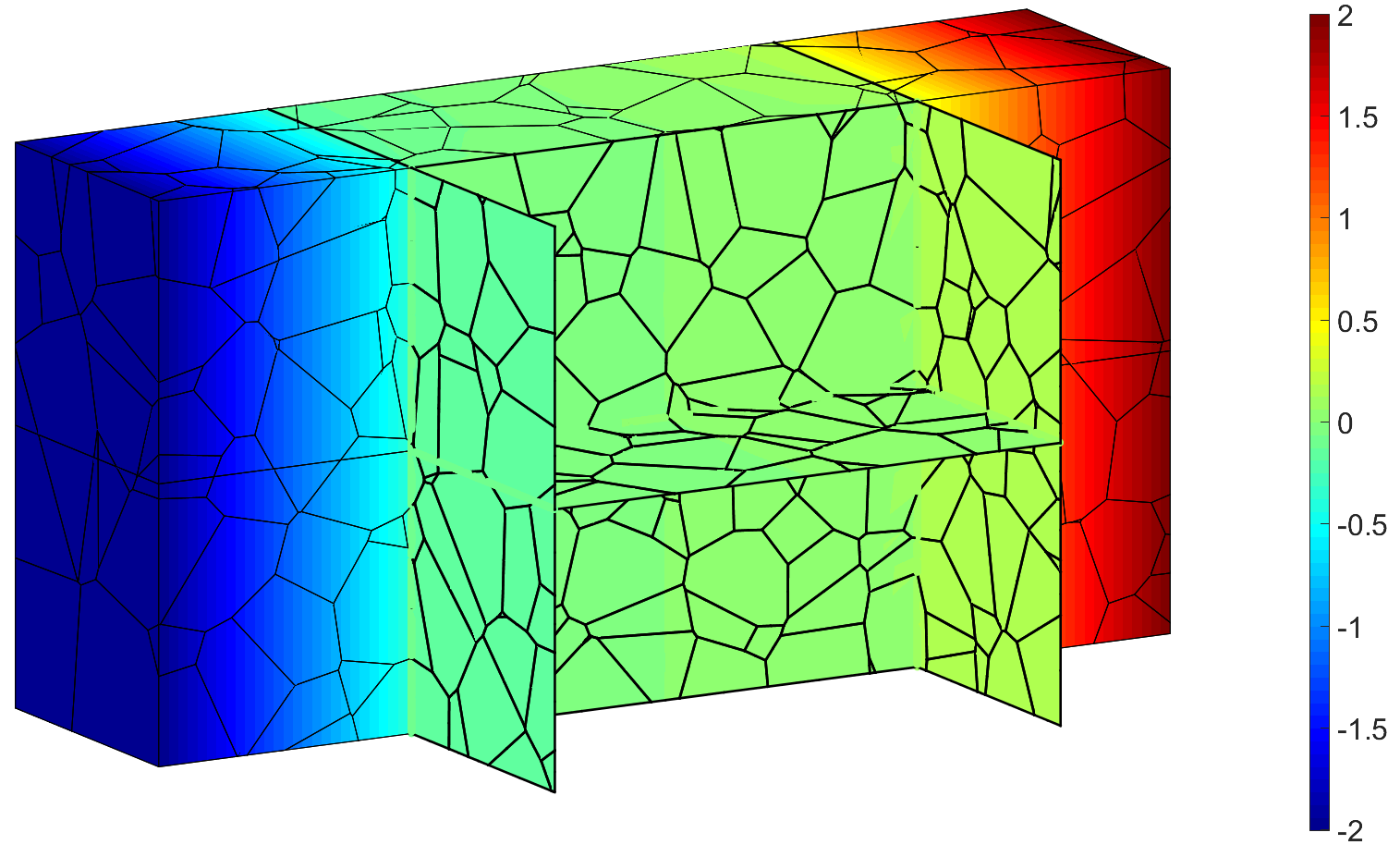}
    \caption{Global pressure field (removed elements for
      visualization)}
    \label{Problem2:Pfield}%
	\end{subfigure}
	\caption{Problem 2}
	\label{Problem2_solutions}
\end{figure}

This problem was chosen to highlight the coupling between entities of
different dimensions and to determine jumps in pressure that arise
from finite normal permeability coefficients. In addition, lower
dimensional objects are given much higher permeability values than
their higher dimensional counterpart, resulting in a flow distribution
that favours traces over fractures, and fractures over matrix.  The
problem consists of a 3D domain
$\Omega^3 = \lbrace -2 \leq x \leq 2,-1 \leq y \leq 1,-1 \leq z \leq 1
\rbrace$, 4 fractures ($\Omega^2_1,...,\Omega^2_4$) that give rise to
5 traces ($\Omega^1_1,...,\Omega^1_5$) that intersect in two points
($\Omega^0_1,\Omega^0_2$). The notation and geometry are shown in
Figure \ref{Problem2:geometry}. Note that trace $\Omega^1_3$ begins
and ends precisely on the surface of $\Omega^2_1$ and $\Omega^2_4$
respectively. Data for this problem is as follows:
\begin{align*}
  &\mathfrak{a}^3 = 10^0, && \eta^2_l = 10 \ (l=1,...,4),\\
  &\mathfrak{a}^2_l = 10^2 \ (l=1,...,4), && \eta^1_l = 10 \ (l=1,...,5),\\
  &\mathfrak{a}^1_l = 10^4 \ (l=1,...,5), && \eta^0_l = 10 \ (l=1,2).
\end{align*}

Boundary conditions are taken as a fixed pressure $P=-2$ on
$\lbrace{x=-2}\rbrace$ and $P=2$ on $\lbrace{x=2}\rbrace$. Homogeneous
Neumann boundary conditions are imposed on all other boundaries for
all entities.  The problem was solved using RT2-VEM elements and a
polyhedral mesh as detailed in Table
\ref{P2table}. 
\begin{table}[]
	\centering
	\begin{tabular}[t]{c|cccc|cc|cc|cc|c}
		\toprule
		Element Type & \multicolumn{4}{c|}{\#Elements} & \multicolumn{2}{c|}{\#3D DOFs} & \multicolumn{2}{c|}{\#2D DOFs} & \multicolumn{2}{c|}{\#1D DOFs}  & \multicolumn{1}{c}{\#0D DOFs}\\
		\midrule
		             & 3D & 2D & 1D & 0D & Flux & Pressure & Flux & Pressure & Flux & Pressure & Pressure \\
		\midrule
		RT2-VEM& 869& 315 &51 &2& 7194 & 48244 & 5391 & 1890 & 162 & 153 & 2\\
		\bottomrule
	\end{tabular}
	\caption{Problem 2: discretization data}
	\label{P2table}
\end{table}

\begin{figure}[!h]
	\centering
	\begin{subfigure}{.8\linewidth}
    \includegraphics[width=\linewidth] {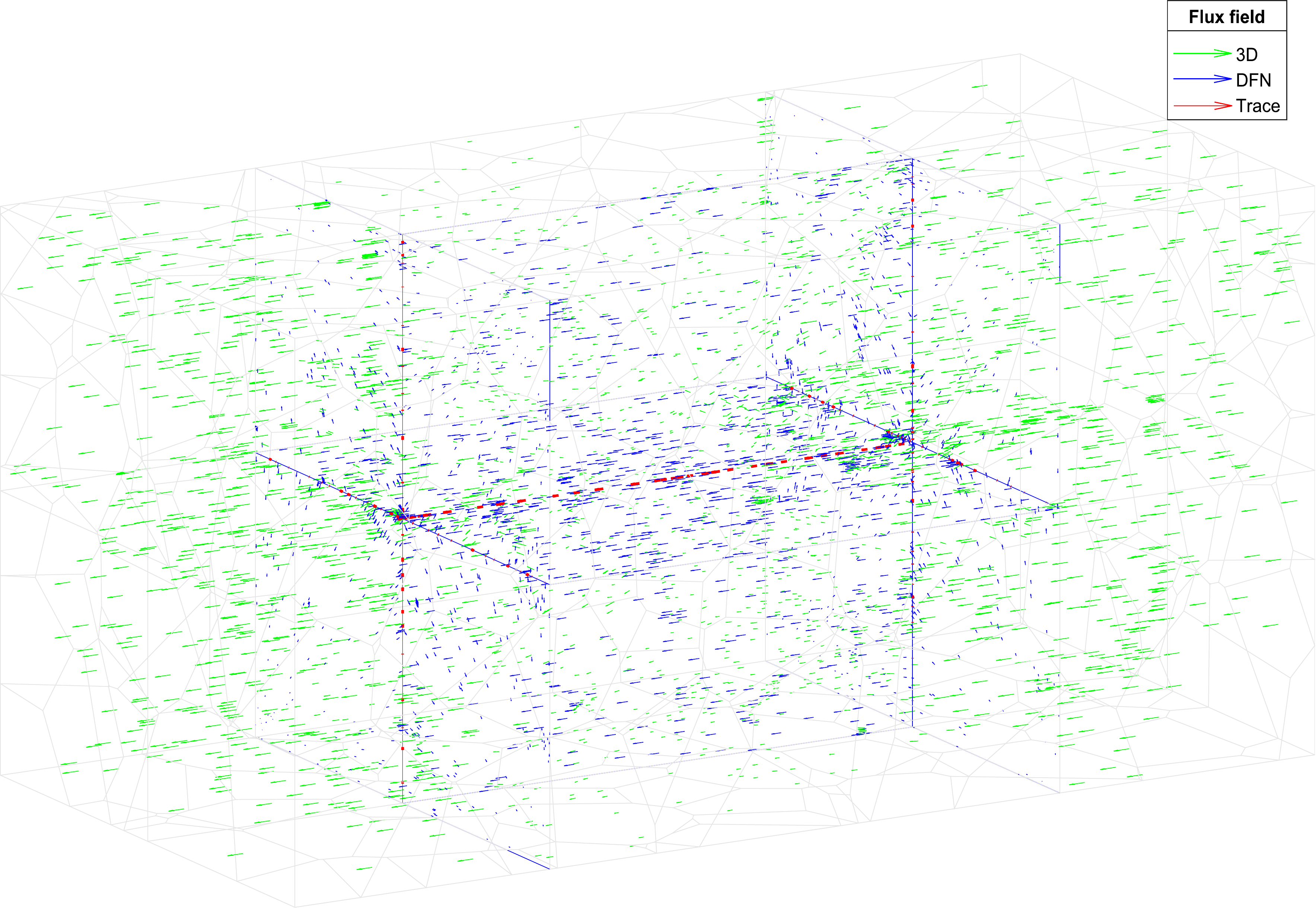}
    \caption{Finite interdimensional normal permeability}
    \label{Problem2_V:finite}
  \end{subfigure}\\
	\begin{subfigure}{.8\linewidth}
    \includegraphics[width=\linewidth]
    {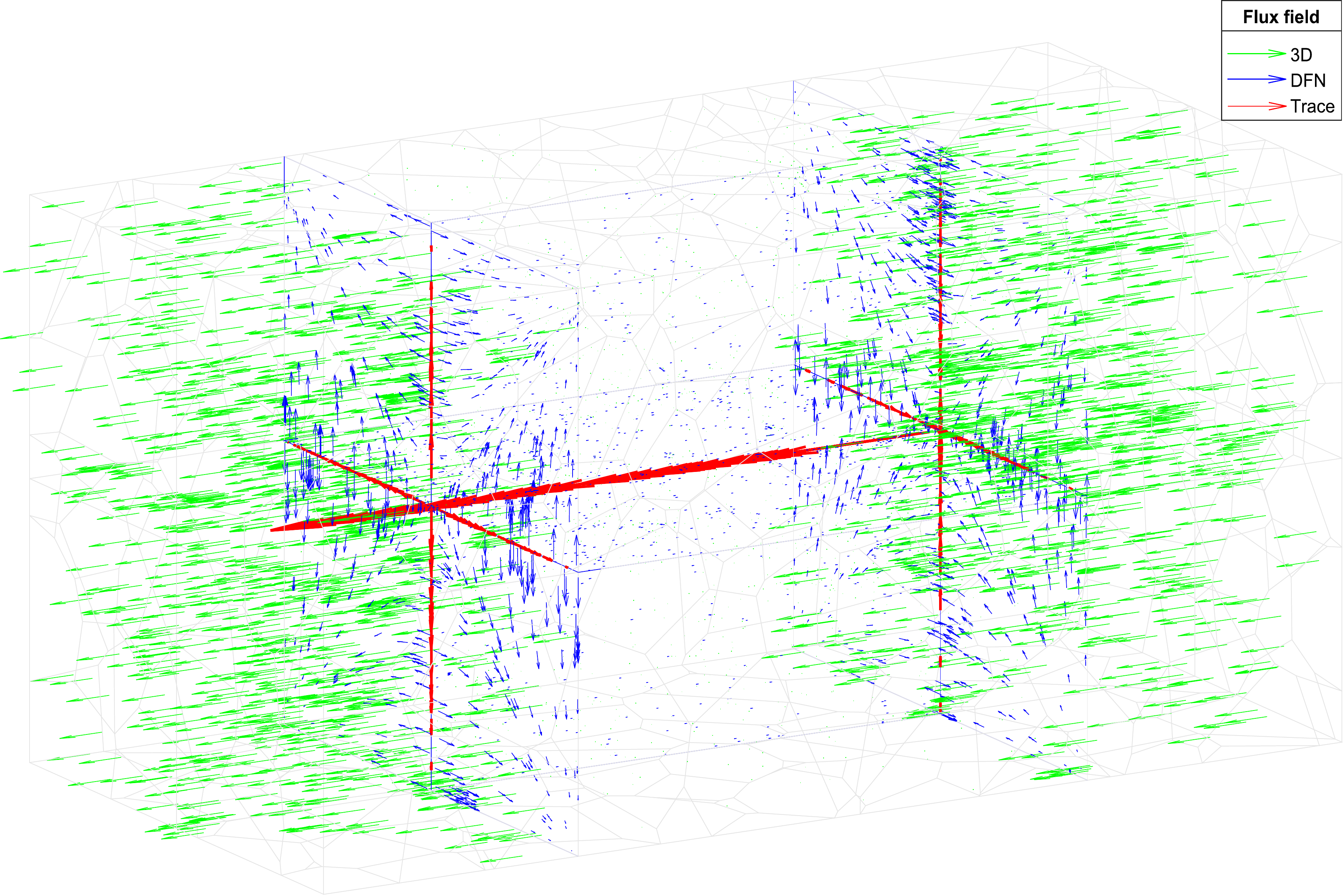}
    \caption{Global pressure continuity}
    \label{Problem2_V:cont}%
  \end{subfigure}
  \caption{Problem 2: Global discrete solutions for the (normalized)
    velocity field.}
  \label{Problem2_V}
\end{figure}

\begin{figure}[!ht]
  \centering
  \begin{subfigure}{.6\linewidth}\centering
    \includegraphics[width=\linewidth] {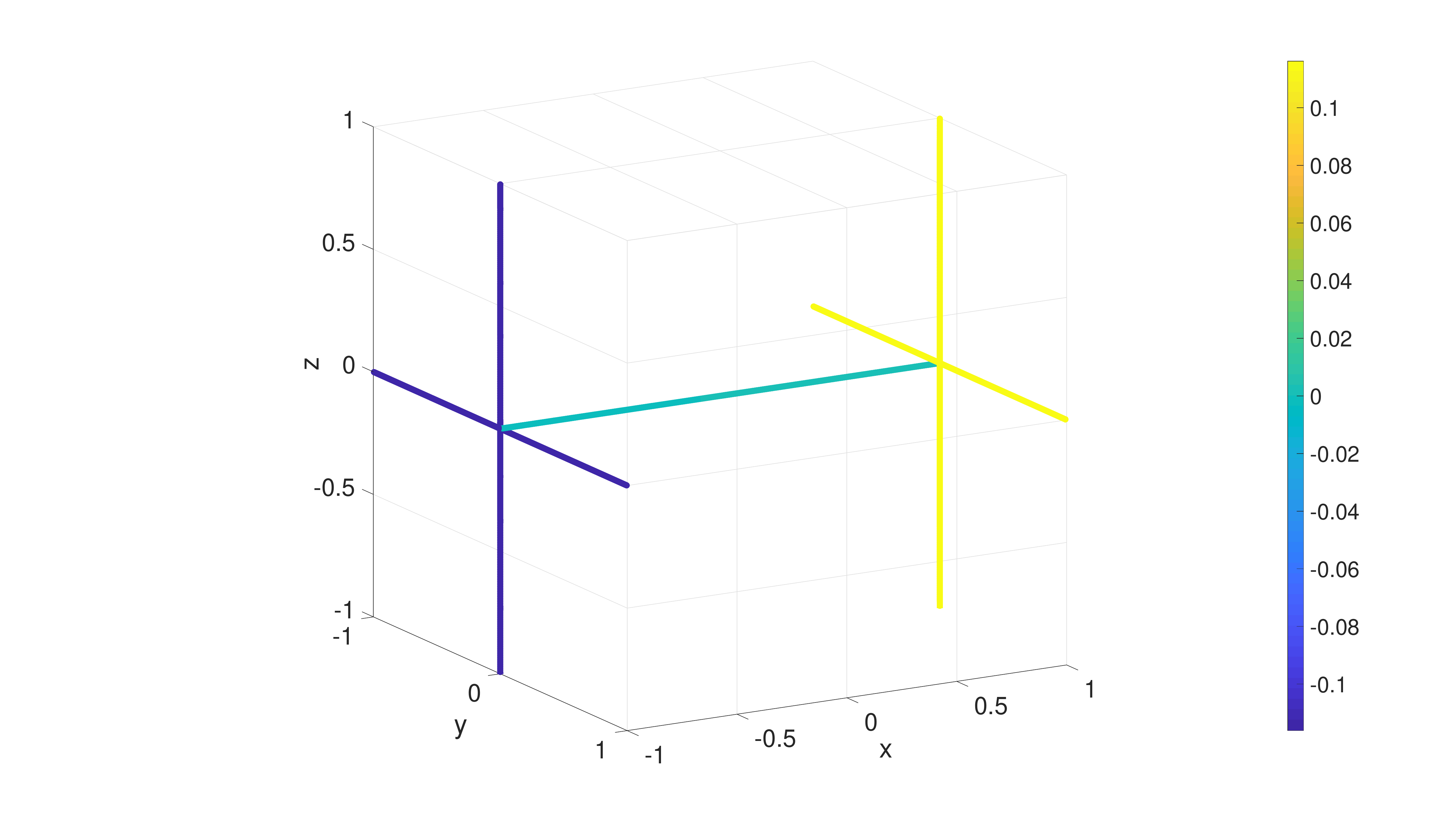}
    \caption{Pressure head on trace network}
    \label{Problem2_solutions2:trace}
  \end{subfigure}
  \\
  \begin{subfigure}{.49\linewidth}\centering
    \includegraphics[height=.5\linewidth] {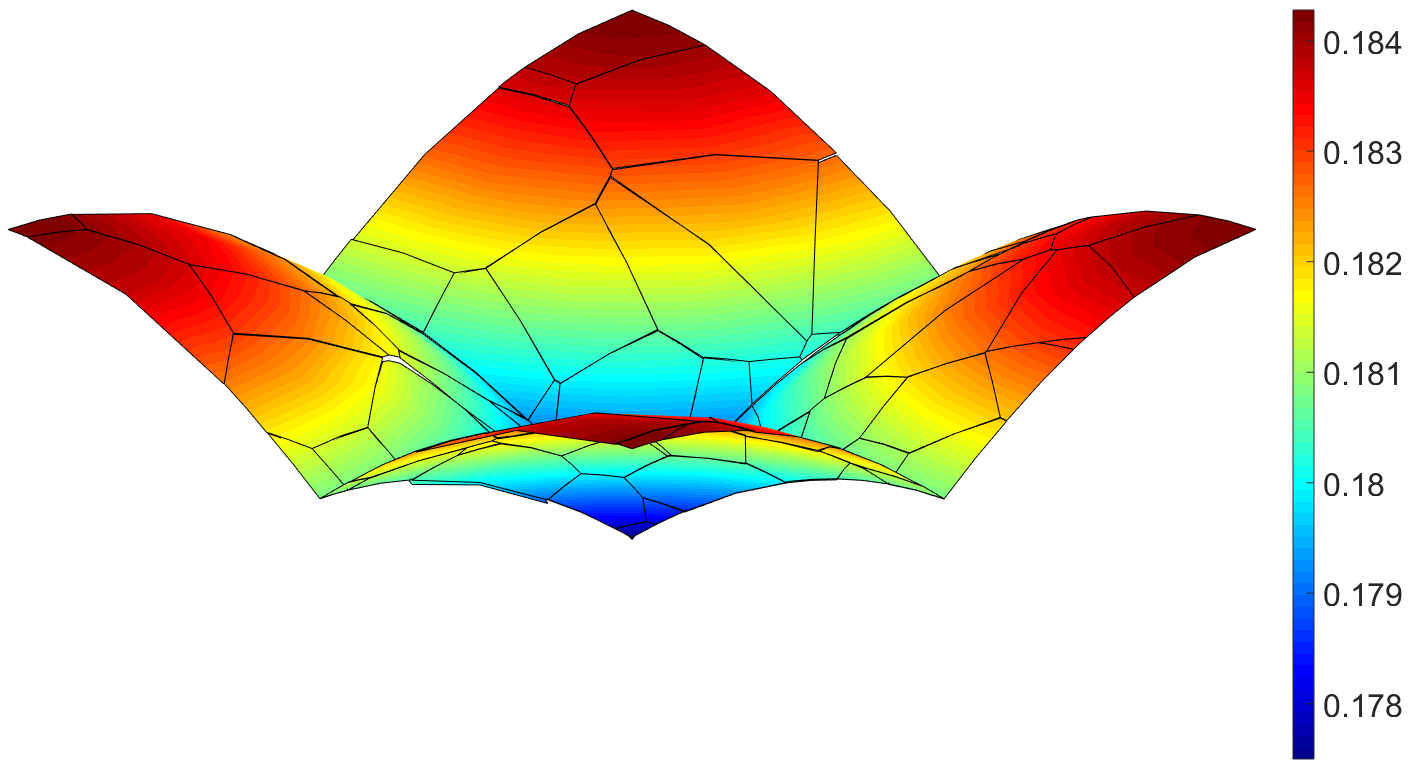}
    \caption{Pressure head on $\Omega^2_1$}
    \label{Problem2_solutions2:F1P} %
  \end{subfigure}
  \begin{subfigure}{.49\linewidth}\centering
    \includegraphics[height=.5\linewidth] {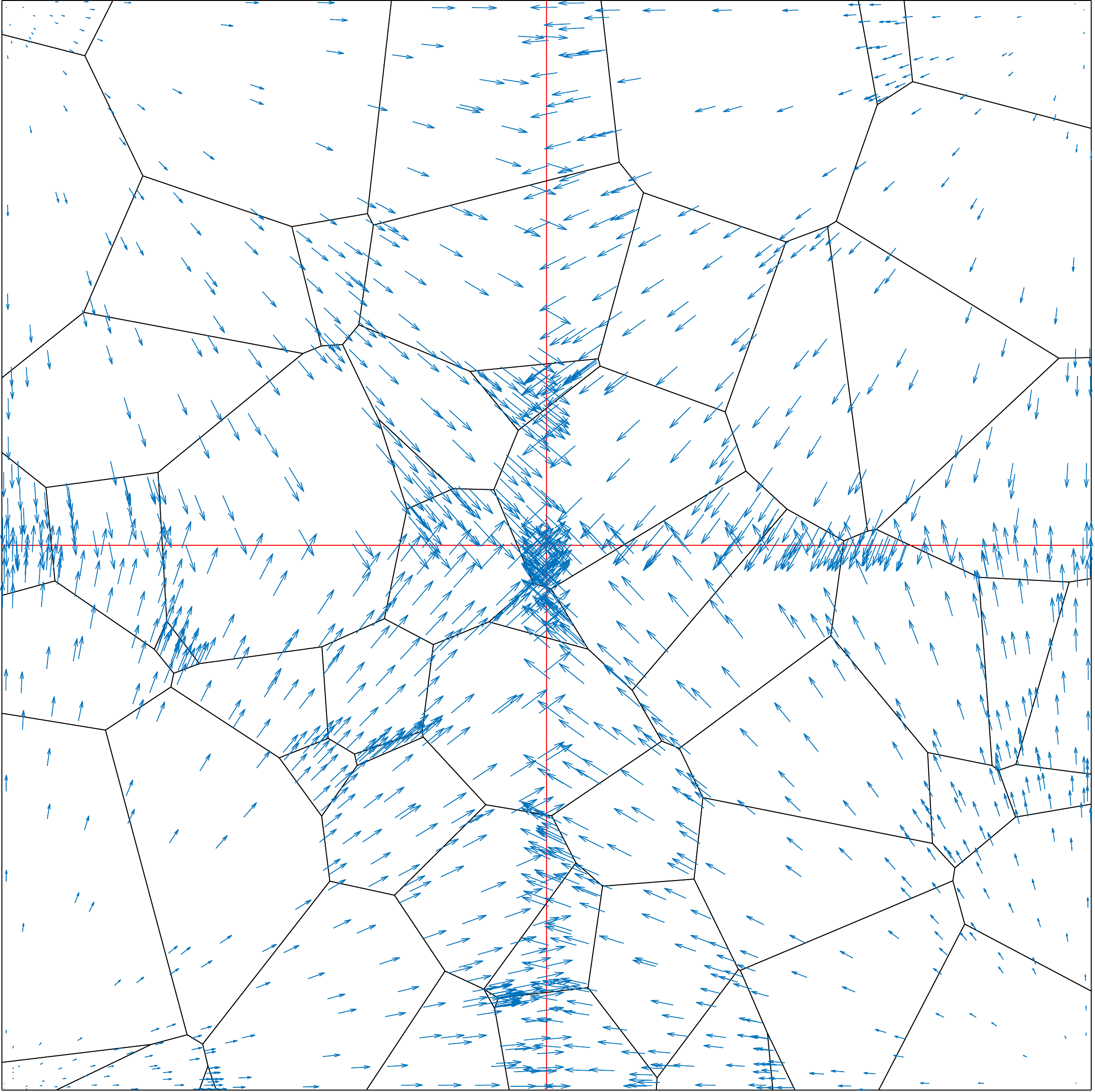}
    \caption{Velocity field head on $\Omega^2_1$}
    \label{Problem2_solutions2:F1V} %
  \end{subfigure}%
\caption{Problem 2: Discrete solution}
\label{Problem2_solutions2} 
\end{figure}

The obtained results present clearly noticeable jumps in pressure (see
Figure \ref{Problem2:Pfield}), as seen across $\Omega^3$ and
$\Omega^2_1$ and $\Omega^2_4$. The effect of the finite
interdimensional permeability in the velocity field is clear: it
penalizes flux exchange since any non-zero flux exchange results in a
drop of pressure inversely proportional to the normal
permeability. For this reason, the velocity field in the mid part of
the rock matrix is not entirely negligible (Figure
\ref{Problem2_V:finite}) despite the much higher permeability of the
fractures and specially the traces, which represent a much more
favourable path between entry and exit point. A second computation is
provided as a comparison in which global continuity is imposed (\ie,
infinite normal permeability) while keeping all other parameters the
same. In this case, the aforementioned effect does not occur and the
flow path is clearly dominated by the lower dimension entities
(specially $\Omega^1_3$), resulting in an almost vanishing flux in the
mid section of the rock matrix and on fractures $\Omega^2_2$ and
$\Omega^2_3$ (Figure \ref{Problem2_V:cont}). Another consequence of
global continuity is a difference in the total flux through the
system. Since Dirichlet boundary conditions are imposed, incoming flux
is dependent on the permeability parameters of the problem and the
computed incoming fluxes are $2.91$ and $7.99$ for the finite
permeability and global continuity situations respectively. As
expected, global pressure continuity across dimensions provides a more
favourable flow path, thus the global domain is more permeable
resulting in larger flux values. Pressure head in the 1D trace network
is given in Figure \ref{Problem2_solutions2:trace} showing pressure
jumps across trace intersections. Finally, discrete solutions for
pressure head and velocity field on fracture $\Omega^2_1$ are shown in
Figures \ref{Problem2_solutions2:F1P} and
\ref{Problem2_solutions2:F1V}.

\newpage
\subsection{Problem 3: Realistic embedded 16 fracture DFN}
\label{16FractureDFN}
\begin{figure}[!h]
	\centering
	\begin{subfigure}{.45\linewidth}
    \includegraphics[width=\linewidth] {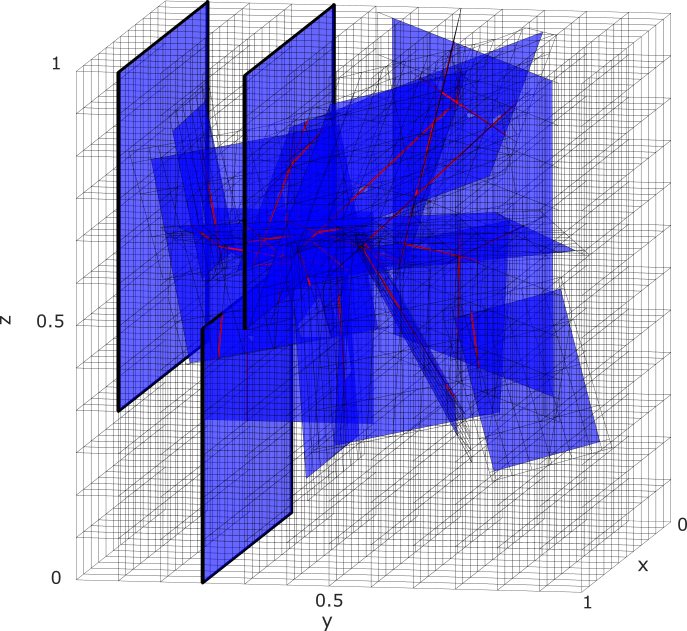}
    \caption{Base mesh and embedded DFN}
    \label{Problem3B_DFN}
	\end{subfigure}\hfill
	\begin{subfigure}{.4\linewidth}
    \includegraphics[width=\linewidth] {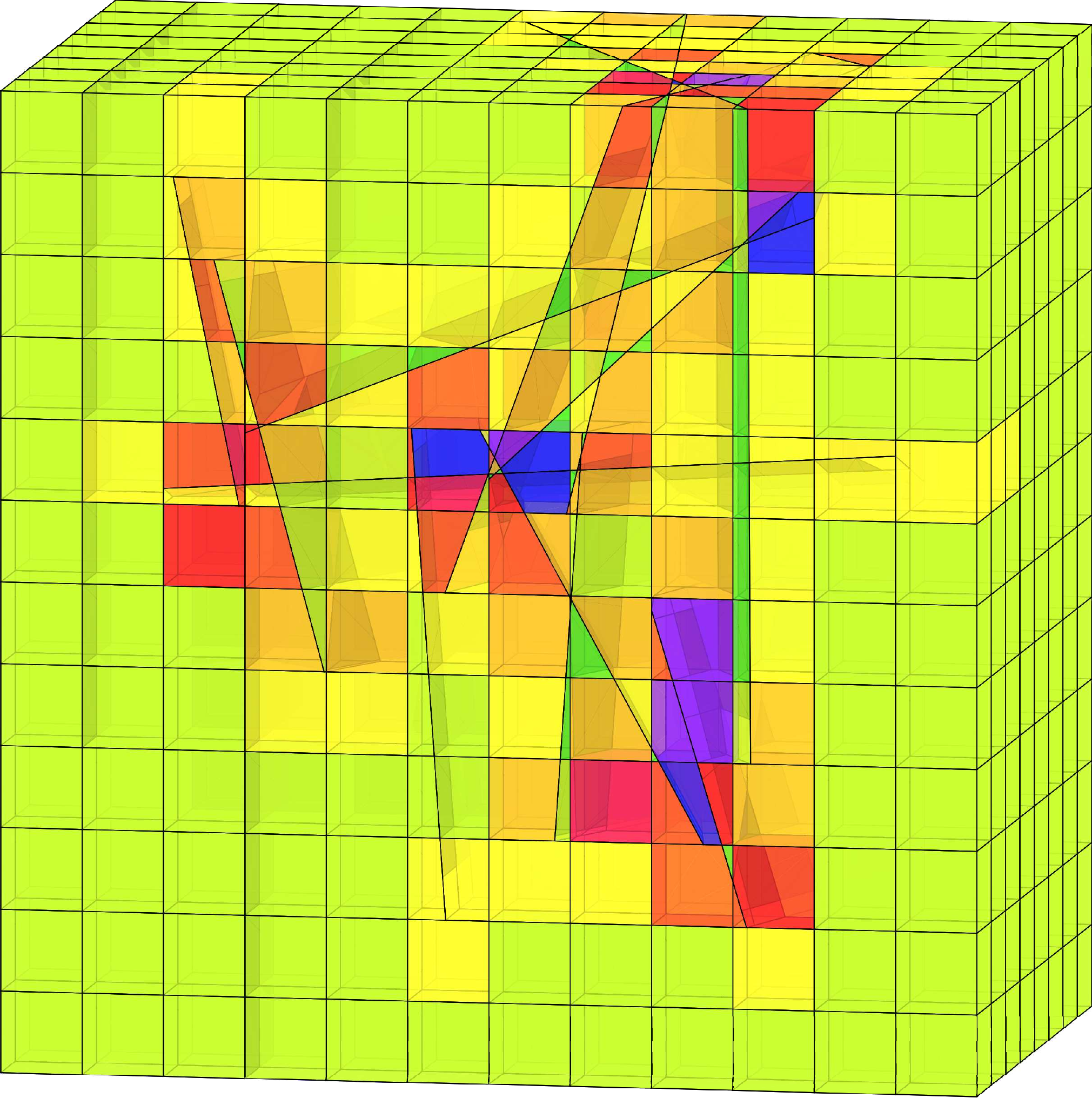}
    \caption{VEM mesh (sliced for visualization)}
    \label{Problem3B_mesh} %
	\end{subfigure}
	\caption{Problem 3: Geometry and discretization}
\end{figure}

In this problem a 16 fracture DFN containing 38 traces and 15
intersections is embedded into a matrix of domain is
$\Omega^3 = \left[0,1\right]^3$ as shown in Figure
\ref{Problem3B_DFN}, where fractures are shaded and traces are
depicted in red. Boundary conditions for the matrix consist of a
constant unitary incoming flux on face $\lbrace y = 1 \rbrace$,
homogeneous Dirichlet on face $\lbrace y = 0 \rbrace$ and homogeneous
Neumann (no-flux condition) on all remaining faces. For the DFN and
the trace network, no-flux condition is imposed on all
boundaries. There are neither source nor sink terms, and the given
parameters are as follows:
\begin{align*}
  &\mathfrak{a}^3 = 1, && \eta^2_l = 1 \ (l=1,...,14), \quad \eta^2_l = 0.01 \ (l=14,15,16)\\
  &\mathfrak{a}^2_l = 10^2 \ (l=1,...,16), && \eta^1_l \rightarrow \infty\ (l=1,...,39),\\
  &\mathfrak{a}^1_l = 10^4 \ (l=1,...,39), && \eta^0_l \rightarrow \infty \ (l=1,...,15).
\end{align*} 
Note that normal permeability between matrix and fractures is finite
and therefore pressure head will not be globally continuous. In
particular, fractures $\Omega^2_{14}$, $\Omega^2_{15}$ and
$\Omega^2_{16}$ (whose boundaries are depicted in a darker shade) have
very low normal permeability and effectively act as flow barriers.

The starting point for the discretization is a polyhedral mesh
comprised of $12^3 = 1728$ cubic elements, which, after successive
fracture additions, results in the final mesh (Figure
\ref{Problem3B_mesh}) of $2734$ polyhedral elements. Even though the
number of fractures may seem low to discretize a network, which can
amount to thousands of fractures, a small number of fractures can
provide a better appreciation of the numerical results. Furthermore,
even in this relatively contained example, the demand for
computational power is high. Besides that, the generation of the final
globally conforming mesh that takes the fractures into account is
challenging as well, resulting in a final mesh that contains a
plethora of very badly shaped elements with many undesirable features
such as: small angles, tiny faces, short edges, large discrepancies in
size between adjacent elements, collapsing nodes, etc (See mesh
details on Figure \ref{Problem3B_mesh}). Nevertheless, it would be
undoubtedly more computationally expensive to solve this problem with
a Finite Element mesh made up of the usual shapes (tetrahedra,
hexahedra, pyramids and wedges) that is globally conforming from the
start. This is due to the fact that a mesh for that situation would
require a very small element size to be of acceptable quality as well
as some type of iterative procedure to obtain matching discretizations
between 3D faces and 2D fractures. Furthermore, VEM has been shown to
be very robust to mesh distortion and the resulting mesh can be
readily handled by the method.

The problem was solved using RT0-VEM, BDM1-VEM, RT1-VEM, BDM2-VEM and
RT2-VEM, whose discretization statistics are provided in Table
\ref{Problem3B_Table}. Note that since 3D RTk-VEM and BDMk-VEM
elements of the same order have an equal number of face DOFs for the
flux variable, they can both be used in conjunction with 2D RTk-VEM
elements for the fractures. The saving in DOFs for BDM elements is due
to a reduction of internal flux DOFs and a lower order approximation
of the 3D pressure field.

\begin{table}[]
	\centering
	\begin{tabular}[t]{ccc|ccc|cc|cc|cc|c}
		\toprule
		\multicolumn{3}{c|}{Element Type} & \multicolumn{3}{c|}{\#Elements} & \multicolumn{2}{c|}{\#3D DOFs} & \multicolumn{2}{c|}{\#2D DOFs} & \multicolumn{2}{c|}{\#1D DOFs}  & \multicolumn{1}{c}{\#0D DOFs}\\
		\midrule
                                      & &     & 3D & 2D & 1D  & Flux & Pressure & Flux & Pressure & Flux & Pressure & Pressure \\
		\midrule
		RT0&RT0&$\mathbb{P}_{1}/\mathbb{P}_{0}$& 2734& 1845 &289& 11584 & 2734 & 4707 & 1845 & 371 & 289 & 15\\
		BDM1&RT1&$\mathbb{P}_{2}/\mathbb{P}_{1}$& 2734& 1845 &289& 42957 & 2734 & 14972 & 5538 & 659 & 578 & 15\\
		RT1&RT1&$\mathbb{P}_{2}/\mathbb{P}_{1}$& 2734& 1845 &289& 51159 & 10936 & 14972 & 5538 & 659 & 578 & 15\\
		BDM2&RT2&$\mathbb{P}_{3}/\mathbb{P}_{2}$& 2734& 1845 &289& 107786 & 10936 & 28919 & 11076 & 951 & 867 & 15\\
		RT2&RT2&$\mathbb{P}_{3}/\mathbb{P}_{2}$& 2734& 1845 &289& 124190 & 27340 & 28919 & 11076 & 951 & 867 & 15\\
		\bottomrule
	\end{tabular}
	\caption{Problem 3: discretization data}
	\label{Problem3B_Table}
\end{table}

\begin{figure}[!h]
	\centering
	\begin{subfigure}{\linewidth}
    \centering
    \includegraphics[height=.38\linewidth] {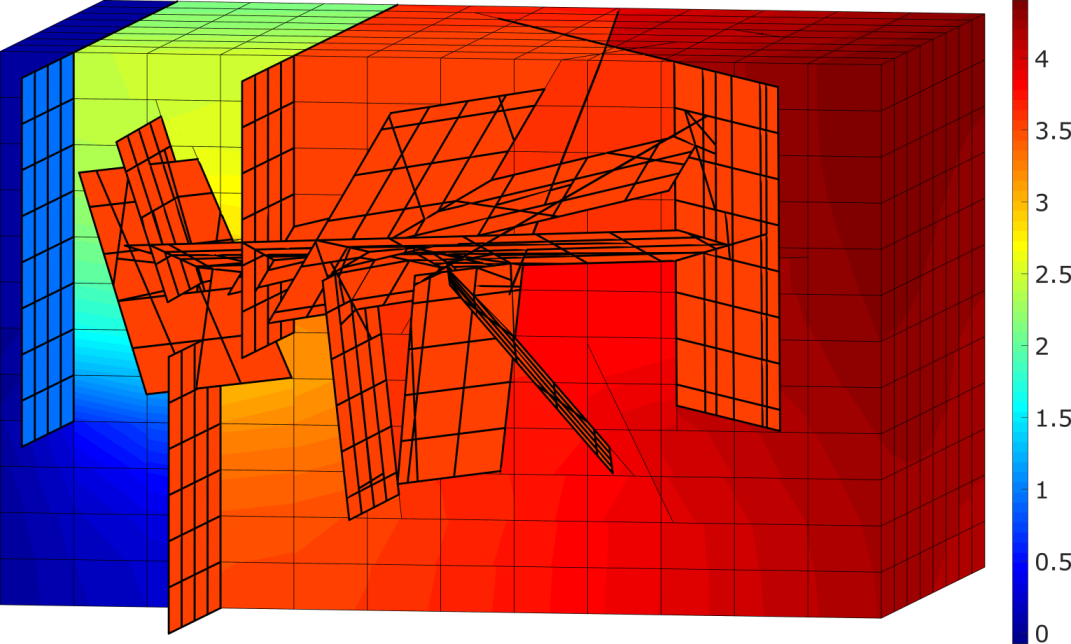}
    \caption{Global pressure field}
	\end{subfigure}
	\\
	\begin{subfigure}{\linewidth}
    \centering
    \includegraphics[height=.5\linewidth]
    {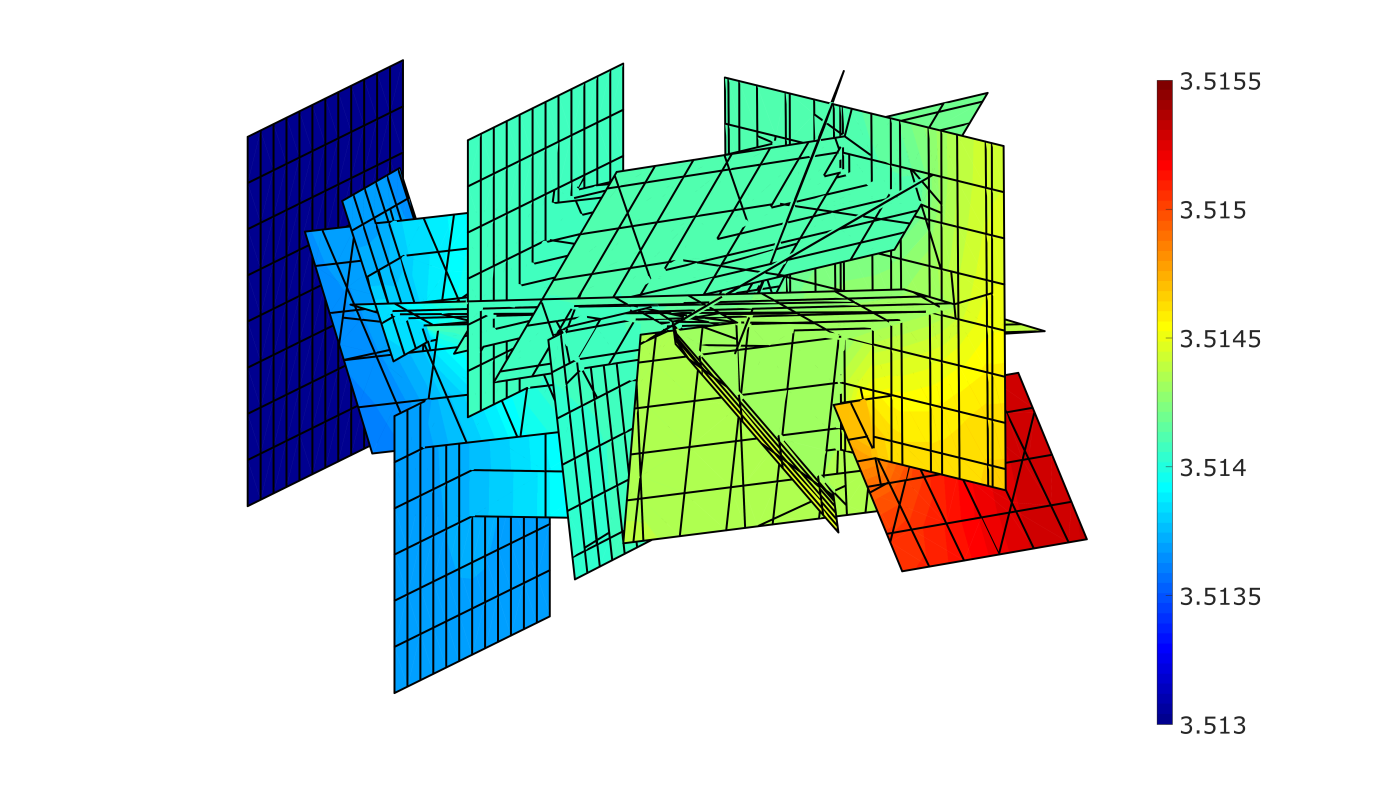}
    \caption{Pressure field on DFN and traces}
	\end{subfigure}
	\caption{Problem 3: Discrete pressure solutions}
	\label{Problem3B_P}
\end{figure}

\begin{figure}[!h]
	\centering \includegraphics[width=.8\linewidth]{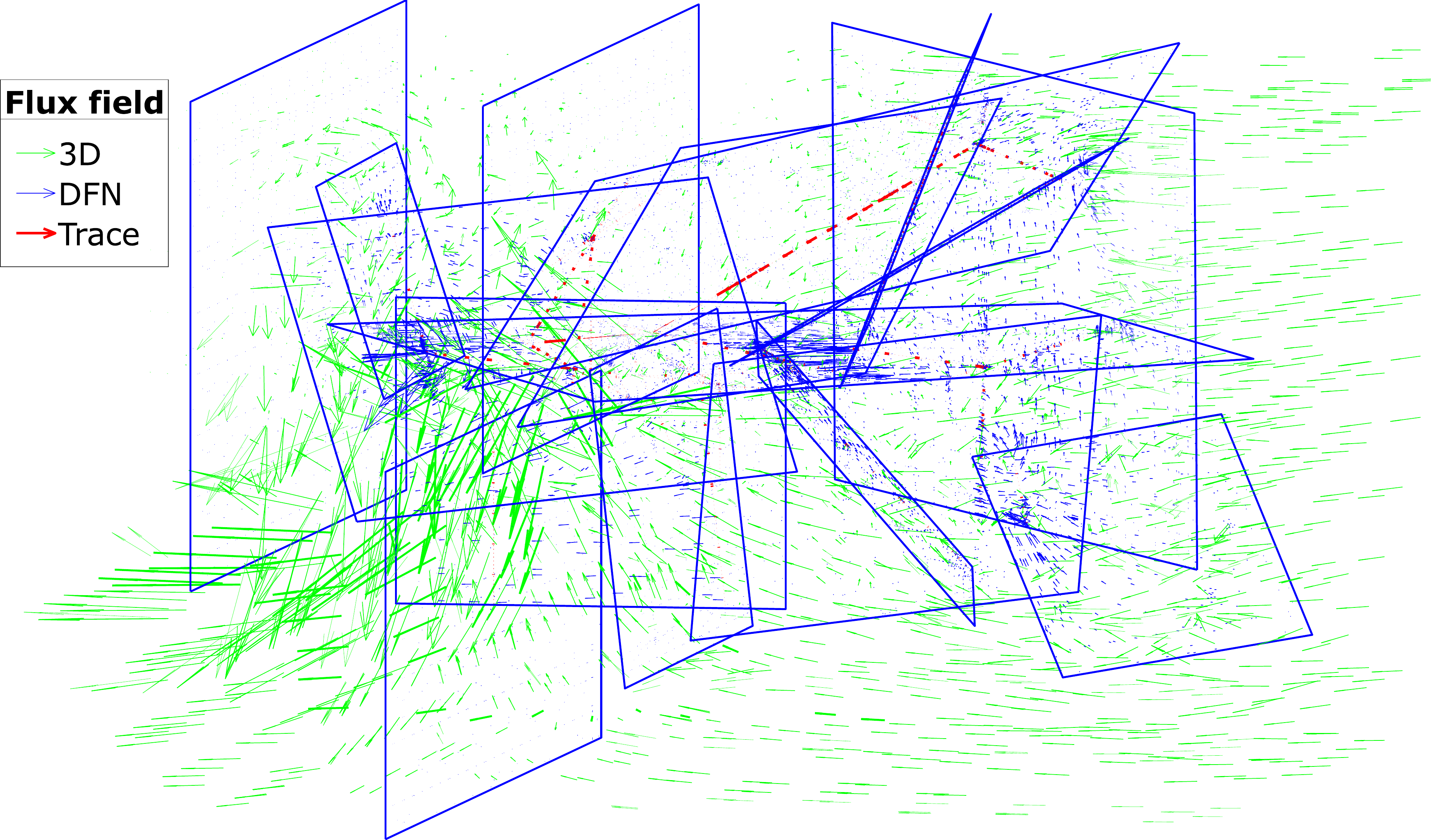}%
	\caption{Problem 3: Discrete velocity field solution}
	\label{Problem3B_V}
\end{figure}

\begin{figure}[htbp]
	\centering
	\begin{subfigure}{.4\linewidth}
    \centering
    \includegraphics[width=\linewidth] {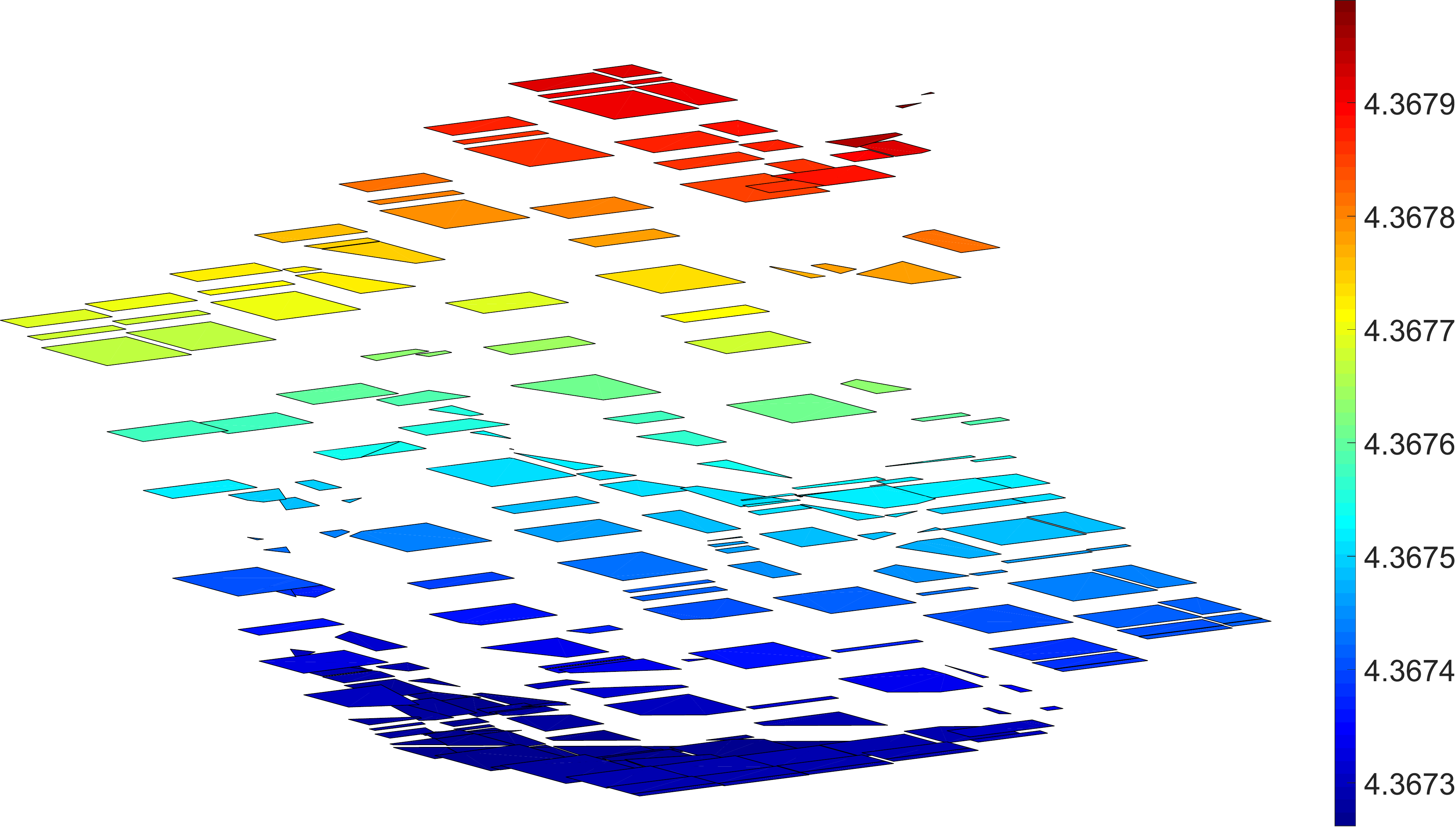}
    \caption{RT0-VEM pressure head}
    \label{Problem3B_F8_O0}
	\end{subfigure}
	\begin{subfigure}{.4\linewidth}
    \centering
    \includegraphics[width=\linewidth] {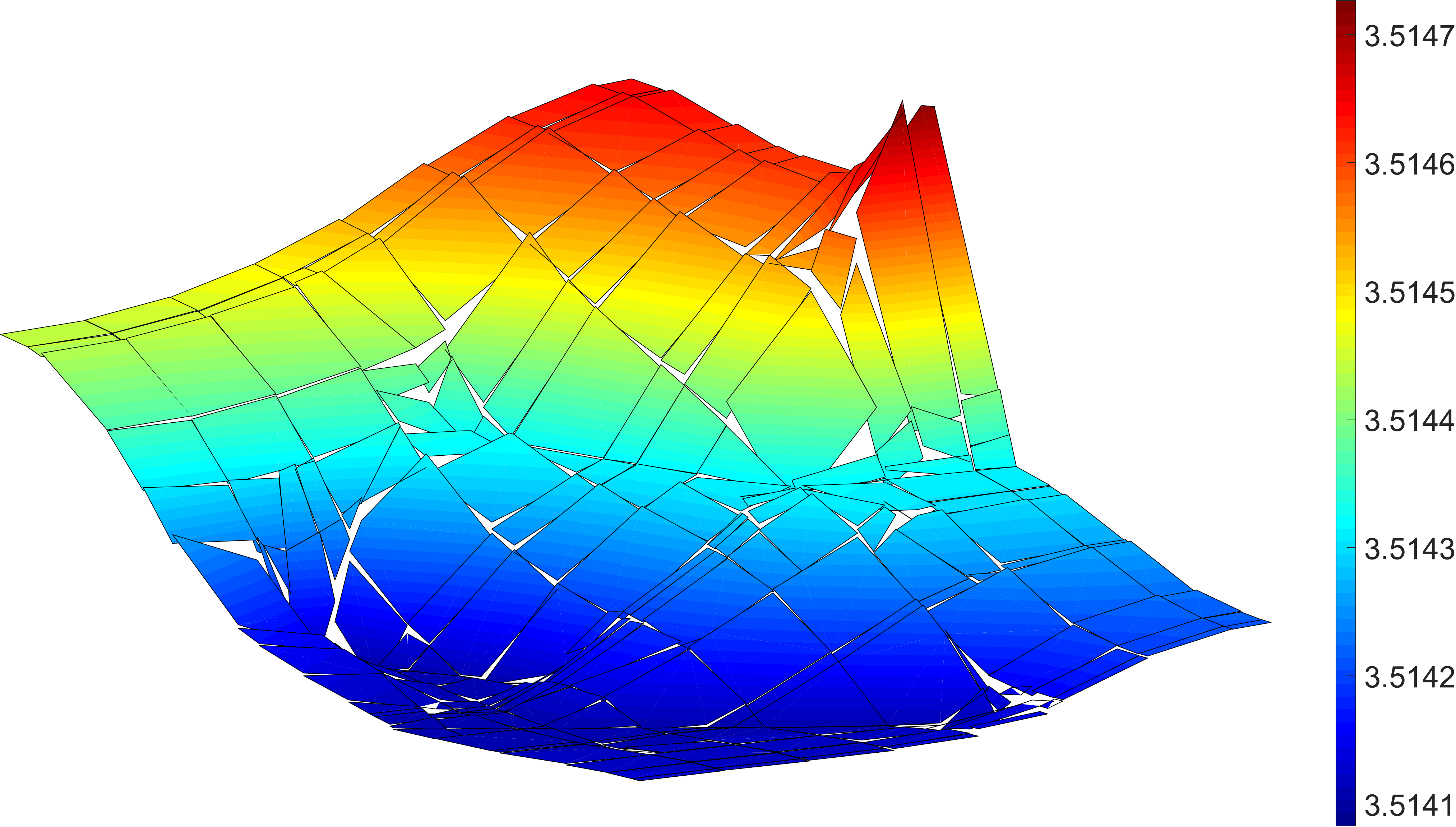}
    \caption{RT1-VEM pressure head}
    \label{Problem3B_F8_O1}%
	\end{subfigure}
	\begin{subfigure}{.4\linewidth}
    \centering
    \includegraphics[width=\linewidth] {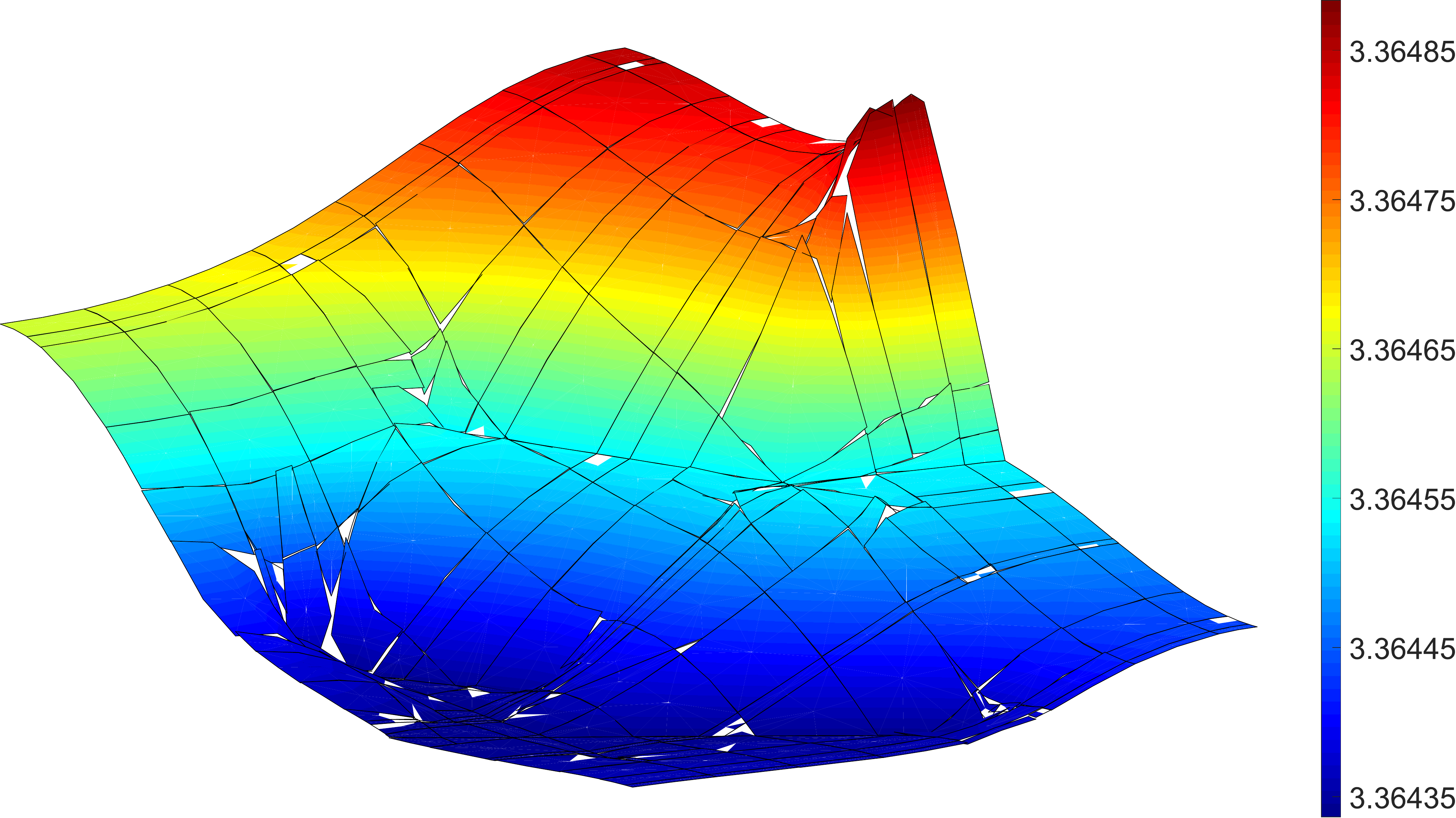}
    \caption{RT2-VEM pressure head}
    \label{Problem3B_F8_O2}%
	\end{subfigure}
	\\
	\begin{subfigure}{.6\linewidth}
    \centering
    \includegraphics[width=\linewidth] {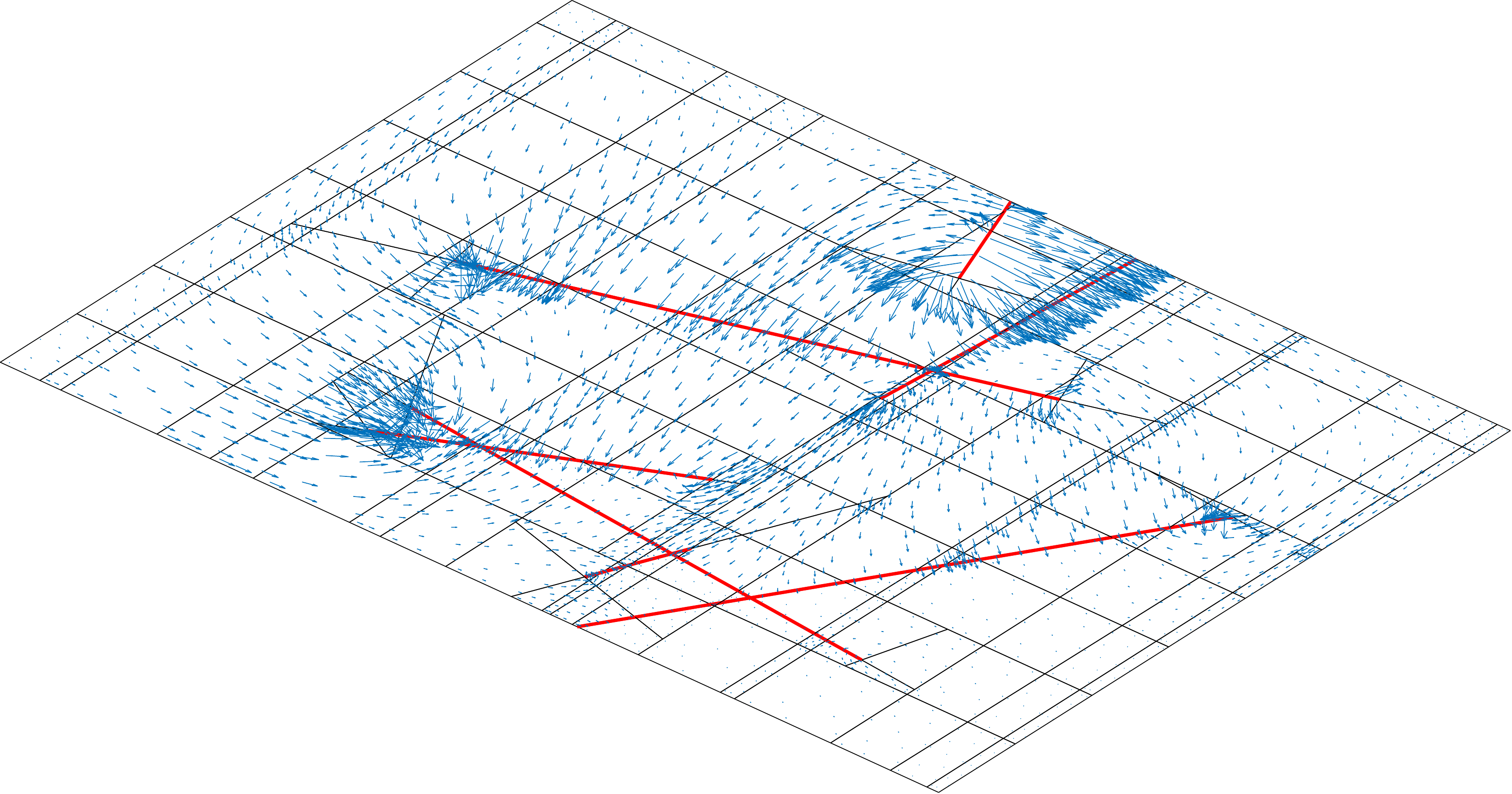}
    \caption{Velocity field (traces in red)}
    \label{Problem3B_F3:Pfield} %
  \end{subfigure}
	\caption{Problem 3: Discrete solutions on $\Omega^2_8$}
	\label{Problem3B_F3}
\end{figure}

In the next Figures the results of the analysis are presented. The
global pressure and the velocity field are shown in Figures
\ref{Problem3B_P} and \ref{Problem3B_V}, shown for the case of RT1-VEM
elements, chosen for their linear approximation of both the flux and
pressure variable. For the pressure, it can be clearly seen that the
field is not globally continuous across fractures, most noticeably
across those fractures with very low normal permeability, and that the
small drop in pressure across the DFN is low due to the relatively
high tangential permeability. The velocity field showcases the
preferred flow paths in the multidimensional domain and the flux
exchanges that take place. Despite not being a very dense, the trace
network provides an advantageous flow path due to the high tangential
permeability value and to the pressure continuity between fractures
and traces as well as on trace intersections. The pressure and
(normalized) velocity field for a particular fracture ($\Omega^2_3$)
presented in \ref{Problem3B_F3} for 2D RT0-, RT1- and RT2-VEM show a
very strong qualitative agreement and a convergence towards and almost
continuous local pressure field. Pressure values on an individual
fracture are of course dependent on the global solution in the
complete domain.

Some conclusions can be drawn from this problem. Firstly, the solution
is highly dependent on the fracture distribution as expected. In fact,
flux is greatest on fractures that provide the least effort path
between boundaries with prescribed boundary conditions. Secondly,
relative permeability between rock matrix and fractures has
considerable influence on the result that can be clearly appreciated
by the fact that only a very small drop in pressure is experienced
across the whole DFN. On the other hand, taking into account the rock
matrix in the model makes it significantly more expensive to solve,
and should be avoided whenever the rock permeability is low enough
such that the problem can be considered as mainly
DFN-dominated. Thanks to global conservation of flux, results are very
reliable for the flux variable even for low order elements. If 1D
tangential permeability is high relative to the other domains, flux
across the trace network is not negligible and should not be
disregarded. A denser network of traces would further increase the
influence of considering 1D elements in the solution. Finally, by
inspecting the meshes used for computing the discrete solutions
(obtained without any mesh improvement technique whatsoever) it can be
seen that mixed VEM show great robustness and are affected very little
by mesh quality.

\section{Conclusions} 
\label{Conclusions}
A mixed Virtual Element formulation for multimensional problems for solving flow problems was explored in this work. A summary of the methodology as well as the main ideas of the method were included and details of its implementation were given a thorough treatment for an arbitrary discretization order. The crucial feature is the use of a polyhedral discretization with no initial requirement of mesh conformity. In fact, the rock matrix is given any mesh, on top of which the fractures that make up the DFN are added. Later, these fracture are incorporated to the global mesh by introducing cuts on the 3D elements. This would require some kind of remeshing-procedure for standard Finite Elements but that is not the case with VEM, than can handle the modified elements without any particular distinction over the original elements. Ultimately, a globally conforming polyhedral VEM mesh is obtained where faces, edges and vertices are automatically shared between inter-dimensional entities. Due to VEM's robustness to mesh distortion, the numerical accuracy is not affected by the irregularity of the mesh.

The complete formulation for problems involving the combination of multidimensional entities was provided, where  2D elements representing embedded fractures are included in a 3D solid matrix. In addition, 1D elements representing traces given by fracture intersections were also included so as to take tangential trace flow into account. These traces constitute a trace network that requires flux balance and adequate pressure conditions on the points defined by traces intersection. Each interaction between entities of different dimensions results in flux exchanges and pressure jumps. 
One key aspect to highlight is that whereas the primal formulation guarantees continuity of pressure head (or pressure head \textit{jump} in the case of non infinite permeability), the mixed formulation has the advantage of being completely mass-conservative by definition of the discrete spaces themselves. This is a more desirable characteristic in flux computations in many cases where the computed flow field is used as the underlying convection field to serve as input for another model (\eg, in transport problems).

Several numerical experiments were performed to assess the validity of the methodology. Starting with a problem with a known polynomial solution, it was shown that the exact solution is recovered within machine accuracy over meshes made up of arbitrary polyhedra. Afterwards, a multidimensional problem was solved to highlight the interactions between inter-dimensional entities and to evaluate the effects of finite inter-dimensional permeability over the flow field. Finally, a small yet realistic DFN was embedded in a 3D matrix for which the flow and pressure fields were determined and comparisons of discrete solutions were made between different interpolation orders.

It is straightforward to generalize the approach to second order elliptic problems. Time-dependent phenomena could be tackled as well and would follow standard procedure once the mass and stiffness matrix for the VEM discretization are obtained. From a computational point of view, the method shows potential for parallelization, since each dimension of the problem can be computed independently (even in parallel themselves) and the coupling between dimensions can be added once the respective stiffness matrices have been computed. However, 3D computations can be orders of magnitude slower than their 2- and 1-dimensional counterparts, mainly due to the much steeper increase in 3D DOFs for mesh refinement and increasing order of interpolation.

In conclusion, the arbitrary order mixed formulation Virtual Element Method presented here displays clear advantages in the study of geometrically complex hybrid dimensional flow problem mainly due to simpler meshing, easier implementation, strongly imposed mass conservation and accurate flux results.

\section*{Acknowledgments}
A.B. and S.S. acknowledge the financial support of MIUR through the
PRIN grant n. 201744KLJL and through the project ``Dipartimenti di
Eccellenza 2018-2022'' (CUP E11G18000350001), and by INdAM-GNCS.


\bibliographystyle{siam}
\bibliography{Bibliography}

\begin{thebibliography}{10}

\bibitem{Ahmed2015}
{\sc R.~Ahmed, M.~Edwards, S.~Lamine, B.~Huisman, and M.~Pal}, {\em
  Control-volume distributed multi-point flux approximation coupled with a
  lower-dimensional fracture model}, Journal of Computational Physics, 284
  (2015), pp.~462 -- 489.

\bibitem{AlHinai2015}
{\sc O.~Al-Hinai, S.~Srinivasan, and M.~F. Wheeler}, {\em Domain decomposition
  for flow in porous media with fractures}, in SPE Reservoir Simulation
  Symposium 23-25 February 2013, Houston, Texas, USA, Society of Petroleum
  Engineers, 2015.

\bibitem{Wheeler2015}
{\sc O.~Al-Hinai, S.~Srinivasan, and M.~F. Wheeler}, {\em Mimetic finite
  differences for flow in fractures from microseismic data}, in SPE Reservoir
  Simulation Symposium, Society of Petroleum Engineers, 2015.

\bibitem{Angot2009}
{\sc {Angot, Philippe}, {Boyer, Franck}, and {Hubert, Florence}}, {\em
  Asymptotic and numerical modelling of flows in fractured porous media},
  ESAIM: M2AN, 43 (2009), pp.~239--275.

\bibitem{Antonietti2016}
{\sc P.~F. Antonietti, L.~Formaggia, A.~Scotti, M.~Verani, and N.~Verzott},
  {\em Mimetic finite difference approximation of flows in fractured porous
  media}, ESAIM: Mathematical Modelling and Numerical Analysis, 50 (2016),
  pp.~809--832.

\bibitem{VEMmixedStress}
{\sc E.~Artioli, S.~de~Miranda, C.~Lovadina, and L.~Patruno}, {\em A
  stress/displacement {V}irtual {E}lement method for plane elasticity
  problems}, Computer Methods in Applied Mechanics and Engineering, 325 (2017),
  pp.~155 -- 174.

\bibitem{VEMbasic}
{\sc L.~Beir\~{a}o~da Veiga, F.~Brezzi, A.~Cangiani, G.~Manzini, L.~D. Marini,
  and A.~Russo}, {\em Basic principles of {V}irtual {E}lement methods}, Math.
  Models and Methods in Applied Sciences,  (2013).

\bibitem{Beirao2016}
{\sc L.~Beir\~{a}o~da Veiga, F.~Brezzi, L.~D. Marini, and A.~Russo}, {\em Mixed
  virtual element methods for general second order elliptic problems on
  polygonal meshes}, ESAIM: Mathematical Modelling and Numerical Analysis, 50
  (2016), pp.~727--747.

\bibitem{divcurlVEM}
{\sc L.~Beir{\~a}o~da Veiga, F.~Brezzi, L.~D. Marini, and A.~Russo}, {\em
  H(div) and {H}(curl)-conforming virtual element methods}, Numer. Math., 133
  (2016), pp.~303--332.

\bibitem{BenedettoThesis}
{\sc M.~F. Benedetto}, {\em The mixed virtual element formulation for three
  dimensional problems with embedded planes}, Master's Thesis,  (2018).

\bibitem{BBBChapter}
{\sc M.~F. Benedetto, S.~Berrone, and A.~Borio}, {\em The {V}irtual {E}lement
  {M}ethod for underground flow simulations in fractured media}, in Advances in
  Discretization Methods, vol.~12 of SEMA SIMAI Springer Series, Springer
  International Publishing, Switzerland, 2016, pp.~167--186.

\bibitem{BBPS}
{\sc M.~F. Benedetto, S.~Berrone, S.~Pieraccini, and S.~Scial\`o}, {\em The
  virtual element method for discrete fracture network simulations}, Comput.
  Methods Appl. Mech. Engrg., 280 (2014), pp.~135 -- 156.

\bibitem{DFNvem1}
{\sc M.~F. Benedetto, S.~Berrone, and S.~Scial\`o}, {\em A globally conforming
  method for solving flow in discrete fracture networks using the virtual
  element method}, Finite Elements in Analysis and Design, 109 (2016), pp.~23
  -- 36.

\bibitem{VEMmixedDFN}
{\sc M.~F. Benedetto, A.~Borio, and S.~Scialò}, {\em Mixed virtual elements
  for discrete fracture network simulations}, Finite Elements in Analysis and
  Design, 134 (2017), pp.~55 -- 67.

\bibitem{BBorth}
{\sc S.~Berrone and A.~Borio}, {\em Orthogonal polynomials in badly shaped
  polygonal elements for the {V}irtual {E}lement {M}ethod}, Finite Elements in
  Analysis \& Design, 129 (2017), pp.~14--31.

\bibitem{BERRONE201714}
{\sc S.~Berrone and A.~Borio}, {\em Orthogonal polynomials in badly shaped
  polygonal elements for the virtual element method}, Finite Elements in
  Analysis and Design, 129 (2017), pp.~14 -- 31.

\bibitem{Berrone2018}
{\sc S.~Berrone, A.~Borio, C.~Fidelibus, S.~Pieraccini, S.~Scial{\`o}, and
  F.~Vicini}, {\em Advanced computation of steady-state fluid flow in discrete
  fracture-matrix models: Fem--bem and vem--vem fracture-block coupling}, GEM -
  International Journal on Geomathematics, 9 (2018), pp.~377--399.

\bibitem{BFPSV18}
{\sc S.~Berrone, C.~Fidelibus, S.~Pieraccini, S.~Scial{\`o}, and F.~Vicini},
  {\em Unsteady advection-diffusion simulations in complex {D}iscrete
  {F}racture {N}etworks with an optimization approach}, Journal of Hydrology,
  566 (2018), pp.~332--345.

\bibitem{BERRONE2017768}
{\sc S.~Berrone, S.~Pieraccini, and S.~Scialò}, {\em Flow simulations in
  porous media with immersed intersecting fractures}, Journal of Computational
  Physics, 345 (2017), pp.~768 -- 791.

\bibitem{BSV}
{\sc S.~Berrone, S.~Scial{\`o}, and F.~Vicini}, {\em Parallel meshing,
  discretization and computation of flow in massive {D}iscrete {F}racture
  {N}etworks}, SIAM J. Sci. Comput., 41 (2019), pp.~C317--C338.

\bibitem{DFNBoon}
{\sc W.~Boon, J.~Nordbotten, and I.~Yotov}, {\em Robust discretization of flow
  in fractured porous media}, SIAM Journal on Numerical Analysis, 56 (2018),
  pp.~2203--2233.

\bibitem{Brenner2016}
{\sc K.~Brenner, M.~Groza, C.~Guichard, G.~Lebeau, and R.~Masson}, {\em
  Gradient discretization of hybrid dimensional {D}arcy flows in fractured
  porous media}, Numerische Mathematik, 134 (2016), pp.~569--609.

\bibitem{VEMmixedBasic}
{\sc {Brezzi, F.}, {Falk, Richard S.}, and {Donatella Marini, L.}}, {\em Basic
  principles of mixed virtual element methods}, ESAIM: M2AN, 48 (2014),
  pp.~1227--1240.

\bibitem{Chave2018}
{\sc F.~Chave, D.~Di~Pietro, and L.~Formaggia}, {\em A hybrid high-order method
  for darcy flows in fractured porous media}, SIAM Journal on Scientific
  Computing, 40 (2018), pp.~A1063--A1094.

\bibitem{DPbook}
{\sc Z.~Chen, G.~Huan, and Y.~Ma}, {\em {Computational Methods for Multiphase
  Flows in Porous Media}}, SIAM, Philadelphia, PA, USA, 2006.

\bibitem{Coulet2019}
{\sc J.~Coulet, I.~Faille, V.~Girault, N.~Guy, and F.~Nataf}, {\em A fully
  coupled scheme using virtual element method and finite volume for
  poroelasticity}, Computational Geosciences,  (2019).

\bibitem{caceresmixed2}
{\sc E.~Cáceres and G.~N. Gatica}, {\em A mixed virtual element method for the
  pseudostress–velocity formulation of the stokes problem}, IMA Journal of
  Numerical Analysis, 37 (2017), pp.~296--331.

\bibitem{VEMbrinkman}
{\sc E.~Cáceres, G.~N. Gatica, and F.~A. Sequeira}, {\em A mixed virtual
  element method for the brinkman problem}, Mathematical Models and Methods in
  Applied Sciences, 27 (2017), pp.~707--743.

\bibitem{caceresmixed1}
\leavevmode\vrule height 2pt depth -1.6pt width 23pt, {\em A mixed virtual
  element method for quasi-newtonian stokes flows}, SIAM Journal on Numerical
  Analysis, 56 (2018), pp.~317--343.

\bibitem{VEMimplementation}
{\sc L.~B. da~Veiga, F.~Brezzi, L.~D. Marini, and A.~Russo}, {\em Virtual
  Element Implementation for General Elliptic Equations}, Springer
  International Publishing, Cham, 2016, pp.~39--71.

\bibitem{MixedDarcyNonMatching}
{\sc {D\'{}Angelo, Carlo} and {Scotti, Anna}}, {\em A mixed finite element
  method for darcy flow in fractured porous media with non-matching grids},
  ESAIM: M2AN, 46 (2012), pp.~465--489.

\bibitem{MixedDassi}
{\sc F.~Dassi and S.~Scacchi}, {\em Parallel solvers for virtual element
  discretizations of elliptic equations in mixed form},
  https://arxiv.org/abs/1901.02330, 00 (2019).

\bibitem{DASSI2019}
{\sc F.~Dassi and G.~Vacca}, {\em Bricks for the mixed high-order virtual
  element method: Projectors and differential operators}, Applied Numerical
  Mathematics,  (2019).

\bibitem{DF99}
{\sc W.~S. Dershowitz and C.~Fidelibus}, {\em Derivation of equivalent pipe
  networks analogues for three-dimensional discrete fracture networks by the
  boundary element method}, Water Resource Res., 35 (1999), pp.~2685--2691.

\bibitem{Faille2016}
{\sc I.~Faille, A.~Fumagalli, J.~Jaffr{\'e}, and J.~E. Roberts}, {\em Model
  reduction and discretization using hybrid finite volumes for flow in porous
  media containing faults}, Computational Geosciences, 20 (2016), pp.~317--339.

\bibitem{X2}
{\sc C.~Fidelibus}, {\em The 2d hydro-mechanically coupled response of a rock
  mass with fractures via a mixed {BEM}-{FEM} technique}, International Journal
  for Numerical and Analytical Methods in Geomechanics, 31 (2007),
  pp.~1329--1348.

\bibitem{formaggia2014}
{\sc L.~Formaggia, A.~Fumagalli, A.~Scotti, and P.~Ruffo}, {\em A reduced model
  for {D}arcy’s problem in networks of fractures}, ESAIM: Mathematical
  Modelling and Numerical Analysis, 48 (2014), pp.~1089--1116.

\bibitem{XFEMreview}
{\sc T.-P. Fries and T.~Belytschko}, {\em The extended/generalized finite
  element method: an overview of the method and its applications}, Internat. J.
  Numer. Methods Engrg., 84 (2010), pp.~253--304.

\bibitem{Fumagalli2018}
{\sc A.~Fumagalli and E.~Keilegavlen}, {\em Dual virtual element method for
  discrete fractures networks}, SIAM Journal on Scientific Computing, 40
  (2018), pp.~B228--B258.

\bibitem{FKS}
{\sc A.~Fumagalli, E.~Keilegavlen, and S.~Scial\`{o}}, {\em Conforming,
  non-conforming and non-matching discretization couplings in discrete fracture
  network simulations}, J. Comput. Phys., 376 (2019), pp.~694--712.

\bibitem{FS13}
{\sc A.~Fumagalli and A.~Scotti}, {\em A numerical method for two-phase flow in
  fractured porous media with non-matching grids}, Advances in Water Resources,
  62 (2013), pp.~454 -- 464.

\bibitem{Fumagalli2019}
{\sc {Fumagalli, Alessio} and {Keilegavlen, Eirik}}, {\em Dual virtual element
  methods for discrete fracture matrix models}, Oil Gas Sci. Technol. - Rev.
  IFP Energies nouvelles, 74 (2019), p.~41.

\bibitem{cubature}
{\sc M.~Gentile, A.~Sommariva, and M.~Vianello}, {\em Polynomial interpolation
  and cubature over polygons}, Journal of Computational and Applied
  Mathematics, 235 (2011), pp.~5232 -- 5239.

\bibitem{Hajibeygi2011}
{\sc H.~Hajibeygi, D.~Karvounis, and P.~Jenny}, {\em A hierarchical fracture
  model for the iterative multiscale finite volume method}, Journal of
  Computational Physics, 230 (2011), pp.~8729 -- 8743.

\bibitem{Gable2015bTransp}
{\sc J.~D. Hyman, S.~Karra, N.~Makedonska, C.~W. Gable, S.~L. Painter, and
  H.~S. Viswanathan}, {\em dfnworks: A discrete fracture network framework for
  modeling subsurface flow and transport}, Computers \& Geosciences, 84 (2015),
  pp.~10 -- 19.

\bibitem{Li2008}
{\sc L.~Li and S.~H. Lee}, {\em Efficient {F}ield-{S}cale {Simulation} of
  {B}lack {O}il in a {N}aturally {F}ractured {R}eservoir {T}hrough {D}iscrete
  {F}racture {N}etworks and {H}omogenized {M}edia}, SPE Reservoir Evaluation \&
  Engineering, 11 (2008).

\bibitem{Lipnikov2014}
{\sc K.~Lipnikov, G.~Manzini, and M.~Shashkov}, {\em Mimetic finite difference
  method}, Journal of Computational Physics, 257 (2014), pp.~1163--1227.

\bibitem{ModelingFractures}
{\sc V.~Martin, J.~Jaffré, and J.~Roberts}, {\em Modeling fractures and
  barriers as interfaces for flow in porous media}, SIAM Journal on Scientific
  Computing, 26 (2005), pp.~1667--1691.

\bibitem{Moinfar2014}
{\sc A.~Moinfar, A.~Varavei, K.~Sepehrnoori, and R.~Johns}, {\em {D}evelopment
  of an {E}fficient {E}mbedded {D}iscrete {F}racture {M}odel for 3{D}
  {C}ompositional {R}eservoir {S}imulation in {F}ractured {R}eservoirs}, SPE
  Journal, 19 (2014).

\bibitem{Neuman05}
{\sc S.~P. Neuman}, {\em Trends, prospects and challenges in quantifying flow
  and transport through fractured rocks}, Hydrogeol. J., 13 (2005),
  pp.~124--147.

\bibitem{NGO2017}
{\sc T.~D. Ngo, A.~Fourno, and B.~Noetinger}, {\em Modeling of transport
  processes through large-scale discrete fracture networks using conforming
  meshes and open-source software}, Journal of Hydrology, 554 (2017), pp.~66 --
  79.

\bibitem{NOETJCP15}
{\sc B.~N{\oe}tinger}, {\em A quasi steady state method for solving transient
  {D}arcy flow in complex 3{D} fractured networks accounting for matrix to
  fracture flow}, J. Comput. Phys., 283 (2015), pp.~205--223.

\bibitem{NOETJCP12}
{\sc B.~N{\oe}tinger and N.~Jarrige}, {\em A quasi steady state method for
  solving transient {D}arcy flow in complex 3{D} fractured networks}, J.
  Comput. Phys., 231 (2012), pp.~23--38.

\bibitem{Nordbotten2019}
{\sc J.~M. Nordbotten, W.~M. Boon, A.~Fumagalli, and E.~Keilegavlen}, {\em
  Unified approach to discretization of flow in fractured porous media},
  Computational Geosciences, 23 (2019), pp.~225--237.

\bibitem{NTTDA92}
{\sc A.~W. Nordqvist, Y.~W. Tsang, C.~F. Tsang, B.~Dverstop, and J.~Andersson},
  {\em A variable aperture fracture network model for flow and transport in
  fractured rocks}, Water Resource Res., 28 (1992), pp.~1703--1713.

\bibitem{P2}
{\sc G.~Pichot, J.~Erhel, and J.~de~Dreuzy}, {\em A generalized mixed hybrid
  mortar method for solving flow in stochastic discrete fracture networks},
  SIAM Journal on scientific computing, 34 (2012), pp.~B86 -- B105.

\bibitem{Qi2005}
{\sc D.~Qi and T.~Hesketh}, {\em An analysis of upscaling techniques for
  reservoir simulation}, Petroleum Science and Technology, 23 (2005),
  pp.~827--842.

\bibitem{Sandve2012}
{\sc T.~Sandve, I.~Berre, and J.~Nordbotten}, {\em An efficient multi-point
  flux approximation method for discrete fracture–matrix simulations},
  Journal of Computational Physics, 231 (2012), pp.~3784 -- 3800.

\bibitem{VMS07}
{\sc M.~Vohral\'{\i}k, J.~Mary\v{s}ka, and O.~Sever\'yn}, {\em Mixed and
  nonconforming finite element methods on a system of polygons}, Applied
  Numerical Mathematics, 51 (2007), pp.~176--193.

\end{thebibliography}

\end{document}